\documentclass[red,11pt,a4paper]{article}

\usepackage{amsmath}
\usepackage{amscd}
\usepackage{amssymb}
\usepackage[pdftex]{color,graphicx,hyperref}
\usepackage{color}
\usepackage{cases}
\usepackage{comment}
\usepackage[normalem]{ulem}
\usepackage{cancel}
\usepackage{cite}
\usepackage{ika}

\usepackage{xcolor}
\usepackage{mathrsfs}
\bf

\newcommand{\dv}{\Div}
\newcommand{\sgn}{{\rm sgn}}
\newcommand{\iintO}[1]{\int_{\mathbb{T}^{2d}} #1 \ \dx\, \dy}
\newcommand{\T}{\mathbb{T}}
\newcommand{\Td}{\mathbb{T}^d}

\newcommand{\iintTO}[1]{\int_0^T\!\!\!\! \int_{\mathbb{T}^{2d}} #1 \ \dx\, \dy\, \dt}

\newcommand{\etat}{{\eta_t}}
\newcommand{\etax}{{\eta_x}}

\def\calD{{\mathcal D}}
\title{A new construction of weak solutions to compressible Navier-Stokes equations}
\author{Nilasis Chaudhuri$^1$, Piotr B. Mucha$^2$, Ewelina Zatorska$^{1}$}
\date{\today}

\begin{document}
\maketitle

\medskip

\begin{center}

1. Department of Mathematics, Imperial College London,\\
South Kensington Campus -- SW7 2AZ, London, UK.\\

2. Institute of Applied Mathematics and Mechanics, University of Warsaw, \\
ul. Banacha 2, 02-097 Warszawa, Poland.

\end{center}

\begin{abstract}
We prove the existence of the weak solutions to the compressible Navier--Stokes system with barotropic pressure $p(\varrho)=\varrho^\gamma$ for $\gamma\geq 9/5$ in three space dimension. The novelty of the paper is the approximation scheme that instead of the classical regularization of the continuity equation (based on the viscosity approximation  $\ep \Delta \varrho$) uses more direct truncation and regularisation of nonlinear terms an the pressure. This scheme is compatible with the Bresch-Jabin compactness criterion for the density. We revisit this criterion and prove, in full rigour, that it can be applied in our approximation at any level.  
\end{abstract}

\maketitle

\section{Introduction}
The current paper deals with the classical compressible barotropic Navier-Stokes system in the periodic domain  $\Omega=\mathbb{T}^d$ with $ d=3$. The unknowns of the system: the mass density $\vr$ and the velocity vector field $\vu$ satisfy the conservation of mass and conservation of momentum  equations, respectively:
\begin{subnumcases}{\label{main}}
	 \pt\vr+\Div (\vr \vu) = 0,\label{main1}\\
	\ptb{\vr\vu}+\Div (\vr \vu \otimes \vu) - \Div \vS(\vu)+ \Grad p(\vr) =\vc{0},\label{main2}
\end{subnumcases}
with the corresponding set of initial data
\begin{align}\label{ini:cond}
\vr(0,x)=\vr_0(x),\quad \vr\vu(0,x)=\vc{m}_0(x).
\end{align} 
In the system above,  $\vS$ denotes the Newtonian stress tensor
\begin{equation}\label{Su}
\vS( \vu) = \mu \left( \Grad \vu + \Grad^t \vu \right) + \xi \Div \vu \, \mathbf{I}, \quad \mu > 0 , \ 2\mu + d\,\xi > 0,
\end{equation}
and $p$ denotes the internal pressure
\begin{equation}\label{prho}
p(\vr)=\vr^\gamma,\quad\quad \gamma\geq \frac32 ~~\mbox{for}\ d=2,\quad \gamma\geq \frac95~~ \mbox{for}\ d=3.
\end{equation}

The purpose of this work is to propose a construction of the weak solutions to the system \eqref{main} that is compatible with the Kolmogorov criterion for compactness of the density. This criterion has been recently resurrected by Bresch and Jabin in \cite{BJ} to prove sequential compactness of the weak solutions to the compressible Navier-Stokes equations with non-monotone pressure and anisotropic stress tensor. The authors refined the classical effective viscous flux argument introduced by Lions \cite{Lions2}, and extended by Feireisl \cite{EF2001}, with direct regularity estimates for the sequence approximating the density. Bresch and Jabin showed that this sequence converges strongly on the support of certain weight satisfying the transport equation. More importantly, this weight vanishes only on vacuum, which yields the strong convergence of the density in the whole domain. This is quite generic technique whose potential has been already demonstrated, for example, to prove the existence of solutions to the multi-phase fluid models \cite{BMZ}, to the model with heterogeneous pressure \cite{BJW} and to the models of cancer growth \cite{VaZa, DPSV}.

Our main theorem is in fact the same as the original result of Lions  \cite[Theorem 7.2]{Lions2}, in the sense that the restriction on the adiabatic exponent $\gamma\geq \frac95$ is required in order to provide the uniform $L^2(0,T; L^2(\T^d))$ bound of the density for $ d=3 $. We state our main theorem for $ d=2,3 $ as follows:
\begin{thm}\label{Main-Thm}
Let $\gamma\geq\frac 32$ if $d=2$ and $\gamma\geq\frac 95$ if $d=3$. Let the initial conditions  \eqref{ini:cond} satisfy 
\eq{\label{ini-cond:hyp}
&\vr_0\in L^\gamma(\T^d),\quad \vr_0\geq 0,\ \mbox{a.e.\ in } \T^d,\ \vr_0\neq 0,\\
&\frac{|\vc{m}_0|^2}{\vr_0}\in L^1(\T^d),\quad \vc{m}_0= 0 \mbox{\ \ a.e.\ on }\{\vr_0=0\}.
}
{Then $(\vr,\vu)$ is the bounded energy weak solution such that $\vr\in L^\infty(0,T; L^\gamma(\T^d))$ and $ \vu \in L^2(0,T;W^{1,2}(\T^d)) $. Moreover, $(\vr,\vu)$ satisfies the renormalized continuity equation and $\vr \in C_{weak}([0,T];L^\gamma(\T^d))$ and $\vr \vu \in C_{weak}([0,T];L^{2\gamma/\gamma+1}(\T^d))$.}
\end{thm}

The proof of the theorem is presented for the more difficult case of $d=3$ only.

Let us conclude this section with a few more remarks concerning the existence of solutions for the compressible Navier-Stokes system.

 First of all, we highlight that the idea of Bresch and Jabin \cite{BJ} was originally proposed to tackle more complex pressure forms.  Indeed, instead of the pressure law \eqref{prho}, one can consider a more general pressure-density relation:
\begin{equation}\label{gen-prho}
    C^{-1}\rho^\gamma - C \leq p(\rho) \leq 
    C \rho^\gamma + C \mbox{ \ \ with \ \ }
    |p'(\rho)|\leq C \rho^{\tilde{\gamma} -1},
\end{equation}
where $C>0$ is a constant, and $ p $ is locally Lipschitz function with $p(0)=0$. Moreover, $\gamma$ and $\tilde{\gamma}$ are related by
\begin{align}\label{gen-gamma}
    \gamma> \lr{\max\{\tilde{\gamma},2\} +1 } \frac{d}{d+2}.
\end{align} 
The extension of our  technique to cover this case can be achieved by a  modification of the main argument of  Bresch and Jabin, see \cite[Section 4.2]{BJ2}.

 Secondly, even if we stick with the monotone pressure law \eqref{prho}, the range of $\gamma$ in Theorem \ref{Main-Thm} can be improved. Indeed, using the  \textit{oscillation defect measure} introduced by Feireisl \cite[Section 6.4]{EF2001}, one can relax the assumptions on $\gamma$ up to $\gamma > 1$ when $d=2$ and $\gamma > \frac 32$ when $d=3$. Unfortunately, we are not able to obtain the same range of $ \gamma$'s by a simple application of the technique based on Bresch and Jabin's idea only.

The paper is organized as follows. In Section \ref{Sec:2} we first introduce our basic level of approximation, for which the existence of solutions is proven using classical fixed point argument in Section \ref{Sec:3}. In Section \ref{Sec:4} we present the details of compactness technique of Bresch and Jabin  \cite{BJ}  for the sequence approximating the density, and we use it to let first approximation parameter to go to zero. The subsequent limit passages are presented in Sections \ref{Sec:5}--\ref{Sec:8}. The proofs of all technical auxiliary lemmas are presented in the Appendix at the end.

\section{Approximation} \label{Sec:2}

The baseline approximation  of  \eqref{main} reads as follows:
\begin{subnumcases}{\label{maine}}
	 \pt\vr+\Div (\vr \vue) +k  \vr^m= 0,\label{maine1}\\[10pt]
	(\delta+[\vr]_\delta)\pt\vu+[\vr \vue]_\delta
	\cdot\Grad\vu - \Div \vS(\vu)+ 
	\Grad p^{M,\lambda}_\ep(\vr) =\vc{0}.\label{maine2}
\end{subnumcases}
We supplement this system with the initial data
\eqh{
\vr(0,x)=\vr_0^\lambda,\quad \vu(0,x)=\vu_0^\lambda
}
such that
\eq{\label{id-approx}
0\leq \vr^\lambda_0 ,\quad
\vr^\lambda_0 \in W^{1,p}(\T^d) \quad \text{and}\quad \vu^\lambda_0 \in W^{2-2/p,p}(\T^d), }
for arbitrary $p>1$, uniformly w.r.t. $\lambda$. In particular $\vr_0^\lambda$ is uniformly bounded in $L^\infty(\T^d)$ independently of all approximation parameters but $\lambda$. Moreover, we assume that the initial data converges strongly to the limiting data in the following sense
\eq{\label{id-approx-conv}
\vr_0^\lambda\to \vr_0\quad \text{strongly in } L^\gamma(\T^d), \quad
\vr_0^\lambda|\vu_0^\lambda|^2\to \frac{|\vc{m}_0|^2}{\vr_0}\quad \text{strongly in } L^1(\T^d).
}


The approximation uses five positive parameters $\ep,\,\delta,\,k,\,M,\,\lambda$ that will vary, and two fixed ones $m,\,\Gamma$ that will be determined later on.
Let us now explain the role of each parameter in order corresponding to subsequent limit passages.
\begin{itemize}
    \item Parameter $\ep$ ($\ep\to 0$) is the regularising parameter in space. We denote $f_\ep=\kappa_\ep\ast f$,
where $\kappa_\ep$ is a mollifying kernel, its choice will be clarified later when the compactness criterion is used. For now we assume that $f_\ep$ is $C^\infty(\T^d)$ function.
\item Parameter $M$ ($M\to\infty$) is the truncation parameter introduced to obtain the a-priori $L^\infty$ bound on the pressure:
\eq{\label{pM}
p^{M,\lambda}(\vr):=\left\{\begin{array}{lll}
\lambda\vr^\Gamma+\vr^\gamma& \mbox{for}& \vr<M\\
\lambda M^\Gamma+M^\gamma& \mbox{for}& \vr\geq M
\end{array}\right.:=\lambda\vr_M^\Gamma+\vr_M^\gamma.
}
\item Parameter $k$ ($k\to0$) provides better a-priori integrability of the density on the level when $\ep=0$ and up to the level of approximation when the so-called Bogovskii estimate can be derived.
\item Parameter $\delta$ ($\delta\to0$) provides another regularisation in space, we denote $[f]_\delta = \kappa_\delta\ast f.$
Note that applying this regularization to the continuity equation \eqref{maine1} we obtain
\eq{\label{maine1O}
 \pt [\vr]_\delta+\Div ( [\vr \vue]_\delta) + k [\vr^m]_\delta= 0.
}

Thanks to this the momentum equation \eqref{maine2} can be rewritten (formally) in the conservative form as
\eq{\label{mom_conserv}
\pt( ([\vr]_\delta +\delta) \vu)+\Div( [\vr \vue]_\delta \otimes\vu) - \Div \vS(\vu)+ \Grad p^{M,\lambda}_\ep(\vr) +k [\vr^m]_\delta \vu
=\vc{0}.
}
This regularization causes that the effective viscous flux {$ \lr{(2\mu+\xi) \Div \vu - p(\vr) }$} is relatively highly regular towards the  finial steps of the approximation.

\item And finally,  parameter $\lambda$ ($\lambda\to 0$) corresponds to artificial pressure $\vr^\Gamma$ that is introduced to improve the integrability of the sequence approximating the density, analogously as in the approximation schemes of Lions \cite{Lions2} and Feireisl \cite{EF2001}.

\end{itemize}

Let us now discuss the key components of the proof and main mathematical difficulties in our approximation scheme. The heuristic arguments presented now are justified rigorously later on in the paper.
At the basic level of approximation, as $\ep >0$ is fixed, the existence or strong solutions follows from the standard existence theory for quasilinear parabolic equations combined with a fixed point argument.

The previous constructions of weak solutions by Lions \cite{Lions2} and by Feireisl \cite{EF2001} use the viscous approximation of continuity equation instead of mollification of transport velocity. The lack of viscous dissipation approximation ($-\ep \Delta \vr$) in our approach makes it fit to apply the Bresch and Jabin criterion for compactness of the density sequence \cite{BJ}, in a direct way. Our approximation is designed in such a way that compactness of density at all following levels of construction is justified using the same argument.

The core of this argument is the analysis of two equivalent formulations of  the effective viscous flux. For the limiting system \eqref{main} they can be written as
\begin{align}
   (2\mu+\xi) \Div \vu - p(\vr)&=(-\Delta)^{-1}\Div( \vr \pt\vu + \vr \vu \cdot\nabla \vu)\label{EVF:form 1}\\
    &=(-\Delta)^{-1}\Div( \pt(\vr \vu) + \Div(\vr \vu \otimes \vu))\label{EVF:form 2},
\end{align}
where we used the equation of continuity to pass from one line to the other.
The first form of the effective viscous flux \eqref{EVF:form 1} is better suited to analyze the convective term $\vr \vu \cdot\nabla \vu$, but as we have no information on the time derivative $\pt\vu$, we take advantage of form \eqref{EVF:form 2}, to work with $ \pt(\vr \vu) $ instead. Its negative Sobolev norm is controlled simply using the  momentum equation.

We implement the above idea in practice by splitting the velocity into two components 
$ \vu= \vu_\eta+\vv_\eta=\pi_\eta\ast\vu+ (\vu-\pi_\eta\ast\vu)$ where $\pi_\eta$ is the time and space mollifier. To show compactness of the effective viscus flux corresponding to $\vu_\eta=\pi_\eta\ast\vu$ we use the  formulation \eqref{EVF:form 1}, i.e.
\begin{equation}\nonumber
    (-\Delta)^{-1}\Div( \vr \pt{\vu_\eta} + \vr \vu \cdot\nabla \vu_\eta).
    \end{equation}
To control the reminder, $\vv_\eta=(\vu -\pi_\eta\ast\vu)$, along with the estimate $	\|\vr|\vv_\eta|\|_{L^1_{t,x}}\leq C \eta^\vt$ from
Lemma \ref{lem:approx}, we use the  fractional integration by parts formula. More precisely, for $\alpha \in (0,1)$, we prove the interpolation inequality (see Lemma \ref{interpolation-lemma}):
\eq{\label{abstract_form}
		&\int_{\mathbb{R}} \int_{\mathbb{R}^d}  \Big((- \Delta)^{-1}\dv \pt \pmb{\varphi} \Big) W(t,x) \dx\; \dt 
		\\
		&\qquad\leq C \Vert \pmb{\varphi}  \Vert_{L^{{q}}}^{\alpha} \Vert \partial_t \pmb{\varphi} \Vert_{L^{{q}^{\prime}}(\R; W^{-1,{q}^\prime}(\R^d) )}^{1-\alpha}\Vert W \Vert_{L^{\bar{q}}}^{1-\alpha} \Vert \partial_t W \Vert_{L^{\bar{q}^{\prime}}(\R; W^{-1,\bar{q}^\prime}(\R^d) )}^{\alpha},
}
where $ \pmb{\varphi} \in L^q(\R^{d+1}) $ with $ \partial_t \pmb{\varphi} \in L^{q^\prime} (\R; W^{-1,q^\prime} (\R^d)) $ and $ W \in L^{\bar{q}}(\R^{d+1} )  $ with $ \partial_t W \in  L^{{\bar{q}}^\prime} (\R; W^{-1,{\bar{q}}^\prime} (\R^d))  $ such that $ 1\leq q, q^\prime, \bar{q}, \bar{q}^\prime \leq \infty $.
Applying this abstract result to $ \pmb{\varphi}\approx \vr\vv_\eta$ , 
and to some generalized, density-dependent weight $W$,  we  show compactness of the effective viscous flux of the reminder using formulation \eqref{EVF:form 2}, i.e. 
$$ (-\Delta)^{-1}\Div\lr{\pt(\vr\vv_\eta) + \Div(\vr \vu \otimes \vv_\eta)}, $$ 
extended suitably from torus to the whole space.
The consecutive levels of our approximation result in less regularity and integrability of the solution and of the weight $W$, making us use the abstract formula \eqref{abstract_form} rather than direct estimate. More intuitive version of this argument is presented at the first level of approximation, for $\ep>0$ fixed, where the approximate solutions are smooth.

\subsection*{Notation and fundamental inequalities}
Throughout the paper we  consider the following notations for the function space:
\begin{itemize}
	\item For $ k_1,k_2 \in \mathbb{Z}$ and $ 1\leq p_1, p_2 \leq \infty $, we denote 
	$$ L^{p_1}(0,T; L^{p_2}(\T^d)) = L^{p_1}_t L^{p_2}_x   \quad\text{ and }\quad  W^{k_1,p_1}(0,T; W^{k_2,p_2}(\T^d)) = W^{k_1,p_1}_t W^{k_2,p_2}_x .$$ 
	Also, we denote $ L^{p_1}(0,T; W^{k_2,p_2}(\T^d))$  as $  L^{p_1}_t W^{k_2,p_2}_x $, and  $ W^{k,2} $ as $ H^{k} $. Moreover, for $ 1\leq p \leq \infty $, we denote $ L^{p}(0,T; L^{p}(\T^d))  = L^p_{t,x}$.
	\item In certain cases, we write  $ u \in  L^{p+0}(\Td) $ to denote that there exists $q >p $ such that $ u \in  L^{q}(\Td)$.
	
	\item By ${\mathcal{M}}$ we denote the maximal operator.  For $f\in W^{1,p}$, we
	 have the following Lagrange inequality:
\eq{\label{Lagrange_ineq}
	    |f(x)-f(y)|\leq C|x-y|\lr{ \mathcal{M}[|\nabla f|](x)
	    +{\mathcal{M}}[|\nabla f|](y)}.
}
	
	\item \textbf{The regularizing kernel:}
	For $h>0$ we define a family of functions $ K_h $ as follows
 \begin{equation}\label{def:Kh}
  {K}_h(z)=
  \left\{ 
  \begin{array}{lr}
   \frac{1}{(|z|+h)^d}& \text{for}\   |z|\leq 1/2, \\[10pt]
   \frac{1}{(1/2 + h)^d} & \text{otherwise}.
  \end{array}
\right.
 \end{equation}
 In particular, they have the following properties:
 \begin{align}
 	&|z||\Grad K_h(z)|\leq CK_h(z), \label{GradK}\\
 	&\|{K}_h\|_{L^1} \sim |\ln h| \text{ as } h\approx 0^+, \label{NormK}
 \end{align}
{where $ C $ is independent of $ h $.}
 We also denote 
  \begin{equation}\label{def_kappa}
\kappa_h=  \Ov K_h=\frac{1}{\|{K}_h\|_{L^1}} {K}_h. 
 \end{equation}
Observe that since $\int_{\T^d} K_h(z) |z| {\rm d}z \leq C$
 uniformly in $h$, therefore
 \begin{equation}\label{small-kapp}
     \int_{\T^d} \kappa_h(z) |z| {\rm d}z \lesssim |\ln h|^{-1},
 \end{equation}
where the symbol ``$\lesssim$" means that the inequality is satisfied up to a constant.

 \item \textbf{Convolutions and commutators:}  First, recall the classical Young inequality 
  \eq{\nonumber
	    \|f \ast g\|_{L^r} \leq \|f\|_{L^p}\|g\|_{L^q} \mbox{\ \ for \ \ } 1+1/r=1/p+1/q.}
For $f \in L^2$, $g\in W^{1,2}$, we will also frequently use
 \begin{equation}\label{comm-est}
     \| [fg]_\delta - [f]_\delta g\|_{L^1} \lesssim 
     |\ln \delta|^{-1} \|f\|_{L^2}\|\nabla g\|_{L^2}.
 \end{equation}
Finally, a combination of the  Lagrange inequality \eqref{Lagrange_ineq} and \eqref{def_kappa}
implies:
 \begin{equation*}
 \int_{\T^d} \kappa_\delta (x-y) (f(y)g(y) -f(y)g(x)) {\rm d}y\lesssim  \int_{\T^d} \kappa_\delta (x-y) |f(y)|
    |x-y| (\mathcal{M}[\nabla g](x) +\mathcal{M}[\nabla g](y)) {\rm d}y.
 \end{equation*}
 
 \medskip
 
To conclude, we state an important lemma which allows us to control convolution of two kernels $ K_{h_1} $ and $ K_{h_2} $ with $ h_1,h_2>0 $ by the sum of them. 
\begin{lemma}\label{Lemma:ConvK}
Given $h_1,h_2>0$, we have
\begin{equation}\label{out2}
 \int_{\Td} \Ov K_{h_1}(x-y) \Ov K_{h_2}(y-z)\, \dy \leq C(\Ov K_{h_1}(x-z)+\Ov K_{h_2}(x-z)),
\end{equation}
which for short will be written as
\begin{equation*}
 \Ov K_{h_1}\ast \Ov K_{h_2} \leq C(\Ov K_{h_1}+\Ov K_{h_2}).
\end{equation*}
\end{lemma}
The proof of lemma is postponed to the Appendix \ref{PL2}.

 \item Throughout the paper we denote by $C$ a generic constant that may change from line to line. It may depend on the initial data, the dimension and the parameters of the approximation.  In certain places we will write $C(\alpha,\beta)$ to highlight the dependence on certain parameters$ (\alpha,\beta) $.
\end{itemize}

\section{Existence of solutions for all parameters fixed}\label{Sec:3}
This section is dedicated to the proof of the following result.

\begin{thm}\label{exi-fix}
 Let $\ep,\,\delta,\, k, M$, and $\lambda$ be fixed. Then there exists a strong solution $(\vr,\vu)$ to the system (\ref{maine})  such that
 \begin{equation*}
  \vu \in W^{1,2}_p((0,T)\times \mathbb{T}^d)= W^{1,p}_t L^{p}_x\cap L^p_tW^{2,p}_x, \qquad \vr \in W^{1,1}_p((0,T)\times\mathbb{T}^d) = W^{1,p}_t L^{p}_x\cap L^p_tW^{1,p}_x 
 \end{equation*}
for any $p<\infty$ and initial data satisfying \eqref{id-approx} .
\end{thm} 
\pf
The proof of existence of solutions  to system \eqref{maine} is obtained via application the Schauder fixed point theorem: 
\textit{``Let ${\cal T}$ be continuous compact operator from $X$ into $X$, where $X$ is abounded convex subset of a Banach space.  Then $T$ has a fixed point."}\\
We  consider the following map:
\begin{equation*}
 {\cal T}(\tilde \vu) = \vu,
\end{equation*}
where $\tilde \vu \in X=L^\infty_t L^2_x\cap L^2_t W^{1,2}_x$ and $\|\tilde \vu\|_{L^\infty_t L^2_x\cap L^2_t W^{1,2}_x} <R$, with $R$ to be be specified later. 
This set is convex in $L^\infty_t L^2_x\cap L^2_t W^{1,2}_x$.

Having $\tilde \vu$ given, we solve the linear (in $\vu$) system with $\tilde{\vu}_\ep = \kappa_\ep \ast {\tilde \vu}$ 
\begin{subnumcases}{\label{maine_fix}}
	 \pt {\vr}+\Div ({\vr \tilde{\vu}_\ep}) +k {\vr^m}= 0,\label{maine_fix1}\\[10pt]
	(\delta + [\vr]_\delta) \pt\vu+ [\vr \tilde{\vu}_\ep]_\delta  \cdot\Grad\vu - \Div \vS(\vu)
	+ \Grad p^{M,\lambda}_\ep (\vr) =\vc{0}.\label{maine_fix2}
\end{subnumcases}
The first equation can be simply solved by the method of characteristics. 
Note that mollification $\kappa_\ep$ makes function ${\tilde\vu}_\ep$ smooth in space. 
Hence, one finds the following estimate (for any $p<\infty$) for the unique solution $\vr=\vr(\tilde\vu_\ep)$ to the continuity equation:
\begin{equation*}
\sup_{t\in (0,T)} ( \|\vr(t)\|_{L^p_x} + \|\nabla \vr(t)\|_{L^p_x}) \leq \left(\|\vr_0^\lambda\|_{L^p_x}+\|\nabla \vr_0^\lambda\|_{L^p_x}+\|\nabla^2 \tilde\vu_\ep\|_{L^1_tL^p_x}\right)
\exp\left\{ \int_0^T \|\nabla {\tilde\vu}_\ep\|_{L^\infty_x} (t)\dt\right\}.
\end{equation*}
The above estimate implies that $\nabla p^{M,\lambda}_\ep (\vr) \in L^p_tL^p_x$ for any $p <\infty$. This information allows us to look at equation \eqref{maine_fix2} as a parabolic equation for $\vu$ where $\vr=\vr(\tilde\vu_\ep)$ is now given. Thanks to H\"older continuity of $[\vr]_\delta$ (at this stage we can use (\ref{maine1O}) in order to provide the H\"older regularity in time), we can apply the  theory from \cite[Section 8.2]{DHP} to prove solvability of  equation \eqref{maine_fix2} with the following estimates
\begin{equation*}
 \|\vu,\pt\vu,\nabla^2 \vu \|_{L^p_{t,x}} \leq C.
\end{equation*}
Thanks to this we have compactness of the operator ${\cal T}$ in $X$.  We now  justify that ${\cal T}$ maps the ball in $L^\infty_t L^2_x\cap L^2_t W^{1,2}_x$ into itself.

Testing \eqref{maine_fix2} by $\vu$, multiplying \eqref{maine_fix1}  regularized in $\delta$, as in \eqref{maine1O} by $\frac{|\vu|^2}{2}$, and summing up the obtained expressions we get the following inequality
\begin{equation*}
 \frac{1}{2} \frac{d}{dt} \intO{ (\delta+[\vr]_\delta) \vert\vu\vert^2} + \frac{k}{2}
 \intO{[\vr^m]_\delta \vert\vu\vert^2} +\intO{\vS(\vu):\nabla \vu}  \leq C(M) \|\Div \vu\|_{L^2_x}.
\end{equation*}
This way we find
\begin{equation*}
 \sup_{t\leq T} \intO{\delta \vu^2} +  \intTO{|\nabla \vu|^2} \leq \intO{\vr_0^\lambda |\vu_0^\lambda|^2} +C(M,T).
\end{equation*}
For fixed $k$ the l.h.s. of above defines the square of the classical norm  $L^\infty_t L^2_x\cap L^2_t W^{1,2}_x$.

Hence we see that the ${\cal T}$ maps a ball into itself, i.e. we have
\begin{equation*}
 \|\vu\|_{L^\infty_t L^2_x\cap L^2_t W^{1,2}_x} \leq R
\end{equation*}
with 
 $R = \sqrt{ \lr{ \intO{\vr_0^\lambda |\vu_0^\lambda|^2} +C(M,T)}}.$
Thus, using Schauder's fixed point theorem we conclude the existence for all given parameters fixed. $\Box$

\section{Passage to the limit $\ep\to 0$}\label{Sec:4} 
In this section, we fix all  parameters other than $ \ep $.  
In order to simplify the notation, we drop the dependence on parameter $\ep$ for the sequence $\{\vr^\ep,\vu^\ep\}_{\ep >0}$ and denote it by $\{\vr,\vu\}$. However, we still use the subscript $\ep$  to denote the regularization in space for $\vu_\ep$ and $p^{M,\lambda}_\ep$.

\medskip

The following strategy of the proof is common for each level of approximation procedure:
\begin{itemize}
    \item Our first goal is to derive uniform bounds w.r.t. $ \ep $, see Section \ref{Sec:4.1}, from which we deduce weak convergence of the approximating sequences.
    \item In order to pass to the limit in the non-linear terms, we need some strong convergence results, in particular, the strong convergence of the density. The required compactness criterion is described in details in Sections \ref{Sec:4.2}--\ref{Ep:comp3}. 
    \item Finally, performing the limit passage $\ep\to 0$, we recover system \eqref{maine} without the parameter $ \ep $. 
\end{itemize}

\subsection{Estimates uniform w.r.t. $\ep$}\label{Sec:4.1}
The uniform estimates are consequences of the basic energy estimate and from an additional estimate of the  `effective viscous flux'.

\subsubsection{Energy estimate}
\begin{lemma}\label{est-ep1}
    Suppose the parameters $\delta,\, k, M \text{ and } \lambda$ are fixed. Then, we have the following uniform bounds w.r.t. $ \ep $:
    \eqh{\vr e^{M,\lambda}(\vr) \in L^\infty_{t}L^1_x, \quad \vu \in L^2_t H^{1}_{x}\cap L^\infty_tL^2_x,\;
    }
    where $(\vr,\vu)$ is a solution obtained in Theorem \ref{exi-fix}, and $\vr e^{M,\lambda}$ is given by
\eq{\label{eM}
\vr e^{M,\lambda}(\vr)=\left\{\begin{array}{lll}
\frac{\lambda}{\gamma-1}\vr^\Gamma+\frac{1}{\gamma-1}\vr^\gamma& \mbox{for}& \vr<M,\\
\frac{\lambda\Gamma}{\Gamma-1}\vr M^{\Gamma-1}- M^\Gamma+\frac{\gamma}{\gamma-1}\vr M^{\gamma-1}- M^\gamma& \mbox{for}& \vr\geq M.
\end{array}\right.
}
\end{lemma}
\pf At first, we note that the equation \eqref{maine1} implies 
\eqh{\Dt \intO{{\vr }} + k\intO{\vr^m} =0.} 
Since the initial data is uniformly bounded w.r.t. $ \ep $, we have $ \vr \in L^\infty_{t}L^1_x$, uniformly in $\ep $. Notice that our approximation  \eqref{maine} keeps the structure, i.e. we can show that the total energy $E$ with 
$$E= \frac{1}{2} \vr \vert \vu \vert^2 + \vr e^{M,\lambda}(\vr)$$
 as a sum of kinetic and internal energies is conserved.

Indeed, taking into account our truncation of pressure $p^{M,\lambda}$ given by \eqref{pM}, we multiply \eqref{maine1O} by $(\vr e_M(\vr))'$ and obtain the followingequation for evolution of  $\vr e^{M,\lambda}$:
\eq{\label{ene}
\pt(\vr e^{M,\lambda}(\vr))+\Div(\vr e^{M,\lambda}(\vr)\vue)+k \vr^m \frac{\lambda\Gamma}{\Gamma-1}\vr_M^{\Gamma-1}+k \vr^m \frac{\gamma}{\gamma-1}\vr_M^{\gamma-1}+p^{M,\lambda}(\vr)\Div\vue=0,}
where used that $(\vr e^{M,\lambda}(\vr))'\vr-\vr e^{M,\lambda}(\vr)=p^{M,\lambda}(\vr)$.

On the other hand, multiplying \eqref{maine2} by $\vu$ and integrating by parts we obtain
\eqh{
&\intO{ (k+ [\vr]_\delta) \pt\lr{\frac{|\vu|^2}{2}}}+
\intO{ [\vr\vue]_\delta \cdot\Grad\lr{\frac{|\vu|^2}{2}}}
\\
&+\intO{\vS(\vu):\Grad\vu}-\intO{p^{M,\lambda}_\ep(\vr)\Div\vu}=0.
}
Now, multiplying the continuity equation \eqref{maine1O} by $\frac{|\vu|^2}{2}$ and adding it to above we get
\eq{\label{enu}
\Dt\intO{ (k+[\vr]_\delta) \frac{|\vu|^2}{2}}
+\intO{\vS(\vu):\Grad\vu}-\intO{p^{M,\lambda}_\ep(\vr)\Div\vu}+\frac{k}{2}
\intO{ [\vr^m]_\delta |\vu|^2}=0.
}
Note that from the properties of convolution it follows that
\eqh{
\intO{p^{M,\lambda}_\ep(\vr)\Div\vu}=\intO{p^{M,\lambda}(\vr)\kappa_\ep\ast(\Div\vu)}=\intO{p^{M,\lambda}(\vr)\Div\vue},
}
therefore, summing up \eqref{enu} with \eqref{ene} integrated over space, we finally obtain
\eq{\label{energy}
&\Dt\intOB{(k+[\vr]_\delta) \frac{|\vu^2|}{2}+\vr e^{M,\lambda}(\vr)}
+\intO{\vS(\vu):\Grad\vu}\\
&\qquad+k \intO{\vr^m \frac{\lambda\Gamma}{\Gamma-1}\vr_M^{\Gamma-1}}+k\intO{ \vr^m \frac{\gamma}{\gamma-1}\vr_M^{\gamma-1}}+\frac{k}{2}
\intO{ [\vr^m]_\delta |\vu|^2}=0.
}
This provides us with uniform bounds w.r.t $\ep$. $\Box$

\subsubsection{Further estimates for $\vr$ and $\Div\vu$}
We start from  improving the estimates for the density.
\begin{lemma}
Let the assumptions of Lemma \ref{est-ep1} be satisfied. Then, uniformly w.r.t. $\ep$ we have:
\eqh{&\vr \in L^\infty_{t}L^{m-1}_x\cap L^{2m-2}_{t,x}.
    }
\end{lemma}
\pf We test (\ref{maine1}) by $\vr^{m-2}$, this leads to:
\begin{equation*}
    \frac{1}{m-1} \frac{d}{dt} \int_{\T^d}
    \vr^{m-1}dx + k\int_{\T^d} \vr^{2m-2}\,\dx
    \leq C\left|\int_{\T^d} \vr^{m-1} \Div \vu\, \dx \right|.
\end{equation*}
Since $\nabla \vu \in L^2_x$, the r.h.s. is bounded. $\Box$

Our next goal is to prove that for fixed $M$, $|p^{M,\lambda}|<C(M)$ implies $\Div\vu$ is in $L^p_t L^\infty_x$. As a consequence of this, we deduce a uniform estimate for density w.r.t. $\ep$ in $L^\infty((0,T)\times\mathbb{T}^d)$. 

\begin{lemma}\label{est-ep2}
Suppose all the hypothesis of Lemma \ref{est-ep1} remains valid. Then we have the following additional estimate for $ \vr$ and $ \Div \vu $: 
\eqh{\Div \vu \in L^p_t L^\infty_x\quad  \text{ and }\quad \vr \in L^\infty_{t,x} \mbox{ \ \ 
for \ \ $ p \in [1,\infty)$. } } Moreover, we obtain $ \vr \vu \in L^2_t L^6_x $.
	\end{lemma}
\pf Instead of working with $\Div\vu$ directly, we focus our attention on the effective viscous flux, denoted by $G$, and given as follows
\eqh{
G= (2\mu+\xi)\Div\vu-p^{M,\lambda}_\ep (\vr).
}
In order to estimate $G$ there is a need of information of the velocity independently from $\ep$. We restate the momentum equation in the following form
\eq{\label{divbs}
 (\delta +[\vr]_\delta) \pt\vu+ [\vr \vue]_\delta  \cdot\Grad\vu - \Div \vS(\vu)
=- \Grad p^{M,\lambda}_\ep(\vr). 
}
First let us note that by definition $p^{M,\lambda}$ is bounded, so the r.h.s. in the above equation belongs to $L^\infty_t W^{-1,\infty}_x $. Thanks to Lemma \ref{est-ep1}
we find that  $[\vr]_\delta $ belongs to  $C^{a}_{t,x}$ for some $a>0$, since $\partial_t 
[\vr]_\delta \in L^{(2m-2)/m}_t C_x$. Since $\delta+[\vr]_\delta \geq \delta >0$ we can show that 
\eqh{\pt \vu \in L^p_t W^{-1,p}_x \mbox{ for arbitrary large}\  p <\infty .} 
The above regularity is a consequence  of the classical maximal regularity results from 
\cite[Section 8]{DHP}. For completeness of the current work, we formulate and prove the relevant result in Lemma \ref{exis-EVF} in Appendix \ref{Lemma16}.

Taking divergence of both sides in \eqref{divbs} we get
\eqh{ (\delta+[\vr]_\delta) \pt\Div\vu - (2\mu+\xi) \lap(\Div\vu)+ \lap p^{M,\lambda}_\ep(\vr) 
=-\Grad [\vr]_\delta \pt \vu-\Div( [\vr\vu_\ep]_\delta \cdot\Grad \vu).}
On the other hand, multiplying the continuity equation by ${p^{M,\lambda}}'(\vr)$, and taking the convolution in space we obtain
\eqh{\pt p^{M,\lambda}_\ep(\vr)+\kappa_\ep \ast \Div(p^{M,\lambda}(\vr)\vu_\ep)+k \kappa_\ep \ast (\vr^m (p^{M,\lambda})'(\vr))+\kappa_\ep \ast \lr{\lr{{(p^{M,\lambda})}'(\vr)\vr-p^{M,\lambda}(\vr)}\Div\vu_\ep }=0,}
and hence we have
\eq{\label{parG}
\frac{\delta+[\vr]_\delta}{2\mu+ \xi} \pt G-\lap G=&-\Grad [\vr]_\delta\cdot \pt \vu-
\Div( [\vr\vu_\ep]_\delta \cdot\Grad \vu)\\
&- \frac{\delta+[\vr]_\delta}{2\mu+ \xi} \left[ \kappa_\ep \ast \Div(p^{M,\lambda}(\vr)\vu_\ep)+k \kappa_\ep \ast \vr^m( {(p^{M,\lambda})}'(\vr))\right. \\
&\hskip4cm \left.+\kappa_\ep \ast \lr{\lr{({p^{M,\lambda}})'(\vr)\vr-p^{M,\lambda}(\vr)}\Div\vu_\ep }\right]
.}
At this point, our goal is to prove that $G \in L^p_t W^{1,p}_x$, hence we need to show that the r.h.s. of \eqref{parG} is in $L^p_tW^{-1,p}_x$ with the norm bounded uniformly w.r.t. $\ep$.

Note that 
$\Grad[\vr]_\delta\cdot \pt \vu$ is in the right space since $[\vr]_\delta$ is smooth in space and bounded in time and $\pt \vu \in L^p_t W^{-1,p}_x$. The same holds true for $\Div([\vr\vu_\ep]_\delta\cdot\Grad \vu)$.
Moreover, we have the uniform in $\ep$ bounds of all terms with $p^{M,\lambda}(\vr)$ and $\vr^m({p^{M,\lambda}})'(\vr)$ as they are equal to $0$ for $\vr>M$. Since $[\vr]_\delta$ is Lipschitz  continuous 
in time and space the other terms on the r.h.s. of (\ref{parG}) are also in $L^p_t W^{-1,p}_x$.
Thus the standard maximal regularity estimate and the fact that $ p $ is arbitrary large imply 
\eqh{G \in L^p_t W^{1,p}_x \subset L^p_t L^\infty_x.}
This in turn yields, using the uniform bound on $p^{M,\lambda}_\ep(\vr)$ w.r.t  $\ep$, that 
$\Div \vu= G +p^{M,\lambda}_\ep(\vr)\in L^p_t L^\infty_x $, which proves the first part of the lemma.

Returning to the estimate for the density, because we have
\begin{equation*}
 \|\Div \vu_\ep\|_{L^p_t L^\infty_x}  \leq  \|\Div \vu\|_{L^p_t L^\infty_x }\leq C(M,k,\delta,\lambda),
\end{equation*}
therefore
\eq{\nonumber
 \inf_{x\in\mathbb{T}^d} {\vr_0^\lambda}\; \exp\left\{ - \int_0^t\|\Div \vu_\ep\|_{L^\infty_x} \, {\rm d}t'\right\} \leq \vr(t,x) \leq \sup_{x\in\mathbb{T}^d} \vr_0^\lambda \; \exp\left\{\int_0^t\|\Div \vu_\ep\|_{L^\infty_x}\,  {\rm d}t'\right\} .} 
So, the density  is bounded uniformly w.r.t. $\ep$. $\Box$

\subsection{Compactness criterion -- formulation}\label{Sec:4.2}

We first recall the classical Kolmogorov compactness criterion, for the proof see \cite[Lemma 3.1]{Belgacem}.
\begin{lemma}\label{lem:compact}
	Let $\{X_n\}_{n=1}^\infty$ be a sequence of functions uniformly bounded in $L^p((0,T)\times\T^d)$
	with $1\leq p<+\infty$. Let $ K_h$ be given by \eqref{def:Kh}.
	If $\{\pt X_n\}_{n=1}^\infty$ is uniformly bounded in $L^r(0,T;W^{-1,r}(\T^d))$ with 
	$r\geq 1$ and  
	\eq{\label{criterion}
		\underset{n}{\lim\sup} \left( \frac{1}{\| K_h\|_{L^1}}  \int_0^T\!\!\!\iintO{
			K_h(x-y)|X_n(t,x)-X_n(t,y)|^p} \, \dt\right) \to 0, \quad \mbox{ as } h\to 0,
	}
	then, $\{X_n\}_{n=1}^\infty$ is compact in $L^p((0,T)\times\T^d)$.
	Conversely, if $\{X_n\}_{n=1}^\infty$ is compact in $L^p([0,T]\times\T^d)$, 
	then \eqref{criterion} holds.
\end{lemma}
\begin{rmk}
The assumptions on the kernel $K_h$ can be relaxed.
It is enough to take a sequence of positive, bounded functions s.t.
	\begin{enumerate}
		\item[i)] $\forall \eta>0$, \ $ \underset{h}{\sup}
	 \intO{K_h(x)\vc{1}_{{\{x:\> |x|\geq\eta\}}}}<\infty$,
		\item[ii)] $\|K_h\|_{L^1(\T^d)}\to+\infty$\ as\ $h\to 0$.
	\end{enumerate}
\end{rmk}

We are not able to apply Lemma \ref{lem:compact} directly into sequence approximating the density. Instead we will first prove the {\emph{weighted compactness}} of $\{\vr^\ep\}_{\ep>0}$. This is the central part of the result by Bresch and Jabin \cite{BJ}. It basically says that the sequence approximating the density $\{\vr^\ep\}_{\ep>0}$ satisfies the following modification of \eqref{criterion}
\eq{\label{criterion0}
\underset{\ep\to0}{\lim\sup} \left(   \frac{1}{\| K_h\|_{L^1}}  \int_0^T\!\!\!\iintO{
			K_h(x-y)|\vr^\ep(t,x)-\vr^\ep(t,y)|(w^\ep(t,x)+w^\ep(t,y))} \, \dt\right) \to 0, \ \mbox{ as } h\to 0
}
involving the weights $w=w^\ep$ satisfying dual to continuity equation, discussed below. Note that in \eqref{criterion0}, we use the notation with $ \ep $ denoting the element of a sequence, playing analogue role to $n$ from the criterion \eqref{criterion}. We will drop this notation when no confusion can arise, and only use $\ep$ to denote the regularisation in space.
\begin{prop}[Proposition 7.2 in \cite{BJ}]
\label{prop:w}
Suppose that ${\cal M}[|\Grad\vu|]$ is  bounded in $L^2_{t,x}$.
Then, for $ \Lambda >0 $, there exists a weight $w$ solving  
\eq{\label{def_w}
\left\{\begin{array}{l}
\pt w+\vue\cdot\Grad w+\Lambda \mathcal{M}[|\Grad\vu|] w=0,\\
w(0,x)=1.
\end{array}\right.
}
Moreover, we have
\begin{itemize}
\item[(i)] For any $(t,x)\in(0,T)\times \T^d$, $0\leq w(t,x) \leq 1$.
\item[(ii)] If we assume moreover that the pair $(\vr,\vue)$ is a solution to the continuity equation:
\eqh{\pt \vr+\Div(\vr\vue)+k\vr^m=0,}
such that $\vr$ is bounded in $L^2_{t,x}$, then there exists $C\geq 0$, such that
\begin{equation}\label{boundlogw}
\underset{t>0}{{\rm ess}\sup} \intO{ \vr |\log w| } \leq C  \Lambda.
\end{equation}
\end{itemize}
\end{prop}

\pf
The only part of this proposition that differs slightly from the proof presented in \cite{BJ} is the proof of estimate \eqref{boundlogw}. Note however, that the equation for $\vr|\log w|$ on the formal level is equal to
\begin{align}\label{eq:rhologw}
    \pt\lr{\vr|\log w|}+\Div(\vr|\log w|\vue)+k\vr^m|\log w|= \vr \Lambda  \mathcal{M}[|\nabla \vu| ],
\end{align}
and the last term on the l.h.s. is nonnegative, so the estimate \eqref{boundlogw} still holds.  $\Box$
\begin{rmk}
    In the next approximation levels, we define the weight $w$ by replacing $ \vue $ with $ \vu $.
\end{rmk}

We stress that the bound (\ref{boundlogw}) is the key point in the restriction of $ \gamma $, which is $\gamma \geq 9/5$ for $ d=3 $. Under this restriction, we  can eventually show that $\vr \in L^2_{t,x}$ uniformly w.r.t. all approximation parameters.

\begin{prop}\label{pro:rho-bdd} Let $\vr_0 \in  L^\infty(\T^d)$, then
\begin{equation}\nonumber
\underset{t>0}{{\rm ess}\sup} \|\vr w\|_{L^\infty} \leq \| \vr_0\|_{L^\infty}.
\end{equation}
\end{prop}
\pf
The product of $\vr w$ satisfies the following equation
\begin{equation*}
 (\vr w)_t + \Div(\vr w \vue) + k\vr^m w + \Lambda  \mathcal{M}[|\nabla \vu|] w \vr=0.
\end{equation*}
And so, we can verify that the transport derivative is bounded from above
\begin{equation*}
 (\vr w)_t + \vue \cdot \nabla (\vr w) = -k \vr^m w - (\Lambda  \mathcal{M}[|\nabla \vu|] +\Div \vue)\vr w  \leq 0,
\end{equation*}
since $\Lambda$ can be taken large enough to have $\Lambda  \mathcal{M}[|\nabla \vu|] \geq |\Div \vue|$. The proof is concluded by integration w.r.t. time. $\Box$

\medskip
\smallskip
\begin{rmk}
We want to prove Theorem \ref{Main-Thm}  for $0\leq \vr_0 \in L^\gamma(\T^d)$. However, for the purpose of approximation we assume $0\leq \vr_0 \in L^\infty(\T^d)$.  Then, in Section \ref{gen-id}, we remove this additional assumption and obtain the result for $0\leq \vr_0 \in L^\gamma(\T^d)$.
\end{rmk}
The purpose of the rest of this section is to prove \eqref{criterion0}. Due to symmetry, it is enough to show the following lemma.
\begin{lemma}\label{Lemma:Kol_Mod}
For the weight $w^\ep $ defined in \eqref{def_w} and $\Ov K_h(x-y)$ defined above, we have
\eq{\label{criterion00}
\underset{\ep\to0}{\lim\sup} \left(  \int_0^T\!\!\!\iintO{
 \Ov K_h(x-y)|\vr^\ep(t,x)-\vr^\ep(t,y)|w^\ep(t,x)} \, \dt\right) \to 0, \ \mbox{ as } h\to 0.
}
\end{lemma}
The proof is presented in two subsections below. 
\subsection{Proof of Lemma \ref{Lemma:Kol_Mod} -- the main estimate}
In the following two subsections we drop the upper index $\ep$ when no confusion can arise. We only keep the lower index $\ep$ which denotes convolution with kernel $\kappa_\ep$.\\

Taking the difference of the equations \eqref{maine1} satisfied by $\vr(x)$
and $\vr(y)$ (we drop the dependence on time), we get
\eq{\label{difference}
&\pt (\vr(x)-\vr(y)) + \dv_x (\vue(x)\lr{\vr(x)-\vr(y)}) +
\dv_y (\vue(y)\lr{\vr(x)-\vr(y)})  \\
&={\dv_x \vue(x) (\vr(x)-\vr(y))
 - \left(\dv_x \vue(x)-\dv_y \vue(y)\right) \vr(x)}\\
&\quad -k \lr{\vr^m(x) -\vr^m(y)} .
}
We want to estimate the modulus of the difference between $\vr(x)$ and $\vr(y)$, because it requires the least integrability of $\vr$. However, at some point of the proof below, we will need to integrate by parts, which requires differentiability of the test function. To enable this,  we modify the standard modulus and signum functions near zero. 
Given $\sigma>0$ we define
\eq{\nonumber
|w|^\sigma= \left\{ 
\begin{array}{lcr}
|w|-\frac{\sigma}{2} & \mbox{ for } & |w| > \sigma, \\
\frac{1}{2\sigma} w^2 & \mbox{ for } & |w|\leq \sigma,
\end{array}
\right.
}
along with
\begin{equation}\label{sigma-bdd}
   \sgn^\sigma:= \partial |w|^\sigma = \left\{ 
    \begin{array}{lcr}
    {\rm sgn \, }w  & \mbox{for} & |w|>\sigma \\
    \frac{w}{\sigma} & \mbox{for} & |w|\leq \sigma
    \end{array}
    \right. 
\mbox{ \ and \ } 
\partial^2 |w|^\sigma = \left\{ 
    \begin{array}{lcr}
    0 & \mbox{for} & |w|>\sigma \\
    \frac{1}{\sigma} & \mbox{for} & |w|\leq \sigma
    \end{array}
    \right. .
\end{equation}
Multiplying \eqref{difference} by $\sgn^\sigma_{xy}:=\sgn^\sigma(\vr(x)-\vr(y))$, we deduce
\eq{\label{eq:drho1}
&\pt |\vr(x)-\vr(y)|^\sigma + \dv_x (\vue(x)|\vr(x)-\vr(y)|^\sigma) +
\dv_y (\vue(y)|\vr(x)-\vr(y)|^\sigma) \\
&+ [(\vr(x)-\vr(y))\sgn^\sigma_{xy} -|\vr(x)-\vr(y)|^\sigma](\dv_x\vue(x)+
\dv_y\vue(y))\\
&\quad  +k  (\vr^m(x) -\vr^m(y))\sgn^\sigma_{xy}  \\
&=\dv_x \vue(x) (\vr(x)-\vr(y))\sgn^\sigma_{xy} 
 - \left[(\dv_x \vue(x)-\dv_y \vue(y))\sgn^\sigma_{xy}\right] \vr(x).
}
\begin{rmk}
To justify this step rigorously, we have to mollify the continuity equations first.
Indeed \eqref{eq:drho1} can be seen as the renormalized equation for $\vr(x)-\vr(y)$. We can justify it, although not always in a pointwise sense, using  the additional regularization of $ |\cdot|^{\sigma} $, the commutator estimates (see for instance \cite[Chap. 3.1]{NS}), and uniform bounds tailored to each level of approximation. 
A similar argument was used for example in \cite[Sec. 5.2]{BMZ}.
\end{rmk}
We introduce the following quantities:
\begin{align}
    &R_h^{\sigma}(t) := \iintO{K_h(x-y) |\vr(x)-\vr(y)|^\sigma w(x)},\label{R}\\
&G_h^{\sigma}(t):= k \iintO{K_h(x-y) (\vr^m(x) -\vr^m(y))\sgn^\sigma_{xy} w(x)}.\label{G}
\end{align}
Multiplying \eqref{eq:drho1} by $K_h(x-y)$ and by $w(x)$, and using the symmetry of $K_h$ we can show that
\eq{\label{eqdtR}
\frac{d}{dt} R_h^{\sigma}(t)+  G_h^{\sigma}(t) = A_1 + A_2 + A_3+A_4+ A_5,
}
where 
\eq{\nonumber
    &A_1 = \iintO{\nabla K_h(x-y) (\vue(x)-\vue(y)) |\vr(x)-\vr(y)|^\sigma w(x)},\\
    &A_2 =  \iintO{ K_h(x-y)|\vr(x)-\vr(y)|^\sigma (\pt w(x) + \vue(x)\cdot\nabla w(x)  +w(x) \dv_x\vue(x))},\\
    &{A_3 = -\iintO{ K_h(x-y)[(\dv \vue(x) - \dv \vue(y))\sgn^\sigma_{xy} ] \vr(x) \,w(x) },}\\
    &A_4=-\iintO{ K_h(x-y)[(\vr(x)-\vr(y))\sgn^\sigma_{xy} -|\vr(x)-\vr(y)|^\sigma](\dv_x\vue(x)+
\dv_y\vue(y))w(x)},\\
&A_5=-\iintO{ K_h(x-y)[(\vr(x)-\vr(y))\sgn^\sigma_{xy} -|\vr(x)-\vr(y)|^\sigma] \dv_x\vue(x) w(x)}.
}
From \eqref{sigma-bdd}, it follows that
\begin{equation}\label{eq:37}
  |(\vr(x)-\vr(y))\sgn^\sigma_{xy} -|\vr(x)-\vr(y)|^\sigma |  \leq \sigma.
\end{equation}
The rest of the proof consists of uniform estimates of $ A_i$'s from the r.h.s. of \eqref{eqdtR}. 

\subsubsection*{Estimate of $A_1$ and $A_2$}
First, recall that due to \eqref{GradK} and symmetry argument we have
\eq{\nonumber
|A_1|&\leq \iintO{|\nabla K_h(x-y)|| \vue(x)-\vue(y)| |\vr(x)-\vr(y)|^\sigma w(x)}\\
&\leq C \iintO{K_h(x-y) \frac{|\vue(x)-\vue(y)|}{|x-y|} |\vr(x)-\vr(y)|^\sigma w(x)}.
}
To continue, we recall the following result.
\begin{lemma}\label{Lagrange}
Let $f \in W^{1,p}(\T^d)$ with $p>1$, then 
\begin{equation}\label{formula:D}
    |f(x)-f(y)|\leq C|x-y|\left( D_{|x-y|} f(x) + D_{|x-y|}f(y)\right) \mbox{ \ a.e. } {x,y \in \T^d,}
\end{equation}
where for $r>0$ we denoted
\begin{equation}\nonumber
D_rf(x)=\frac{1}{r}\int_{B(0,r)} \frac{|\nabla f(x+z)|}{|z|^{d-1} } {\rm d}z.
\end{equation}
\end{lemma}
The proof of this lemma is postponed to the Appendix \ref{Lemma11}. 

Applying this result, we can further estimate $A_1$ as follows
\eq{\label{A1:est}
|A_1|\leq& C\iintO{K_h(x-y) \lr{D_{|x-y|}[\vue(x)]+D_{|x-y|}\vue(y)} |\vr(x)-\vr(y)|^\sigma w(x)}\\
\leq& C\iintO{K_h(x-y) \lr{D_{|x-y|}\vue(x)-D_{|x-y|}\vue(y)} |\vr(x)-\vr(y)|^\sigma w(x)}\\
&+C\iintO{K_h(x-y) {\cal M}[|\Grad\vue(x)|] |\vr(x)-\vr(y)|^\sigma w(x)}\\
\leq& C\iintO{K_h(x-y) \lr{D_{|x-y|}\vue(x)-D_{|x-y|}\vue(y)} |\vr(x)-\vr(y)|^\sigma w(x)}\\
&+C\iintO{K_h(x-y) {\cal M}[|\Grad\vu(x)|] |\vr(x)-\vr(y)|^\sigma w(x)}.
}
 Recall that the transport equation for the weight is given by \eqref{def_w}, and so the term $A_2$ is equal to
\eqh{A_2= \iintO{ K_h(x-y)|\vr(x)-\vr(y)|^\sigma (w(x) \dv_x\vue(x)-\Lambda {\cal M}|\Grad\vue| w)}.}
Since ${\cal M}[|\Grad\vu|]$ is bounded in $L^2_{t,x}$, and  $\Lambda$ can be taken arbitrarily large, the last term in \eqref{A1:est} can be absorbed by  $A_2$, and we conclude that
\eq{\label{a1a2}
|A_1|+ A_2\leq C\iintO{K_h(x-y) \lr{D_{|x-y|}\vue(x)-D_{|x-y|}\vue(y)} |\vr(x)-\vr(y)|^\sigma w(x)}.
}

\subsubsection*{Estimate of $A_3$}
Applying the operator $(-\Delta)^{-1} \dv$ to both sides of {\eqref{maine2}}, we deduce that
\begin{equation}\nonumber
(2\mu + \xi) \dv \vu = p^{M,\lambda}_\ep -(- \Delta)^{-1}\dv( (\delta+[\vr]_{\delta}) \pt \vu + [\vr\vue]_\delta \cdot\nabla \vu) = p^{M,\lambda}_\ep + \calD
\end{equation}
is satisfied in the sense of distributions, where we denoted
\eq{\nonumber
\calD &:=-(- \Delta)^{-1}\dv(( \delta+[\vr]_\delta) \pt \vu + [\vr\vue]_\delta \cdot\nabla \vu)\\
&=-(- \Delta)^{-1}\dv( \pt( (\delta+ [\vr]_\delta) \vu)+\Div( [\vr \vue]_\delta \otimes\vu)  +k [\vr^m]_\delta \vu).}
From now on, we also denote $\kappa_\ep^{(2)}:=\kappa_\ep\ast \kappa_\ep$.  The term $A_3$ takes the form:

\eq{\nonumber
A_{3}
=& - \frac{1}{2\mu+\xi} \iintO{ K_{h}(x-y)  \kappa_\ep\ast \kappa_\ep\ast\lr{p^{M,\lambda}(x)-p^{M,\lambda}(y)}\sgn^\sigma_{xy} \vr(x) w(x)}\\
& - \frac{1}{2\mu+\xi} \iintO{ K_{h}(x-y) \kappa_\ep\ast \Big(\calD(x)-\calD(y)\Big) 
\sgn^\sigma_{xy} \vr(x)   w(x)}\\
=& - \frac{1}{2\mu+\xi} \iintO{ K_{h}(x-y)  \lr{\kappa_\ep^{(2)}\ast p^{M,\lambda}(x)-p^{M,\lambda}(x)}\sgn^\sigma_{xy} \vr(x)  w(x)}\\
& +\frac{1}{2\mu+\xi} \iintO{ K_{h}(x-y)  \lr{\kappa_\ep^{(2)}\ast p^{M,\lambda}(y)-p^{M,\lambda}(x)}\sgn^\sigma_{xy} \vr(x) w(x)}\\
& - \frac{1}{2\mu+\xi} \iintO{ K_{h}(x-y) \kappa_\ep\ast \Big(\calD(x)-\calD(y)\Big)  
\sgn^\sigma_{xy} \vr(x) w(x)}
\\
=&A_{3}^{(i)}+A_{3}^{(ii)}+A_{3}^{(iii)}.
}
Let us now estimate each of the terms on the r.h.s..

\medskip

\noindent{\emph{Estimate of $A_3^{(i)}$.}~} Using the estimate $|p^{M,\lambda}(x)-p^{M,\lambda}(y)|\leq C |\vr(x)-\vr(y)|$, we evaluate
\eq{\label{a31}
(2\mu+\xi) A_3^{(i)}=&- \iintO{ K_{h}(x-y)  \lr{ \kappa^{(2)}_\ep\ast p^{M,\lambda}(x)-p^{M,\lambda}(x)}\sgn^\sigma_{xy} \vr(x) w(x)}
\\
&\leq  C\|K_h\|_{L^1} \intO{   | \kappa^{(2)}_\ep\ast p^{M,\lambda}(x)-p^{M,\lambda}(x)|  w(x)}\\
&\leq C\|K_h\|_{L^1}\intO{  \lr{\int_{\mathbb{T}^{2d}} \kappa_\ep(x-r)\kappa_\ep(r-z)|p^{M,\lambda}(z)-p^{M,\lambda}(x)| \, {\rm d}r\, {\rm d}z}   w(x)}\\
&\leq  C\|K_h\|_{L^1}   \iintO{  {\Ov K_\ep(x-y)}|\vr(y)-\vr(x)|   w(x)}\\
&\leq C\|K_h\|_{L^1}   \iintO{  {\Ov K_\ep(x-y)}|\vr(y)-\vr(x)|^\sigma   w(x)} + C\sigma \|K_h\|_{L^1},
}
where to get the third inequality we used \eqref{def_kappa}  along with \eqref{out2}, and to get the last inequality we used \eqref{eq:37}.

\medskip
\noindent{\emph{Estimate of $A_3^{(ii)}$.}~}We have a sequence of inequalities
\eq{\label{a32}
(2\mu+\xi) A_3^{(ii)}&= \iintO{ K_{h}(x-y)  \lr{ \kappa^{(2)}_\ep\ast p^{M,\lambda}(y)-p^{M,\lambda}(x)}\sgn^\sigma \vr(x) w(x)}
\\
&\leq C \iintO{ K_{h}(x-y) \lr{ \int_{\mathbb{T}^{2d}}\kappa_\ep(y-r)\kappa_\ep(r-z) 
|p^{M,\lambda}(z)- p^{M,\lambda}(x)| \, {\rm d} r\, {\rm d} z  } w(x)}\\
&\leq C\|K_h\|_{L^1}\int_{\mathbb{T}^{3d}}{ \Ov K_{h}(x-y)  \Ov K_\ep(y-z)|p^{M,\lambda}(z) -p^{M,\lambda}(x)|   w(x)}\,\dx\, \dy\, {\rm d}z \\
&\leq C  \|K_h\|_{L^1} \int_{\mathbb{T}^{2d}}{ \lr{ \Ov K_{h}(x-z)+\Ov K_{\ep}(x-z)} |p^{M,\lambda}(z) -p^{M,\lambda}(x)|  w(x)}\,\dx\, {\rm d}z\\
&\leq C \|K_h\|_{L^1} \left(\int_{\mathbb{T}^{2d}}{ \Ov K_{h}(x-y) |\vr(y) -\vr(x)|^\sigma  w(x)}\,\dx\, \dy \right. \\
& \left.\hspace{5cm}+ \int_{\mathbb{T}^{2d}}{ \Ov K_{\ep}(x-y) |\vr(y) -\vr(x)|^\sigma   w(x)}\,\dx\, \dy +\sigma \right).
}

\noindent{\emph{Estimate of $A_3^{(iii)}$.}~} 
At this level of approximation, the regularity of  solutions is so high that we can integrate by parts in time to solve the problem described at the end of Section \ref{Sec:2}. At this step we can avoid application of the  formula \eqref{abstract_form}, but we still need to go through additional regularisation of the velocity vector field. We have the following result.
\begin{lemma}\label{lem:approx}
	Let $\pi_\eta$ denote the time and space mollifier. 
	Let $\vr \in L^\infty_tL^\gamma_x$, $\vu\in L^2_t W^{1,2}_x$, 
	$\partial_t (\vr\vu) \in L^{q}_tW^{-1,q}_x$ with $ q>1 $.
	Then there exists a positive constant $\vt$ depending on $ \gamma $ and $ q $, such that
	\eq{\nonumber
		\|\vr|\vu-\pi_\eta\ast\vu|\|_{L^1_{t,x}}\leq C \eta^\vt.
	}
\end{lemma}
The proof of this Lemma is presented in the Appendix \ref{Lemma12}
(it can also be found in the proof of Lemma 8.3 in \cite{BJ}).

\medskip

Let us first rewrite the time integral of $A_{3}^{(iii)}$ in the following way:
\eq{\label{intD-D}
\int_0^T A_3^{(iii)} {\rm d}t &=- \frac{1}{2\mu+\xi} \iintTO{ K_{h}(x-y) \kappa_\ep\ast \Big(\calD(x)-\calD(y)\Big)  
\sgn^\sigma_{xy} \vr(x) w(x)}  \\
&=- \frac{1}{2\mu+\xi}
\iintTO{ K_{h}(x-y) \kappa_\ep\ast \Big(\calD(\vu(x))-\calD(\vu(y))\Big) \sgn^\sigma_{xy}
\vr(x)  w(x)},
}
where, with the slight abuse of notation, we used the same symbol $\cal D$ to denote
\eq{\label{def:D}
\calD(\bar\vu) &=(- \Delta)^{-1}\dv(\pt( ([\vr]_\delta +\delta) \bar\vu)+\Div( [\vr \vue]_\delta \otimes\bar\vu)  +k [\vr^m]_\delta \vu)\\
&=
(- \Delta)^{-1}\dv( \pt( ([\vr]_\delta +\delta) \bar\vu))+(- \Delta)^{-1}\dv(  \Div( [\vr \vue]_\delta \otimes\bar\vu) )
+k(- \Delta)^{-1}\dv( [\vr^m]_\delta \bar\vu)\\
&:=\calD_1(\bar\vu)+\calD_2(\bar\vu) +\calD_3(\bar\vu).}
For abbreviation, we introduce the notation:
\eq{
 W=W(t,x,y)= \sgn^\sigma_{xy}  \vr(x)   w(x).}
Note that  $W\in L^\infty(0,T;L^\infty(\T^{2d}))$ uniformly w.r.t. $\ep$, but not w.r.t. $M$.
Note also, that if we replace in \eqref{def:D} $\bar\vu$ by $\pi_\eta\ast \vu$ then the integral \eqref{intD-D} is controlled. Indeed, using the continuity equation and then the Lagrange inequality \eqref{Lagrange_ineq} we may write
\eq{\label{a33a}
&\iintTO{ K_{h}(x-y)\kappa_\ep\ast  \Big(\calD(\pi_{\eta}\ast\vu(x))-\calD(\pi_{\eta}\ast\vu(y))\Big) W}\\
&=\iintTO{ K_{h}(x-y)\kappa_\ep\ast  \Big((- \Delta_x)^{-1}\dv_x( (\delta+[\vr]_\delta)\pt\pi_{\eta}\ast\vu)\\
	& \hspace{5cm}-(- \Delta_y)^{-1}\dv_y( (\delta+ [\vr]_\delta)\pt\pi_{\eta}\ast\vu)\Big) W}\\
&\qquad+\iintTO{ K_{h}(x-y) \kappa_\ep\ast \Big((- \Delta_x)^{-1}\dv_x( [\vr \vue]_\delta\cdot\Grad_x\pi_{\eta}\ast\vu)\\
	& \hspace{5cm}-(- \Delta_y)^{-1}\dv_y( [\vr \vue]_\delta\cdot\Grad_y\pi_{\eta}\ast\vu)\Big) W}\\
&\leq C\|W\|_{L^\infty_{t,x,y}} h^\beta\eta^{-\alpha},
}
for some $\alpha,\beta>0$. Thus, the problem is to control the integral \eqref{intD-D} for the residue, i.e. when $\bar\vu$ in \eqref{def:D} is replaced by
\begin{equation}\nonumber
\vve=\vu-\pi_{\eta}\ast \vu.
\end{equation}
Note that $\vve$ is uniformly bounded in $L^2_tW^{1,2}_x$, and thanks to Lemma \ref{lem:approx} we actually know that 
\eq{\label{eta_teta}
\|(\delta + [\vr]_\delta ) \vve \|_{L^1_tL^1_x}\leq C \eta^\vt
}
for some $ \vt>0 $.
Let us now control the component of \eqref{intD-D} corresponding to $\calD_1(\vv_\eta)$, we have
\begin{align}\label{Ep:D1}
\begin{split}
	&\iintTO{ K_{h}(x-y) \kappa_\ep\ast  \Big(\calD_1(\vve(x))-\calD_1(\vve(y))\Big) W}\\
	&=
	\iintTO{ K_{h}(x-y)\kappa_\ep\ast  \Big((- \Delta_x)^{-1}\Div_x( \pt (([\vr]_\delta +\delta)\vve))\\
	& \hspace{5cm}-(- \Delta_y)^{-1}\Div_y( \pt (([\vr]_\delta +\delta)\vve))\Big) W}\\
	&=
	\left. \int_{\T^{2d}} \!\!K_{h}(x-y)\left[ \kappa_\ep\ast \Big((- \Delta_x)^{-1}\Div_x( ([\vr]_\delta +\delta)\vve)-(- \Delta_y)^{-1}\Div_y(([\vr]_\delta +\delta)\vve)\Big)\right] W\, \dx\,\dy \right|^{t=T}_{t=0}\\
	&-\iintTO{ \!\!K_{h}(x-y) \kappa_\ep\ast \Big((- \Delta_x)^{-1}\Div_x( ([\vr]_\delta +\delta)\vve)-(- \Delta_y)^{-1}\Div_y( ([\vr]_\delta +\delta)\vve)\Big) \pt W}.
\end{split}
\end{align}
Let us explain why the integrals at the end points $t=0,T$ vanish. First note that  $W$ is bounded in $L^\infty_{t,x,y}$. 
On the other hand, the operator $-(- \Delta)^{-1}\dv$ ``gains" one derivative, so we can estimate
\eqh{&|\kappa_\ep \ast (- \Delta_x)^{-1}\Div_x( ([\vr]_\delta +\delta)\vve)-\kappa_\ep \ast(- \Delta_y)^{-1}\Div_y(([\vr]_\delta +\delta)\vve)| \\
&\leq C|x-y|\lr{\mathcal{M} |\kappa_\ep \ast (-\Delta)^{-1} \nabla \Div \lr{([\vr]_\delta +\delta)\vve(x)}|+\mathcal{M}|
\kappa_\ep \ast (-\Delta)^{-1} \nabla \Div \lr{([\vr]_\delta +\delta)\vve(y)}|} .}
The energy estimate at this stage implies that $\delta
\int_{\T^d} \vert \vu\vert^2 \dx$ is controlled uniformly in time, moreover, $[\vr]_\delta $ is bounded. Since $\kappa_\ep \ast (-\Delta)^{-1} \nabla \Div$
is the operator of order zero,  we conclude that $\mathcal{M}|([\vr]_\delta +\delta)\vve(y)|\in L^\infty(0,T,L^2(\T^d))$.
Next, we observe that $\int K_h (z) |z| {\rm d}z$ is bounded by a constant uniformly in $h$, and so,
the whole boundary term is controlled uniformly in $h$ and $\eta$.  This reasoning applies to both boundary terms at $t=0$ and $t=T$. {In the final step we will divide this constant by  $|\ln h|$ and show that it disappears in the limit.}

Let us now focus  on the second term in \eqref{Ep:D1}:
\eq{\label{D1V}
&\iintTO{ K_{h}(x-y) \kappa_\ep\ast \Big((- \Delta_x)^{-1}\dv_x( ([\vr]_\delta +\delta)\vve)-(- \Delta_y)^{-1}\dv_y( ([\vr]_\delta +\delta)\vve)\Big) \pt W}.
}
Computing $\pt W$, we obtain
\eq{\label{Wt}
 \pt W=&\partial_t \sgn^\sigma_{xy} \vr(x) w(x) +\sgn^\sigma_{xy}  \pt(\vr(x) w(x))  .
}
Using the continuity equation and the transport equation for weights \eqref{def_w}, 
we obtain
\begin{equation*}
	(\vr w)_t + \Div(\vr w \vue) + k\vr^m w + \Lambda  \mathcal{M}[|\nabla \vu|] \vr w =0.
\end{equation*}
Then, tedious computations lead to 
\eq{\label{pt_sign}
\pt \sgn^\sigma_{xy}=&-\dv_x(\vue(x) \sgn_{xy}^\sigma) -
\dv_y(\vue(y) \sgn_{xy}^\sigma)  +(\sgn_{xy}^\sigma -  \vr(x)\partial \sgn^\sigma_{xy} )
\dv_x \vue(x) \\
& + (\sgn^\sigma_{xy} + \vr(y) \partial \sgn^\sigma_{xy}) \dv_y \vue(y)-k \partial \sgn_{xy}^\sigma (\vr^m(x)- \vr^m (y) ) \vr(x) w(x).
}
\begin{rmk}
    This is the crucial place explaining why we needed to smooth the sign function out. Indeed, the last term of the r.h.s of \eqref{pt_sign} is controlled by $\sigma^{-1}$.
\end{rmk}
Coming back to \eqref{Wt}, we get
\begin{align*}
	\pt W 
	=& - \Div_x \left(\vr(x) w(x)\vue(x) \sgn^\sigma_{xy} \right)-
	\dv_y(\vue(y) \sgn_{xy}^\sigma)  \vr(x)w(x) \\
	&- k \vr^m(x) w(x) \sgn^\sigma_{xy}- \Lambda  \mathcal{M}[|\nabla \vu|]  \vr(x) w(x) \sgn^\sigma_{xy} - \vr(x)\dv_x\vue(x) \partial \sgn_{xy}^\sigma \vr(x)w(x) \\
	&+ (\sgn^\sigma_{xy} + \vr(y) \partial \sgn^\sigma_{xy}) \dv_y \vue(y) \vr(x) w(x) -k \partial \sgn_{xy}^\sigma (\vr^m(x)- \vr^m (y) ) \vr(x) w(x).
\end{align*}
We now introduce the following notation: 
\eq{\label{ABC:defn}
	&\vc{A}(t,x,y)= -\vr(x) w(x)\vue(x) \sgn^\sigma_{xy} ,\\
	&\vc{B}(t,x,y)= - \vue(y) \sgn_{xy}^\sigma  \vr(x)w(x) ,\\
	&\mathcal{C}(t,x,y)= \Big(- k \vr^{m-1}(x)  \sgn^\sigma_{xy}- \Lambda  \mathcal{M}[|\nabla \vu|]  \sgn^\sigma_{xy} - \vr(x)\dv_x\vue(x) \partial \sgn_{xy}^\sigma\\
	&\quad \quad + (\sgn^\sigma_{xy} + \vr(y) \partial \sgn^\sigma_{xy}) \dv_y \vue(y) -k \partial \sgn_{xy}^\sigma (\vr^m(x)- \vr^m (y) )\Big) \vr(x) w(x)}
allowing us to write
\eq{\pt W=\dv_x \vc{A} +\dv_y \vc{B} +\mathcal{C}\label{ABC}.
}
Using the uniform bounds w.r.t. $ \ep $, especially  $ L^\infty_{t,x} $ bound for the density and $L^2_{t,x}$ bound for ${\cal M}[|\Grad\vu|]$, we justify that 
\begin{align*}
	&\vc{A} \in L^2(0,T;L^6(\T^d_x;L^\infty(\T^d_y)),\; \vc{B} \in L^2(0,T;L^6(\T^d_y;L^\infty(\T^d_x)),\\
	&\mathcal{C}\in L^2(0,T;L^2(\T^{d}_x;L^\infty(\T^d_y))
	+L^2(0,T;L^2(\T^{d}_y;L^\infty(\T^d_x))).
\end{align*}
Taking this into account, we obtain in a short form that
\eq{\nonumber
\pt W(t,x,y)\in L^2(0,T;W^{-1,3}(\T^d_x \times \T^d_y))= L^2(0,T;W^{-1,3}(\T^{2d})).
}
In particular, it follows that $C$ can be written as a  divergence of  $L^3$ vector fields  modulo some constant, so, with slight abuse of notation, we can write
\eq{ \label{boundABC}
\partial_t W=\Div_x \vc{A}+\Div_y\vc{B}+\{\mathcal{C}\},\quad
\text{where}\quad
\vc{A},\ \vc{B}\in L^2_tL^3_{x,y} ,\quad  \{\mathcal{C}\}=\int_{\T^{2d}} \mathcal{C}\, \dx\, \dy.}

With this observation at hand, we return to the estimate of \eqref{D1V}, and we have
\eqh{
	&\iintTO{ K_{h}(x-y) \kappa_\ep\ast \Big((- \Delta_x)^{-1}\dv_x( ([\vr]_\delta +\delta)\vve)-(- \Delta_y)^{-1}\dv_y( ([\vr]_\delta +\delta)\vve)\Big) \pt W}\\
	&=\iintTO{\! K_{h}(x-y) \kappa_\ep\ast \Big((- \Delta_x)^{-1}\dv_x( ([\vr]_\delta +\delta)\vve)-(- \Delta_y)^{-1}\dv_y( ([\vr]_\delta +\delta)\vve)\Big)\\
	&\hspace{12cm}( \Div_x\vc{A}+ \Div_y \vc{B})}\\
	&=-\iintTO{\!\! \Grad_x K_{h}(x-y) \kappa_\ep\ast \Big((- \Delta_x)^{-1}\dv_x( ([\vr]_\delta +\delta)\vve)-(- \Delta_y)^{-1}\dv_y( ([\vr]_\delta +\delta)\vve)\Big) \vc{A}}\\
	&\quad -\iintTO{ K_{h}(x-y) \Grad_x(- \Delta_x)^{-1}\dv_x(\kappa_\ep\ast  ([\vr]_\delta +\delta)\vve)\cdot \vc{A}}\\
	&\quad +\iintTO{ \Grad_y K_{h}(x-y) \kappa_\ep\ast \Big((- \Delta_x)^{-1}\dv_x( ([\vr]_\delta +\delta)\vve)-(- \Delta_y)^{-1}\dv_y( ([\vr]_\delta +\delta)\vve)\Big) \vc{B}}\\
	&\quad +\iintTO{ K_{h}(x-y) \Grad_y(- \Delta_y)^{-1}\dv_y(\kappa_\ep\ast  ([\vr]_\delta +\delta)\vve)\cdot \vc{B}}\\
	& = \Sigma_{i=1}^4 J_i.
}
To estimate $J_1$,  we first recall that $|x-y||\Grad K_h(x-y)|\leq CK_h$. Therefore, introducing the notation $$ \mathscr{S}=\kappa_\ep\ast (- \Delta)^{-1}\dv( ([\vr]_\delta +\delta)\vve),$$ we obtain 
\eqh{
	|J_1|\leq &\iintTO{ \vert \Grad_x K_{h}(x-y)\vert\;  \vert \mathscr{S}_x-\mathscr{S}_y\vert \vc{A}}\\
	\leq&\iintTO{K_{h}(x-y)\frac{|\mathscr{S}_x-\mathscr{S}_y|}{|x-y|} |\vc{A}|}.
}
Using again the formula \eqref{formula:D} from Lemma \ref{Lagrange}, we further deduce
\eqh{
	|J_1|\leq& C\iintTO{K_{h}(x-y)\lr{{\cal M}[|([\vr]_\delta +\delta)\vve(x)|]+{\cal M}[|([\vr]_\delta +\delta)\vve(y)|] }|\vc{A}|}.
}
Recall that  $\vve$ is  bounded in $ L^2_tL^6_x $ uniformly w.r.t. $ \ep $. Moreover, we can interpolate with \eqref{eta_teta} to obtain 
\eqh{
\|\vve\|_{L^p_tL^q_x}\leq C\eta^\vt
}
for $1\leq p< 2, \ 1\leq q<6$. Since $\frac{1}{6}+\frac{1}{3}<1$, the bound \eqref{boundABC} for $\vc{A}$ is enough to conclude that
\begin{align}\nonumber
	|J_1|\leq C \eta^\vt \|K_h\|_{L^1_x},
\end{align}
for some generic constant $\vt>0$.\\
The estimate of $J_2$ is similar, and it gives the same control in terms of $h$. Finally, the estimates of $ J_3 $ and $ J_4 $ follow the same argument with $x$ replaced by $y$.

Now we check the part of $A_3^{(iii)}$ in \eqref{intD-D} corresponding to $ \mathcal{D}_2(\vv_\eta) $ defined in \eqref{def:D}.  We  use the H\"older inequality and \eqref{eta_teta} (remembering that $\delta$ is fixed), to deduce
\eqh{
&\iintTO{ K_{h}(x-y) \kappa_\ep\ast \Big(\calD_2(\vve(x))-\calD_2(\vve(y))\Big) W}\\
&=\iintTO{ K_{h}(x-y) \kappa_\ep\ast  \Big((- \Delta_x)^{-1}\dv_x  \Div_x([\vr\vue]_\delta\otimes \vve)\\
& \hspace{5cm}-(- \Delta_y)^{-1}\dv_y  \Div_y([\vr\vue]_\delta\otimes \vve)\Big) W }\\
&\leq C \eta^\vt\|K_h\|_{L^1_x}.
}
Next, we can estimate the term related to $\mathcal{D}_3(\vv_\eta)$. It is also quite simple at this level, since $(- \Delta)^{-1}\dv $ gains one derivative. Therefore, \eqref{eta_teta} implies
\begin{align*}
    \iintTO{ K_{h}(x-y) \kappa_\ep\ast \Big(\calD_3(\vve(x))-\calD_3(\vve(y))\Big) W}
   \leq  C(M,\delta,\sigma)\eta^\vt\|K_h\|_{L^1_x}.
\end{align*}
Finally, we add up all the estimates for $A_3^{(iii)}$ and obtain
\eq{\label{a33b}
\int_0^T A_3^{(iii)} {\rm d}t &=- \frac{1}{2\mu+\xi} \iintTO{ K_{h}(x-y) \kappa_\ep\ast \Big(\calD(x)-\calD(y)\Big)  
\sgn^\sigma_{xy} \vr(x) w(x)}\\
&\leq C \|K_h\|_{L^1_x} \left( h^\beta \eta^{-\alpha} +
     \eta^\vt \right).
}

\subsubsection*{{Estimate of $A_4$ and $A_5$}}

In order to estimate the terms $A_4$ and $A_5$ we note that by integrability of $\nabla \vue$ in $L^2_{t,x}$, the definition of  $|\cdot|^\sigma$, and the inequality \eqref{eq:37} we obtain
\begin{equation}\label{a4a5}
    \int_0^T( |A_4|+|A_5| ) \;{\rm d}t \leq C \sigma,
\end{equation}
where the constant $C$ depends only on the initial data.

\subsection{Proof of Lemma \ref{Lemma:Kol_Mod} -- conclusion}\label{Ep:comp3}
Summarising the previous section, we can now integrate \eqref{eqdtR} w.r.t. time and use estimates \eqref{a1a2}, \eqref{a31}, \eqref{a32},\eqref{a33a},\eqref{a33b} and \eqref{a4a5} to obtain 
\eq{\label{St-S0-1}
&
 R_{h}^{\sigma}(T) - R_{h}^{\sigma}(0)\\
&\leq C    \iintTO{  \Ov K_{h}(x-y)\vert{\vr(y)-\vr(x)}\vert^\sigma   w(x)}\\
&\quad+ C \iintTO{  \Ov K_{\ep}(x-y)\vert{\vr(y)-\vr(x)}\vert^\sigma   w(x)} \\
&\quad +C\iintTO{\Ov K_h(x-y) \lr{D_{|x-y|}\vue(x)-D_{|x-y|}\vue(y)} \vert{\vr(x)-\vr(y)}\vert^\sigma w(x)}\\
&\quad+ C(h_0^{\beta}\eta^{-\alpha}+\eta^\vt)+C\sigma .
}
To handle the third term we apply the following result for the components of vector $\vu$.
\begin{lemma}\label{Lemma:Du}
Let $u \in W^{1,2}(\T^d)$ then
\begin{equation}\nonumber
\int_{\T^d} \frac{{\rm d}z}{(|z|+h)^d} \| D_{|z|}u(\cdot) - 
D_{|z|}u(\cdot + z)\|_{L^2} \leq C|\ln h|^{1/2} \|u\|_{W^{1,2}(\Td)}.
\end{equation}
\end{lemma}
The proof of this lemma is postponed to the Appendix \ref{Lemma13}, see also the original paper of Bresch and Jabin \cite[Lemma 6.3]{BJ}.
\medskip

Changing variables for $z=x-y$ in \eqref{St-S0-1} we have
\eq{\nonumber
&\iintO{\Ov K_h(x-y) \lr{D_{|x-y|}\vue(x)-D_{|x-y|}\vue(y)} |\vr(y)-\vr(x)|^\sigma w(x)}\\
&=\int_{\T^{2d}}{\Ov K_h(z) \lr{D_{|z|}\vue(y+z)-D_{|z|}\vue(y)} |\vr(y)-\vr(y+z)|^\sigma w(z+y)}\, {\rm d}z\,\dy\, \\
&\leq C(M) \int_{\T^d}{\Ov K_h(z) \|D_{|z|}\vue(\cdot+z)-D_{|z|}\vue(\cdot)\|_{L^2_x}}\, {\rm d}z\,\\
&\leq C(M) |\log h|^{-1/2}\|\vu(t)\|_{W^{1,2}_x}.
}
Let us now re-introduce the notation of the  $\ep-$labelled  sequences $\vr^\ep$ and $w^\ep$, in particular, we denote
\eq{
R_{h}^{\ep,\sigma}(T)=\iintO{\Ov K_h(x-y) |\vr^\ep(x)-\vr^\ep(y)|^\sigma w^\ep(x)}.
}
 Taking $\lim\sup_{\ep\to0}$ in \eqref{St-S0-1} we thus get
\eq{\label{St-S02}
&\limsup_{\ep\to0} R_{h}^{\ep,\sigma}(T) -\limsup_{\ep\to0} R_{h}^{\ep,\sigma}(0)\\
&\leq C \limsup_{\ep\to0}   \iintTO{  \Ov K_{h}(x-y)|\vr^\ep(y)-\vr^\ep(x)|^\sigma   w^\ep(x)} \\
&\quad+ C \limsup_{\ep\to0}  \iintTO{  \Ov K_{\ep}(x-y)|\vr^\ep(y)-\vr^\ep(x)|^\sigma   w^\ep(x)} \\
&\quad +C(M) |\log h|^{-1/2}+ C(h^{\beta}\eta^{-\alpha}+\eta^\vt+\sigma),
}
where $C$ depends now also on $T$.
Observe that the second term on the r.h.s. can be estimated as follows
\eq{\nonumber
& \limsup_{\ep\to0} \iintTO{  \Ov K_{\ep}(x-y) |\vr^\ep(y)-\vr^\ep(x)|^\sigma    w^\ep(x)}\\
&\leq \int_0^T\limsup_{\ep\to0}\iintO{  \Ov K_{\ep}(x-y) 
|\vr^\ep(y)-\vr^\ep(x)|^\sigma  w^\ep(x)}\,{\rm d} t\\
&\leq  \int_0^T\sup_{h\leq h_0} \limsup_{\ep\to0}\iintO{  \Ov K_{h}(x-y) |\vr^\ep(y)-\vr^\ep(x)|^\sigma   w^\ep(x)}\,{\rm d} t,
}
for some fixed $1>h_0>0$. Due to the assumption on the initial data \eqref{id-approx} we have
\eq{\nonumber
\limsup_{\ep\to0} R^{\ep,\sigma}_{h}(0)\to 0,
}
as $h\to 0$.  So, taking $\sup_{h\leq h_0}$ on both sides of \eqref{St-S02} we get
\eq{ \label{jaktuwywalicM} 
&\sup_{h\leq h_0}\limsup_{\ep\to0} R_{h}^{\ep,\sigma}(T) \\
&\leq C  \int_0^T  \sup_{h\leq h_0} \limsup_{\ep\to0}  \iintO{  \Ov K_{h}(x-y) |\vr^\ep(y)-\vr^\ep(x)|^\sigma   w^\ep(x)} \,{\rm d} t\\
&\quad +C |\log h_0|^{-1/2}+ C(h_0^{\beta}\eta^{-\alpha}+\eta^\vt+\sigma) +o(h_0)\\
&\leq C   \int_0^T \sup_{h\leq h_0}\limsup_{\ep\to0} R_{h}^{\ep,\sigma}(t)\,{\rm d}t  +C(M) |\log h_0|^{-1/2}+ C(h_0^{\beta}\eta^{-\alpha}+\eta^\vt+\sigma) +o(h_0).
}
Applying the Gr\"onwall's inequality we get
\eq{ \label{jtw2}
\limsup_{\ep\to0} R_{h_0}^{\ep,\sigma}(T)& \leq\sup_{h\leq h_0}\limsup_{\ep\to0} R_{h}^{\ep,\sigma}(T) \\
&\leq \lr{C(M)|\log h_0|^{-1/2}+ C(h_0^{\beta}\eta^{-\alpha}+\eta^\vt+\sigma) +o(h_0) }e^{CT}.
}
Choosing $\eta=h_0^{\frac{\beta}{2\alpha}}$ we see that the r.h.s. tends to $0$ with $h_0\to0$. {Finally, since in the l.h.s of (\ref{jtw2}) $C$ above does not depend on $\sigma$}, we  pass to the limit $\sigma \to 0$, we
obtain
$$
\limsup_{\ep\to0} \int_{0}^{T} \!\!\!\int_{\T^{2d}} \Ov K_{h_0}(x-y) |\vr^\ep(x)-\vr^\ep(y)| w^\ep(x\,)\dx\,\dy\,\dt \to 0 \mbox{ \ as \ } h_0 \to 0.
$$ 
This is precisely the statement of Lemma \ref{Lemma:Kol_Mod}, the proof is now complete. $\Box$

\subsection{Kolmogorov criterion -- removal of the weights}\label{Ep:comp4}

Our goal now is to remove the weight $w^\ep(x)$ from \eqref{criterion00}, or its equivalent form \eqref{criterion0}, where both of the weights $w^\ep(x), w^\ep(y)$ are present. We prove the following result.
\begin{lemma}\label{Lemma:Kol_rem_weight}
	The sequence $ \{\vr^\ep \}_{\ep>0} $ of approximative solutions to system (\ref{maine}) satisfies
	\eq{\label{KolCri-rem}
		\underset{\ep\to0}{\lim\sup} \left(  \int_0^T\!\!\!\iintO{
			\Ov K_h(x-y)|\vr^\ep(t,x)-\vr^\ep(t,y)|} \, \dt\right) \to 0, \ \mbox{ as } h\to 0.
	}
\end{lemma}
\pf
Let $\zeta<1$. We define $\Omega_\zeta = \{x: w^\ep(t,x) \leq \zeta\}$ and denote by 
$\Omega_\zeta^c$ its complementary. 
\eq{\nonumber
&\underset{\ep\to0}{\lim\sup} \left(  \int_0^T\!\!\!\iintO{
 \Ov K_h(x-y)|\vr^\ep(t,x)-\vr^\ep(t,y)|} \, \dt\right)\\
 &=
 \underset{\ep\to0}{\lim\sup} \left(  \int_0^T\!\!\!\int_{x\in\Omega^c_\zeta \cup y\in\Omega^c_\zeta}
 \Ov K_h(x-y)|\vr^\ep(t,x)-\vr^\ep(t,y)| \dx\dy\, \dt\right)\\
&\quad + \underset{\ep\to0}{\lim\sup} \left(  \int_0^T\!\!\!\int_{x\in\Omega_\zeta \cap y\in\Omega_\zeta}
 \Ov K_h(x-y)|\vr^\ep(t,x)-\vr^\ep(t,y)| \, \dx\,\dy\, \dt\right)\\
 &=I_1+I_2.
 }
 The first term is easy, note that
\eq{
I_1&= \underset{\ep\to0}{\lim\sup} \left(  \int_0^T\!\!\!\int_{x\in\Omega^c_\zeta \cup y\in\Omega^c_\zeta}
 \Ov K_h(x-y)|\vr^\ep(t,x)-\vr^\ep(t,y)| \dx\dy\, \dt\right)\\
 &\leq \underset{\ep\to0}{\lim\sup} \frac{1}{\zeta} \left(  \int_0^T\!\!\!\int_{x\in\Omega^c_\zeta \cup y\in\Omega^c_\zeta}
 \Ov K_h(x-y)|\vr^\ep(t,x)-\vr^\ep(t,y)| (w^\ep(t,x){+w^\ep(t,y)})\dx\dy\, \dt\right) .
}
For the term $I_2$ we use that both $\zeta$ and $\omega$ are less than 1, and so $|\log w|\geq |\log\zeta|$, and so
\eq{\label{Remove-w-I2}
I_2&=  \underset{\ep\to0}{\lim\sup} \left(  \int_0^T\!\!\!\int_{x\in\Omega_\zeta \cap y\in\Omega_\zeta}
 \Ov K_h(x-y)|\vr^\ep(t,x)-\vr^\ep(t,y)| \, \dt\right)\\
 &\leq \underset{\ep\to0}{\lim\sup} \left(  \int_0^T\!\!\!\int_{x\in\Omega_\zeta \cap y\in\Omega_\zeta}
 \Ov K_h(x-y)|\vr^\ep(t,x)-\vr^\ep(t,y)|\frac{|\log w^\ep(t,x)|}{|\log\zeta|}  \dx\, \dy\, \dt\right)\\
 &\leq \sup_{\ep>0} \frac{2}{|\log\zeta|}\iintTO{ \overline{K}_h(x-y)\vr^\ep(t,x)|\log w^\ep(t,x)| }\\
 & \leq \sup_{\ep >0}  \frac{2}{|\log\zeta|} \int_0^T\!\!\!\int_{\T^d}{ \vr^\ep(t,x)|\log w^\ep(t,x)| }\\
 & \leq \frac{C}{|\log\zeta|},
}
By estimate \eqref{boundlogw} from Proposition \ref{prop:w} we conclude the proof. $\Box$

\subsection{Limit passage $\ep\to0$ in the equations}
Thanks to  compactness of the sequence approximating the density, we are able to pass to the limit in the pressure term. After letting $\ep\to 0$ we end up with the approximate system:
\begin{subnumcases}{\label{maine2epboth}}
	\pt\vr+\Div (\vr \vu) +k  \vr^m= 0,\label{maine1ep}\\[10pt]
	\pt((\delta+[\vr]_\delta)\vu)+\Div([\vr \vu]_\delta
	\otimes\vu) - \Div \vS(\vu)+ \Grad p^{M,\lambda}(\vr)+k [\vr^m]_\delta \vu =\vc{0},\label{maine2ep}
\end{subnumcases}
satisfied in the sense of distributions. Note that at this stage the density is pointwisely  bounded, thus the renormalized  approximate continuity equation is satisfied automatically.

\section{Compactness and passage to limit for $M\to\infty$}\label{Sec:5} 
In this section  the parameters $k,\; \delta\; $ and $\lambda$ are fixed, and we study the limit passage when $ M \rightarrow \infty$. We skip the upper index $M$ when no confusion can arise, in particular, we denote $ \vr^{M}=\vr $ and $ \vu^{M}=\vu  $. 

\subsection{Uniform estimates}
To  establish uniform estimates w.r.t. $ M $, we need to impose certain conditions on $ m, \Gamma $, and $ \gamma $. They will be introduced in several places below, where relevant.

To begin, we take the limit $\ep\to0$ in energy equality \eqref{energy} to obtain
\eq{\label{energy:M}
&\intO{\lr{(\delta+[\vr]_\delta) \frac{|\vu|^2}{2}+\vr e^{M,\lambda}(\vr)}(\tau)}
+\inttO{\vS(\Grad \vu):\Grad\vu}\\
&\qquad+k \inttO{\vr^m \frac{\lambda\Gamma}{\Gamma-1}\vr_M^{\Gamma-1}}+k\inttO{ \vr^m \frac{\gamma}{\gamma-1}\vr_M^{\gamma-1}}+\frac{k}{2}
\inttO{ [\vr^m]_\delta |\vu|^2}\\
&\leq \intO{\lr{(\delta+[\vr_{0}^{\lambda}]_\delta) \frac{|\vu_{0}^{\lambda}|^2}{2}+\vr_0^{\lambda} e^{M,\lambda}(\vr_{0}^{\lambda})}}
}
for a.e. $ \tau \in (0,T)$. From the approximation of initial data \eqref{id-approx}, we infer that the r.h.s of \eqref{energy:M} is uniformly bounded w.r.t. $ M $.
This change of equality to inequality is a consequence of lower semi-continuity of convex functions. Note however, that the resulting inequality does not provide us with useful estimates for $\vr$ only for the truncation $\vr_M$. Therefore, we use the continuity equation directly. Multiplying \eqref{maine1ep} by $\vr$ and integrating by parts we get
\eq{\label{55}
\frac12\Dt\intO{\vr^2}+k\intO{\vr^{m+1}}+\frac12\intO{\vr^2\Div\vu}=0.
}
If $m+1\geq4$ this equality gives us uniform in $M$ bound for $\vr$ in $L^{m+1}((0,T)\times \T^d)$. For reasons that will become clear in the subsequent limit passages we also assume
\begin{align}\label{m, Gamma : M,m}
	m+1 > \Gamma \geq m .
\end{align} 
Summarising the above, for fixed $ k, \lambda$ and $ \delta $, we have the following uniform bounds w.r.t. $ M $: 
\begin{align}\label{est:M}
	\begin{split}
	&\|\vr\|_{L^\infty_t L^2_x} +\|\vr\|_{L^{m+1}_{t,x}}\leq C,\\
	&\|\vu\|_{L^2_tW^{1,2}_x}+\|\vu\|_{L^\infty_tL^2_x}\leq C,\\
	&\|[\vr]_\delta \vert \vu \vert^2\|_{L^\infty_tL^1_x}\leq C,\\
	& \|\vr \vu\|_{L^2_tL^{\frac{3}{2}}_x} + \|\vr \vu\|_{L^{\infty}_tL^1_x}\leq C,\\
	& \|[\vr]_\delta \vu\|_{L^\infty_tL^\frac{4}{3}_x} +  \|[\vr]_\delta \vu\|_{L^{m+1}_t L^\frac{2(m+1)}{m+3}_x}\leq C. 
		\end{split}
\end{align}
 Note that if $m+1>\Gamma>\gamma$ we may pass to the limit $M\to\infty$ in the pressure term of \eqref{maine2ep}. However, to identify the limit, we need to go through the Kolmogorov compactness criterion again. 
 
 \subsection{Compactness criterion -- the main estimate} \label{cc:M}
 The only changes in the main estimate from the previous section appear in the estimate of $A_3$. We discuss the details below.
 
First, since there is no regularisation $\kappa_\ep$, we only need to control
\eq{\label{A3_sum:M}
A_{3}
=& - \frac{1}{2\mu+\xi} \iintO{ K_{h}(x-y)  \lr{p^{M,\lambda}(x)-p^{M,\lambda}(y)} \sgn^\sigma_{xy}\; \vr(x) w(x)}\\
& - \frac{1}{2\mu+\xi} \iintO{ K_{h}(x-y)  \Big(\calD(\vu(x))-\calD(\vu(y))\Big) \sgn^\sigma_{xy} \; \vr(x) w(x)}\\
=&A_3^{(i)}+A_3^{(ii)},
}
where the definition of $ \mathcal{D}(\vu) $ is given in \eqref{def:D}.

The first term $ A_3^{(i)} $ has a good sign due to the monotonicity of $p^{M,\lambda}$, so it can be moved to the l.h.s of the main estimate. The difficulty is thus only in the second term $A_3^{(ii)} $. Recall that the estimate of an analogous term  in the previous step (term $A_3^{(iii)}$ in \eqref{intD-D}) required the $L^\infty$ bound for the density, which is no longer available. 

The main idea here is to use again Lemma \ref{lem:approx} along with Lemma \ref{interpolation-lemma} formulated below. To do this, we again consider the time-space regularization  of $ \vu $ denoted by $ \pi_{\eta}\ast\vu $. We split
\begin{align*}
	\vu = \lr{\vu-\pi_{\eta} \ast\vu} +  \pi_{\eta}\ast\vu  = \vc{v}_{\eta} + \pi_{\eta}\ast\vu , 
\end{align*}
consequently
\begin{align*}
    -({2\mu+\xi}) A_3^{(ii)}=&\iintO{ K_{h}(x-y)  \Big(\calD(\pi_{\eta}\ast\vu(x))-\calD(\pi_{\eta}\ast\vu (y))\Big) \sgn^\sigma_{xy} \; \vr(x) w(x)}\\
    +&\iintO{ K_{h}(x-y)  \Big(\calD(\vc{v}_{\eta}(x))-\calD(\vc{v}_{\eta}(y))\Big) \sgn^\sigma_{xy} \; \vr(x) w(x)}\\
    =&I_{1}+ I_{2}.
\end{align*}
{\emph{Estimate of $I_{1}$.}} Using the continuity equation for $ [\vr]_{\delta} $, we have
\begin{align}\label{est:M:I2}
\begin{split}
	 &\iintTO{ K_{h}(x-y)  \Big(\calD( \pi_{\eta}\ast\vu(x))-\calD(\pi_{\eta}\ast\vu(y) )\Big) \sgn^\sigma_{xy} \; \vr(x) w(x)}\\
	 &=\iintTO{ K_{h}(x-y) \Big((- \Delta_x)^{-1}\dv_x( [\vr]_\delta^\delta\pt\pi_{\eta}\ast\vu)-(- \Delta_y)^{-1}\dv_y( [\vr]_\delta^\delta\pt\pi_{\eta}\ast\vu)\Big) W(t,x,y)}\\
	 &+\iintTO{ K_{h}(x-y)  \Big((- \Delta_x)^{-1}\dv_x( [\vr \vu]_\delta\cdot\Grad_x\pi_{\eta}\ast\vu)\\
	 &\hspace{4cm}-(- \Delta_y)^{-1}\dv_y( [\vr \vu]_\delta\cdot\Grad_y\pi_{\eta}\ast\vu)\Big) W(t,x,y)}\\
	 &\leq C\|W(t,x,y)\|_{\infty} h^\beta\eta^{-\alpha}
	 \end{split}
\end{align}
for some $ \alpha, \beta >0 $, where we used an abbreviate notation $ [\vr]_\delta^\delta= \delta+[\vr]_\delta.$

\bigskip

\noindent {\emph{Estimate of $I_{2}$.}} We first estimate the term \[\iintTO{ K_{h}(x-y)  \Big(\calD_{1}(\vc{v}_{\eta} (x)) -\calD_{1}(\vc{v}_{\eta}  (y))\Big) \sgn^\sigma_{xy} \; \vr(x) w(x)},\] 
with $\calD_{1}$ as defined in \eqref{def:D}.

This estimate is similar to \eqref{Ep:D1} from the previous section. We might further split:
\eq{&\iintTO{ K_{h}(x-y)  \Big(\calD_1(\vve(x))-\calD_1(\vve(y))\Big) W(t,x,y)}\\
	&=
	\iintTO{ K_{h}(x-y) \Big((- \Delta_x)^{-1}\Div_x( \pt ([\vr]_\delta^\delta\vve))-(- \Delta_y)^{-1}\Div_y( \pt ([\vr]_\delta^\delta\vve))\Big) W(t,x,y)}\\
	&=\iintTO{ K_{h}(x-y) (- \Delta_x)^{-1}\Div_x( \pt ([\vr]_\delta^\delta\vve)) W(t,x,y)}\\
	&\quad - \iintTO{ K_{h}(x-y) ((- \Delta_y)^{-1}\Div_y( \pt ([\vr]_\delta^\delta\vve)) W(t,x,y)} \\
	&={I}_{2,1}+{I}_{2,2}.\label{M:D1}
}
Unfortunately, we do not have the sufficient regularity to perform integration by parts as done in \eqref{Ep:D1}. Instead of this, we state an abstract interpolation lemma, mentioned in the begining of our paper.
\begin{lemma}\label{interpolation-lemma}
	Let $ \pmb{\varphi} \in L^q(\R^{d+1}) $ with $ \partial_t \pmb{\varphi} \in L^{q^\prime} (\R; W^{-1,q^\prime} (\R^d)) $ and let $ W \in L^{\bar{q}}(\R^{d+1} )  $ with $ \partial_t W \in  L^{{\bar{q}}^\prime} (\R; W^{-1,{\bar{q}}^\prime} (\R^d))  $ such that $ 1\leq q, q^\prime, \bar{q}, \bar{q}^\prime \leq \infty $ with either one of the following relations satisfied:
	\begin{align}\label{il:a1}
		\begin{split}
			&(1) \; \qquad\frac{1}{q^\prime}+\frac{1}{\bar{q}}<1<	\frac{1}{q}+\frac{1}{\bar{q}^\prime },\\
			& (2)\; \qquad\frac{1}{q^\prime}+\frac{1}{\bar{q}}>1>	\frac{1}{q}+\frac{1}{\bar{q}^\prime }.
		\end{split}
	\end{align}
	Then there exists $ 0<\alpha<1 $ such that 
\eq{
		\int_{\mathbb{R}} \int_{\mathbb{R}^d}  \Big((- \Delta_x)^{-1}\dv_x \pt \pmb{\varphi} \Big) W(t,x)\, \dx\; \dt 
		\leq C \Vert \pmb{\varphi}  \Vert_{L^q_{t,x}}^{\alpha} \Vert \partial_t \pmb{\varphi} \Vert_{L^{{q}^{\prime}}_tW^{-1,{q}^\prime}_x}^{1-\alpha}\Vert W \Vert_{L^{\bar{q}}_t}^{1-\alpha} \Vert \partial_t W \Vert_{L^{\bar{q}^{\prime}}_t W^{-1,\bar{q}^\prime}_x}^{\alpha}.
}
\end{lemma}
 The above lemma is a key result concerning the compactness of the effective viscous  flux. The proof for the whole space case is presented in the Appendix \ref{Lemma15}. To adopt these considerations to the periodic setting we de efine a partition of the unity for $\T^d$,  say $\{\psi_k\}$, such that each of supports of $\psi_k$ in contained in one map of the torus. Then functions $\psi_k\pmb{\varphi}$ are defined over the space $\R^d$ by the canonical embedding. Since $\psi_k$ are smooth, the functions  $\psi_k\pmb{\varphi}$ have the required regularity $W^{-1,q'}_x$. To extend the problem to the whole time line, we can cut the time endpoints by a suitable compactly supported function, which will give us assumptions from the lemma up to an arbitrarily small correction.
	
\medskip

To apply Lemma \ref{interpolation-lemma}, we recall
\begin{align*}
\pt W=\dv_x \vc{A} +\dv_y \vc{B} +\mathcal{C},
\end{align*}
where 
$ \vc{A},\;\vc{ B} $ and $ \mathcal{C} $ are given by \eqref{ABC:defn}.
We point out that in our previous estimates for $ \mathcal{C} $ in \eqref{ABC} we used the uniform $ L^\infty $ bound of the density which we no longer have. Therefore, here we use
\begin{align}\label{M:ABC}
	\begin{split}
	&\vc{A} \in L^2(0,T;L^6(\T^d_x;L^\infty(\T^d_y)),\; \vc{B} \in L^2(0,T;L^6(\T^d_y;L^\infty(\T^d_x)),\\
	& \mathcal{C} \in L^1(0,T; L^1 (\T^d_x; L^\infty (\T^d_y)))+ L^1(0,T; L^1 (\T^d_y; L^\infty (\T^d_x)))
\end{split}
\end{align}
with
\begin{align*}
	&\Vert \mathcal{C} \Vert_{L^1_tL^1_x L^\infty_y+ L^1_t L^1_y L^\infty_x}
	 \leq C\left(k, \frac{1}{\sigma} \right)\; \Vert \vr w \Vert_{L^\infty_{t,x}} \left(\Vert \vr \Vert_{L^2_t L^2_x} \Vert \nabla_x \vu  \Vert_{L^2_t L^2_x} + \Vert \vr \Vert_{L^2_t L^2_y} \Vert \nabla_x \vu  \Vert_{L^2_t L^2_y}  \right).
\end{align*}
Using $m$-dependent density estimates one can verify that  $ \mathcal{C} $ is little bit more integrable, i.e.  
$$\mathcal{C} \in L^2(0,T; L^q (\T^d_x; L^\infty (\T^d_y)))+ L^2(0,T; L^q(\T^d_y; L^\infty (\T^d_x))), $$ 
with $ 1<q=\min\{\frac{m+1}{m}, \frac{2(m+1)}{m+3}\} $, although 
we continue with the $L^1_tL^1_x L^\infty_y+ L^1_t L^1_y L^\infty_x$-norm of $\mathcal{C}$.

Next, we introduce 
 \begin{align}\label{M:w1}
	\mathscr{W}_1(t,x)= \int_{\T^d} K_h(x-y) W(t,x,y) \; \dy,\quad \mathscr{W}_2(t,y)= \int_{\T^d} K_h(x-y) W(t,x,y) \; \dx 
\end{align} 
and rewrite $I_{2,1}$ and $I_{2,2}$ from \eqref{M:D1} as 
\eqh{
I_{2,1}= \int_0^T\!\!\! \int_{\Td} (- \Delta_x)^{-1}\Div_x( \pt ([\vr]_\delta^\delta\vve)) \mathscr{W}_1\, \dx \; \dt ,\\
I_{2,2}= \int_0^T\!\!\! \int_{\Td} (- \Delta_y)^{-1}\Div_y( \pt ([\vr]_\delta^\delta\vve)) \mathscr{W}_2\, \dy \; \dt.
}
Note that we have
\begin{align*}
	\Vert 	 \mathscr{W}_1  \Vert_{L^\infty_{t,x}} +\Vert 	 \mathscr{W}_2  \Vert_{L^\infty_{t,y}} \leq \Vert K_h \Vert_{L^1_x} \Vert W  \Vert_{L^\infty_{t,x,y}} .
\end{align*}
Let us now focus on the estimate of $I_{2,1}$. A direct calculation yields
\begin{align*}
	\pt \mathscr{W}_1= & \int_{\T^d} K_h(x-y) \pt W(t,x,y) \; \dy \\
	&=  \int_{\T^d} \lr{(K_h(x-y) \Div_x \vc{A}(t,x,y)+K_h(x-y) \Div_y \vc{B}(t,x,y) + K_h(x-y) \mathcal{C}(t,x,y)) }\; \dy\\
	&= \Div_x  \lr{\int_{\T^d} K_h(x-y) \vc{A}(t,x,y) \dy }+  \int_{\T^d} K_h(x-y) \mathcal{C}(t,x,y) \; \dy .
\end{align*}
Now, using the bound from \eqref{M:ABC}, we obtain the uniform bound of $ \int_{\T^d} K_h(x-y) \vc{A}(t,x,y) \,\dy $ in $  L^2_tL^6_x $  and $  \int_{\T^d} K_h(x-y) \mathcal{C}(t,x,y) \,  \dy$ in $ L^1_{t,x}  $ with the following estimates:
\begin{align*}
&	\left\Vert  \int_{\T^d} K_h(x-y) \vc{A}(t,x,y)\; \dy \right\Vert_{  L^2_tL^6_x } \leq \Vert K_h \Vert_{L^1_x} \Vert \vc{A}\Vert_{ L^2_tL^6_xL^\infty_y}, \\
&\left\Vert \int_{\T^d} K_h(x-y) \mathcal{C}(t,x,y) \;\dy  \right \Vert_{L^1_{t,x}} \leq \Vert K_h \Vert_{L^1} \Vert \mathcal{C} \Vert_{L^1_tL^1_xL^\infty _y+ L^1_tL^1_yL^\infty_x}.
\end{align*}
Noticing that \begin{align*}
   \left\Vert  \Div \lr{\int_{\T^d} K_h(x-y) \vc{A}(t,x,y)\; \dy} \right\Vert_{  L^2_tW^{-1,6}_x } \leq  \left\Vert  \int_{\T^d} K_h(x-y) \vc{A}(t,x,y)\; \dy \right\Vert_{  L^2_tL^6_x} ,
\end{align*}
and 
\begin{align*}
     L^1_{t,x} \hookrightarrow L^1_tW^{-1,1}_x ,
\end{align*}
 we have 
\begin{align*}
\Vert 	\pt \mathscr{W}_1 \Vert_{ L^2_tW^{-1,6}_x +L^1_tW^{-1,1}_x}\leq \Vert K_h \Vert_{L^1_x}\lr{\Vert \vc{A}\Vert_{ L^2_tL^6_xL^\infty_y}+ \Vert \mathcal{C} \Vert_{L^1_t L^1_xL^\infty_y+ L^1_tL^1_y L^\infty_x}}.
\end{align*}
To summarise, we can write that
\eq{\nonumber
\mathscr{W}_1 \in  L^\infty_{t,x} \text{ with } \pt  \mathscr{W}_1 \in L^1_t W^{-1,1}_x,}
with the norms uniformly bounded w.r.t. $M$.

 On the other hand, recalling the momentum equation
\eq{\nonumber
	\pt( [\vr]_\delta^\delta \vu)=-\Div( [\vr \vu]_\delta \otimes\vu) + \Div \vS(\vu)- \Grad p^{M,\lambda}(\vr) -k [\vr^m]_\delta \vu
} and using uniform estimates \eqref{est:M}, we deduce the following estimates for the r.h.s of above:
\begin{align}\label{mom-est-M}
	\begin{split}
		&\Vert [\vr \vu]_\delta \otimes\vu \Vert_{L^2_t L^{\frac{12}{11}}_x} \leq \Vert [\vr \vu]_\delta  \Vert_{L^\infty_x L^{\frac{4}{3}}_x} \Vert \vu \Vert_{L^2_t L^{6}_x}\leq C ,\\
		&\Vert \vS(\vu)\Vert_{L^2_t L^2_x} \leq C,\\
		&\Vert p^{M,\lambda}(\vr) \Vert_{L^q_t L^q_x} \leq \Vert \vr^\Gamma \Vert_{L^q_t L^q_x} \leq \Vert \vr \Vert_{L^{m+1}_{t,x}} \leq C, \quad\text{ for }\  1<q<\frac{m+1}{\Gamma},\\
		&\Vert k [\vr^m]_\delta \vu \Vert_{L^q_t L^q_x} \leq  \Vert [\vr^m]_\delta\Vert_{ L^{\frac{m+1}{m}}_{t}L^\infty_x} \Vert \vu \Vert_{L^\infty_t L^{2}_x}\leq C, \quad \text{for some } q>1.
	\end{split}
\end{align}
The last estimate holds because of the space regularization $ [\cdot]_{\delta} $,  i.e. $ \vr^m \in L^{\frac{m+1}{m}}_{t,x} $ implies $   [\vr^m]_\delta \in { L^{\frac{m+1}{m}}_t L^\infty_x} $. 
Therefore, all the terms on the r.h.s. of the momentum equation are bounded in $  L^{p}_t W^{-1,p}_{x}$ for some $ p>1 $, henceforth $ \pt ([\vr]_\delta^\delta \vu) \in  L^{p}_t W^{-1,p}_x$. Since we already have $ \vr \vu \in L^{\infty}_t L^{\frac{4}{3}}_x $, we deduce that 
\begin{align*}
	&[\vr]_\delta^\delta \vve \in L^{\frac{4}{3}}_t L^{\frac{4}{3}}_x 
	\text{ with } \pt ([\vr]_\delta^\delta \vve) \in  L^{p}_t W^{-1,p}_x \text{ for } p>1.
\end{align*}

To estimate $I_{2,1}$ we use the interpolation Lemma \ref{interpolation-lemma} with 
\begin{align*}
	q=\frac{4}{3}, \ q^\prime=p>1, \ \bar{q}=\infty \text{ and } \bar{q}^\prime=1, 
\end{align*}  
we  conclude that there exists $ 0<\alpha<1 $ such that
\begin{align}\nonumber
	\begin{split}
		\vert I_{2,1} \vert&\leq  \left\vert \int_0^T\!\!\! \int_{\Td} (- \Delta_x)^{-1}\Div_x( \pt ([\vr]_\delta^\delta\vve)) \mathscr{W}_1 \dx \; \dt \right\vert \\
		& \leq \Vert [\vr]_\delta^\delta \vve \Vert_{L_{t,x}^{\frac{4}{3}}}^{\alpha} \Vert \pt ([\vr]_\delta^\delta \vve) \Vert_{L^{p}_tW^{-1,p}_{x}}^{1-\alpha}\Vert \mathscr{W}_1 \Vert_{L^{\infty}}^{1-\alpha} \Vert \partial_t \mathscr{W}_1 \Vert_{L^{1}_tW^{-1,1}_x}^{\alpha}\\
		&\leq C\left(k, \frac{1}{\sigma} \right) \eta^\vt \Vert K_h \Vert_{L^1_x}.
	\end{split}
\end{align}
A similar estimate for $ {I}_{2,2} $ is achieved by considering $\mathscr{W}_2(t,y)$ in place of $\mathscr{W}_1(t,x)$.

Now, we also check that by the same arguments as in the previous section we have
\begin{align*}
	&\iintO{ K_{h}(x-y) \Big(\calD_{2}(\vve(x))-\calD_{2}(\vve(y))\Big) W(t,x,y)}\\
	&=\iintO{ K_{h}(x-y)  \Big((- \Delta_x)^{-1}\dv_x(  \Div_x([\vr\vu]_\delta\otimes \vve))-(- \Delta_y)^{-1}\dv_y(  \Div_y([\vr\vu]_\delta\otimes \vve))\Big) W(t,x,y)}\\
	&\leq C \eta^\vt\|K_h\|_{L^1_x}.
\end{align*}
The estimate of the part corresponding to $ \mathcal{D}_3 $ is even simpler at  in this level. It suffices to combine the last estimate of \eqref{mom-est-M}  with \eqref{eta_teta}, to obtain the desired result.\par 

Summarising our estimates from above, the analogue of \eqref{St-S0-1} is now 
\begin{align}\nonumber
	\begin{split}
			& \iintO{\Ov K_h(x-y) \vert{\vr(x)-\vr(y)}\vert^\sigma (w(x)+w(y)) (T)}\\
		&\leq \iintO{\Ov K_h(x-y) \vert{\vr_0(x)-\vr_0(y)}\vert^\sigma w_0(x)}\\
		&\quad +C   \iintTO{  \Ov K_{h}(x-y)\vert{\vr(y)-\vr(x)}\vert^\sigma   w(x)}\\
		&\quad +C\iintTO{\Ov K_h(x-y) \lr{D_{|x-y|}\vu(x)-D_{|x-y|}\vu(y)} \vert{\vr(x)-\vr(y)}\vert^\sigma w(x)}\\
		&\quad+ C\left(k, \frac{1}{\sigma} \right)(h^{\beta}\eta^{-\alpha}+\eta^\vt) +C \sigma  .
	\end{split}
\end{align}
In the above inequality,  we cannot use the uniform $L^\infty$ bound for the density in order to treat the third term. Instead we write
\begin{align*}
	&\iintO{\Ov K_h(x-y) \lr{D_{|x-y|}\vu(x)-D_{|x-y|}\vu(y)} |\vr(y)-\vr(x)|^\sigma w(x)}\\
	&\leq \iintO{\Ov K_h(x-y) \vert \lr{D_{|x-y|}\vu(x)-D_{|x-y|}\vu(y)} \vert \lr{\vr(y)+ \vr(x) +\sigma} w(x)}\\
	&\leq \iintO{\Ov K_h(x-y) \vert \lr{D_{|x-y|}\vu(x)-D_{|x-y|}\vu(y)} \vert \vr(y) w(x)}\\
	&\;+ \iintO{\Ov K_h(x-y) \vert \lr{D_{|x-y|}\vu(x)-D_{|x-y|}\vu(y)} \vert  \vr(x)w(x)}\\
	&\;+ \sigma \iintO{\Ov K_h(x-y) \vert \lr{D_{|x-y|}\vu(x)-D_{|x-y|}\vu(y)} \vert    w(x)}= \sum_{i=1}^3 \mathcal{J}_i.
\end{align*}
For the term $ \mathcal{J}_1 $, we use the change of variables for $z=x-y$ to get  
\begin{align*}
	&\iintO{\Ov K_h(z) \vert \lr{D_{|z|}\vu(x)-D_{|z|}\vu(y)} \vert \vr(y) w(y+z)}\\
	& \leq C \int_{\T^d}{\Ov K_h(z) \|D_{|z|}\vu(\cdot+z)-D_{|z|}\vu(\cdot)\|_{L^2_y} \Vert \vr \Vert _{L^2_y}} \, {\rm d}z\,\\
	&\leq C |\log h|^{-1/2}\|\vu(t)\|_{W^{1,2}_x}\Vert \vr(t) \Vert _{L^2_y} .
\end{align*}
Since $ \vr \in L^2_{t,x} $ and $ \vu \in L^2_t W^{1,2}_{x} $, we have the following estimate
\begin{align*}
	&\iintTO{{\Ov K_h(x-y) \vert \lr{D_{|x-y|}\vu(x)-D_{|x-y|}\vu(y)} \vert \vr(y) w(x)}} \leq C |\log h|^{-1/2}
\end{align*}
Similarly, for $ \mathcal{J}_2 $ and $ \mathcal{J}_3 $, an analogous computation yields
\begin{align*}
	&\iintTO{{\Ov K_h(x-y) \vert \lr{D_{|x-y|}\vu(x)-D_{|x-y|}\vu(y)} \vert \vr(x) w(x)}} \leq C |\log h|^{-1/2}
\end{align*}
and 
\begin{align*}
	\sigma \iintO{\Ov K_h(x-y) \vert \lr{D_{|x-y|}\vu(x)-D_{|x-y|}\vu(y)} \vert    w(x)} \leq C \sigma  |\log h|^{-1/2}
\end{align*}
Therefore, combining all the above estimates we deduce
\begin{align}\nonumber
	\begin{split}
	&\iintTO{\Ov K_h(x-y) \lr{D_{|x-y|}\vu(x)-D_{|x-y|}\vu(y)} |\vr(y)-\vr(x)|^\sigma w(x)}\\
	&\leq C(1+\sigma)  |\log h|^{-1/2}.
\end{split}
\end{align}
Hence, bringing back the $\vr^M$ notation, we write the equivalent inequality of \eqref{jaktuwywalicM} in the following form
\begin{align}\label{M-lim}
	\begin{split}
			& \iintO{\Ov K_h(x-y) \vert{\vr^M(x)-\vr^M(y)}\vert^\sigma (w(x)+w(y)) (T)}\\
		&\leq \iintO{\Ov K_h(x-y) \vert{\vr_0(x)-\vr_0(y)}\vert^\sigma w_0(x)}\\
		&\quad +C   \iintTO{  \Ov K_{h}(x-y)\vert{\vr^M(y)-\vr^M(x)}\vert^\sigma   w(x)}\\
		&\quad +C(1+\sigma) |\log h|^{-1/2} + C\left(k, \frac{1}{\sigma} \right)(h^{\beta}\eta^{-\alpha}+\eta^\vt) +C \sigma  .
	\end{split}
\end{align}

We repeat a similar argument as done in Sections \ref{Ep:comp3}.  At first, we consider $\limsup_{M\rightarrow \infty}$ in \eqref{M-lim} and then we use a Gr\"onwall's inequality to get 
\begin{align}\label{M-lim-2}
	\begin{split}
			& \limsup_{M\rightarrow \infty} \iintO{\Ov K_h(x-y) \vert{\vr^M(x)-\vr^M(y)}\vert^\sigma (w(x)+w(y)) (T)}\\
		&\leq \left(C(1+\sigma) |\log h|^{-1/2} + C\left(k, \frac{1}{\sigma} \right)(h^{\beta}\eta^{-\alpha}+\eta^\vt) +C\sigma +o(h) \right) e^{CT} .
	\end{split}
\end{align}
Now, we choose $ \eta = h^{\frac{\beta}{2 \alpha}}$ and first perform the limit passage $ h \rightarrow 0$. After this step, the constant $C$ that involves $\frac{1}{\sigma}$ disappears. Therefore, the final limit passage $ \sigma \rightarrow 0 $ gives us the desired result. Finally, we perform the removal of weight as described in Section \ref{Ep:comp4}.

With strong convergence of the density at hand and the estimates \eqref{est:M} uniform in $M$, passage to the limit $M\to\infty$ is straightforward, provided that $m+1> \Gamma$.
After letting $M\to\infty$ our approximate system reduces to:
\begin{subnumcases}{\label{maine2Mboth}}
	 \pt\vr+\Div (\vr \vu) +k  \vr^m= 0,\label{maine1M}\\[10pt]
	\pt((\delta+[\vr]_\delta)\vu)+\Div([\vr \vu]_\delta\otimes\vu) - \Div \vS(\vu)+ \Grad p^\lambda(\vr) +k[\vr^m]_\delta \vu=\vc{0},\label{maine2M}
\end{subnumcases}
and it is satisfied in the sense of distributions along with the renormalized continuity equation.

\begin{rmk}
In the case of the general pressure law \eqref{gen-prho}, the term $A_3^{(i)}$ in \eqref{A3_sum:M} does not have the desired sign. Therefore, one needs the relation \eqref{gen-gamma} to be able to  estimate this term. For a more detailed discussion, see \cite[Section 4.2]{BJ2}.
\end{rmk}

\section{Compactness and passage to the limit $k\to 0$} \label{Sec:6} 
In this section we fix parameters $ \lambda \; $ and $\delta$ and consider $ k \rightarrow 0$. Our abbreviate notation is now $ \vr^{k}=\vr $ and $ \vu^{k}=\vu  $ for $k>0$. 

Note that estimate following from  \eqref{55} that we used before is not uniform in $k$, and therefore needs to be replaced. To do so, we derive an additional pressure estimate, called the Bogovskii estimate. 

\subsection{Uniform estimates}
In order to pass with $k\to0$ we need to obtain some uniform in $k$ estimates for the density.
Letting  $M\to\infty$ in the energy inequality \eqref{energy:M} and using lower semi-continuity of convex functions, for a.e. $\tau \in (0,T) $ we have:
\eq{\label{energy:k}
&\intO{\lr{(\delta +[\vr]_\delta) \frac{|\vu|^2}{2}+\vr e^{\lambda}(\vr)}(\tau)}
+\inttO{\vS(\vu):\Grad\vu}+ k\frac{\lambda\Gamma}{\Gamma-1}\inttO{\vr^{m+\Gamma-1}}\\
&\qquad+k\frac{\gamma}{\gamma-1}\inttO{  \vr^{m+\gamma-1}}\leq\intO{\lr{(\delta +[\vr_0^{\lambda}]_\delta) \frac{|\vu_0^{\lambda}|^2}{2}+\vr_0^{\lambda} e^{\lambda}(\vr_0^{\lambda})}},
}
where the r.h.s of the above inequality is uniformly bounded w.r.t. $ k$.
Therefore, we have the following estimates:
\begin{align}\label{est:k}
	\begin{split}
		&\|\vr\|_{L^\infty_t L^\Gamma_x}+ k^{\frac{1}{m+\Gamma-1}} \|\vr\|_{ L^{m+\Gamma+1}_{t,x}}\leq C,\\
		&\|\vu\|_{L^2_tW^{1,2}_x}+\delta^\frac{1}{2}\|\vu\|_{L^\infty_tL^2_x}\leq C,\\
		&\|[\vr]_\delta \vert \vu \vert^2\|_{L^\infty_tL^1_x}\leq C,\\
		& \|[\vr]_\delta \vu\|_{L^\infty_tL^\frac{2\Gamma}{\Gamma+1}_x} +   \|[\vr]_\delta \vu\|_{L^{m+1}_tL^\frac{2(m+1)}{m+3}_x}\leq C. 
	\end{split}
\end{align}
 The first estimate  allows us to prove that the limit of the extra term in the continuity equation $k\vr^m$ is zero, but it
 offers only an $L^1$ estimate for the pressure. The better estimate of the pressure will come from the so-called Bogovskii estimate, derived below.
 
 \subsection{The Bogovskii estimate}
 Since  $ \Gamma $ is large, we can test the momentum equation \eqref{maine2M} by a function:
$$\vcg{\varphi}=\psi\Grad\lap^{-1}[\vr], \quad\psi\in C^{\infty}_c[0,T).$$
After testing we obtain
\eq{\label{loc_pressure}
\intTO{\psi p(\vr)\vr}=&(\xi+\mu)\intTO{\psi\vr\Div\vu}\\
&+\intTO{\psi\mu\Grad\vu:\Grad^2\lap^{-1}[\vr]}\\
&-\intTO{\psi([\vr\vu]_\delta\otimes\vu):\Grad^2\lap^{-1}[\vr]}\\
&+\intTO{\psi \pt(([\vr]_\delta+\delta )\vu)\Grad\lap^{-1}[\vr]}\\
&+k\intTO{\psi[\vr^m]_\delta\vu\cdot\Grad\lap^{-1}[\vr]}
=\sum_{i=1}^5 I_i.
}
For $\Gamma\geq3$, the first three integrals can be controlled from the energy estimate, so we focus only on $I_4$ and $I_5$. To estimate $I_4$ we need to integrate by parts w.r.t. time. Then, the most restrictive term coming from this operation is equal to
\eq{\nonumber
&-\intTO{\psi ([\vr]_\delta+\delta)\vu\Grad\lap^{-1}[\pt\vr]}\\
&=\intTO{\psi([\vr]_\delta+\delta)\vu\Grad\lap^{-1}[\Div(\vr\vu)]}
+\intTO{\psi([\vr]_\delta+\delta)\vu\Grad\lap^{-1}[k\vr^m]}.
}
Again, the first term is standard for Navier--Stokes system, so we only focus on the second one, we have
\eq{\nonumber
\left|\intTO{\psi ([\vr]_\delta+\delta)\vu\Grad\lap^{-1}[ k\vr^m]}\right|
\leq C_{\psi}\|[\vr]_\delta\|_{L^\infty_{t,x}}\|\vu\|_{L^2_tL^6_x}\|\Grad\lap^{-1}[k\vr^m]\|_{L^\infty_tL^{\frac65}_x}\\
\leq C_{\psi}k\|[\vr]_\delta\|_{L^\infty_{t,x}}\|\vu\|_{L^2_tL^6_x}\|\vr^m\|_{L^\infty_tL^1_x},
}
which is bounded provided that $m\leq \Gamma$ as specified in \eqref{m, Gamma : M,m}. 
The estimate of $I_5$ from \eqref{loc_pressure} is easy under these conditions. 
After all, the  Bogovskii estimate provides another uniform in $k$ bound for the density:
\begin{align*}
	\|\vr\|_{L^{\Gamma+1}_{t,x}}\leq C.
\end{align*}
For more detailed proof of Bogovskii estimate for Navier-Stokes equations, see \cite[Chapter 5]{F}.

With this estimate, passage to the limit  $k\to0$ in the continuity and in the momentum equation boils down to justification of the compactness criterion for the sequence $\vr^k$ approximating the density.  The proof that the Kolmogorov compactness criterion works again at this level is presented in the subsection below.

\subsection{Compactness criterion}\label{Compact:k}

We will now discuss the main differences in the proof of the  main estimate from Section \ref{cc:M}. They are related only to the estimate of the term $I_2$ in the estimate of 
 $A^{(ii)}_3$, as the estimate for $I_1$ stays the same as in \eqref{est:M:I2}. 
 Even for the estimate of $I_2$ we only need to focus on one of the terms related to $\calD_1(\vv_\eta)$, i.e.
 \eqh{I_{2,1}= \int_0^T\!\!\! \int_{\Td} (- \Delta_x)^{-1}\Div_x( \pt ([\vr]_\delta^\delta\vve)) \mathscr{W}_1\, \dx \; \dt,
  }
  where
  \eqh{\mathscr{W}_1(t,x)= \int_{\T^d} K_h(x-y) W(t,x,y) \dy .}
 To control this term we first derive the estimate of $ W  $ and $ \pt W$ uniform w.r.t. $ k $. 
Recall the definition of $\vc{A},\, \vc{B}$, and $\mathcal{C}$ from \eqref{ABC:defn} and the expression 
\eqh{\pt W=\dv_x \vc{A} +\dv_y \vc{B} +\mathcal{C}.
}
First, we notice that 
\begin{align*}
	&\|\vc{A}\|_{L^2_tL^6_xL^\infty_y}+\|\vc{B}\|_{L^2_tL^6_yL^\infty_x}\leq C.
\end{align*}
Let us now inspect the term $\mathcal{C}$ which is given as 
\begin{align*}
    \mathcal{C}(t,x,y)= &\Big(- k \vr^{m-1}(x)  \sgn^\sigma_{xy}- \Lambda  \mathcal{M}[|\nabla \vu|]  \sgn^\sigma_{xy} - \vr(x)\dv_x\vu(x) \partial \sgn_{xy}^\sigma\\
	&+ (\sgn^\sigma_{xy} + \vr(y) \partial \sgn^\sigma_{xy}) \dv_y \vu(y) -k \partial \sgn_{xy}^\sigma (\vr^m(x)- \vr^m (y) )\Big) \vr(x) w(x).
\end{align*}
Here, we have the following estimates for the terms involving $ k $:
\begin{align*}
	&\Vert k \vr^{m-1} \sgn^\sigma_{xy} \vr(x) w(x) \Vert_{L^1_{t,x,y}} \leq C(\Gamma) k^{\frac{\Gamma}{m+\Gamma-1}} \Vert k^{\frac{1}{m+\Gamma-1}} \vr \Vert_{ L^{m+\Gamma-1}_{t,x}  }^{m-1} \Vert \vr w \Vert_{L^\infty_{t,x}},\\
	&\Vert k \partial \sgn_{xy}^\sigma (\vr^m(x)- \vr^m (y) )  \vr(x) w(x) \Vert_{L^1_{t,x,y}} \leq k \frac1\sigma \Vert \vr w \Vert_{L^\infty_{t,x}}.
\end{align*}
For the other terms we use the uniform estimates \eqref{est:k}, to obtain
\begin{align}\label{k-C}
\begin{split}
	\Vert \mathcal{C} \Vert_{L^1_tL^1_x L^\infty_y+ L^1_tL^1_yL^\infty_x}
	& \leq C\left( \frac{1}{\sigma} \right)\; \Vert \vr w \Vert_{L^\infty_{t,x}} \left(\Vert \vr \Vert_{L^2_t L^2_x} \Vert \nabla_x \vu  \Vert_{L^2_t L^2_x} + \Vert \vr \Vert_{L^2_t L^2_y} \Vert \nabla_x \vu  \Vert_{L^2_t L^2_y}  \right)\\
	&+ C\lr{\frac1\sigma} k^\beta \Vert \vr w \Vert_{L^\infty_{t,x}}\lr{1+\Vert k^{\frac{1}{m+\Gamma-1}} \vr \Vert_{L^{m+\Gamma-1}_{t,x}}^{m-1}+\Vert k^{\frac{1}{m+\Gamma-1}} \vr \Vert_{L^{m+\Gamma-1}_{t,y}}^{m-1}},
\end{split}
\end{align}
for some $ \beta >0$. 
Hence, for the term $\mathscr{W}_1$ defined in \eqref{M:w1} we have 
\begin{align*}
	&\Vert 	 \mathscr{W}_1  \Vert_{L^\infty_{t,x}} \leq \Vert K_h \Vert_{L^1_x} \Vert W  \Vert_{L^\infty_{t,x,y}} ,\\
 &\Vert 	\pt \mathscr{W}_1 \Vert_{ L^2_tW^{-1,6}_x +L^1_tW^{-1,1}_x} \leq \lr{C+ k^\beta}\Vert K_h \Vert_{L^1_x}.
\end{align*}
The time derivative of the momentum might be estimated from the momentum equation
\eq{\nonumber
	\pt( [\vr]_\delta^\delta \vu)=-\Div( [\vr \vu]_\delta \otimes\vu) + \Div \vS(\vu)- \Grad p^{\lambda}(\vr) -k [\vr^m]_\delta \vu.
} 
To estimate the terms on the r.h.s. of the above equation we use the uniform estimates \eqref{est:k} to get
\begin{align}\nonumber
	\begin{split}
		&\Vert [\vr \vu]_\delta \otimes\vu \Vert_{L^2_t L^{\frac{12}{11}}_x} \leq \Vert [\vr \vu]_\delta  \Vert_{L^\infty_x L^{\frac{4}{3}}_x} \Vert \vu \Vert_{L^2_t L^{6}_x}\leq C ,\\
		&\Vert \vS(\vu)\Vert_{L^2_t L^2_x} \leq C,\\
		&\Vert p^{\lambda}(\vr) \Vert_{L^q_t L^q_x} \leq \Vert \vr^\Gamma \Vert_{L^q_t L^q_x} \leq \Vert \vr \Vert_{L^{\Gamma+1}_{t,x}} \leq C, \text{ for } 1<q<\frac{\Gamma+1}{\Gamma},\\
		&\Vert k [\vr^m]_\delta \vu \Vert_{L^q_t L^q_x} \leq k^\beta  \Vert [\vr^m]_\delta\Vert_{ L^{\frac{m+\Gamma+1}{m}}_{t}L^\infty_x} \Vert \vu \Vert_{L^\infty_t L^{2}_x}\leq C, \ \text{for some } q>1,
	\end{split}
\end{align}
where $ \beta >0 $ is a generic constant. As a result we deduce that 
\begin{align*}
	[\vr]_\delta^\delta \vve \in L^{\frac{2\Gamma}{\gamma+1}}_t L^{\frac{2\Gamma}{\Gamma+1}}_x \quad
	\text{ with }\quad \pt ([\vr]^\delta_\delta \vve) \in  L^{p}_t W^{-1,p}_{x} \text{ for } p>1. 
\end{align*}
Hence using interpolation Lemma \eqref{interpolation-lemma}, with 
\begin{align*}
	q=\frac{2\Gamma}{\Gamma+1},\  q^\prime=p>1, \ \bar{q}=\infty \text{ and } \bar{q}^\prime=1, 
\end{align*}  
we conclude that there exists $ 0<\alpha<1 $ such that 
\begin{align}\nonumber
	\begin{split}
		\vert {I}_{2,1} \vert&\leq  \left\vert \int_0^T \!\!\!\int_{\Td} (- \Delta_x)^{-1}\Div_x( \pt ([\vr]_\delta^\delta\vve)) \mathscr{W}_1 \dx \; \dt \right\vert \\
		& \leq \Vert [\vr]_\delta^\delta \vve \Vert_{L_{t,x}^{q}}^{\alpha} \Vert \pt ([\vr]_\delta^\delta \vve) \Vert_{L^{p}_tW^{-1,p}_x}^{1-\alpha}\Vert \mathscr{W}_1 \Vert_{L^{\infty}_{t,x}}^{1-\alpha} \Vert \partial_t \mathscr{W}_1 \Vert_{L^{1}_tW^{-1,1}_x}^{\alpha}\\
		&\leq C(1+ k^\beta) \eta^\vt \Vert K_h \Vert_{L^1_x} .
	\end{split}
\end{align}
The estimate of $I_{2,2}$ is done the same way as in the previous step.

With this estimate at hand, the proof of compactness of the sequence approximating the density follows the same steps as in the previous section. We skip the details here, although it is worth mentioning that the constant $C\lr{k,\frac{1}{\sigma}}$ appearing in \eqref{M-lim-2} has a particular form $ (1+ k^\beta) C\lr{\frac{1}{\sigma}} $ that follows from \eqref{k-C}.

\subsection{Limit passage $k\to0$}
Clearly, passage to the limit $k\to0$ involves only a difficulty in passing to the limit in the pressure term and in the additional $k$-dependent terms in the continuity equation and momentum equation.
After compactness criterion is proven, and thanks to the Bogovskii estimate, identification of the weak limit in the pressure term is possible. For the term in the continuity equation, we note that
\eq{\nonumber
\left|\intTO{k\vr^m\phi}\right|=k^{\frac{\Gamma-1}{m+\Gamma-1}}\left|\intTO{\lr{k^{\frac{1}{m+\Gamma-1}}\vr}^m\phi}\right|\\
\leq Ck^{\frac{\Gamma-1}{m+\Gamma-1}}
\| k^{\frac{1}{m+\Gamma-1}}\vr\|^m_{L^{m+\Gamma-1}_{t,x}}\|\phi\|_{L^{\frac{m+\Gamma-1}{\Gamma-1}}_{t,x}}
\leq Ck^{\frac{\Gamma-1}{m+\Gamma-1}},}
and so it disappears in the limit $k\to0$ for any $\phi\in C^\infty_c((0,T)\times\Td)$. 

Next, for the term in the momentum equation, our first observation is
$  k^{\frac{1}{m+\Gamma-1}} \vr \in L^{m+\Gamma+1}_t $ implies $ k^{\frac{1}{m+\Gamma-1}} [\vr]_\delta \in L^{m+\Gamma+1}_t L^\infty_x $. Using the fact
\[ \Vert [\vr^m]_\delta \Vert_{L^\infty_x} \leq C(\delta)  \Vert \vr^m \Vert_{L^1_x} \]
along with a simple computation, we deduce
\begin{align*}
	\left|\intTO{k[\vr^m]_\delta \vu\cdot \pmb{\phi}}\right|\leq C(\delta) k^{\frac{\Gamma-1}{m+\Gamma-1}}  \Vert k^{\frac{1}{m+\Gamma-1}} \vr \Vert_{ L^{m+\Gamma-1}_{t,x}  }^m \Vert \vu \Vert_{L^2_t W^{1,2}_{x}},
\end{align*}
for any $\pmb{\phi}\in C^\infty_c((0,T)\times\Td)$. Thus this term also disappears in the limit $k\to 0$. 

Ultimately, our approximate system after letting $k\to 0$ is reduced to
\begin{subnumcases}{\label{maine1kboth}}
	 \pt\vr+\Div (\vr \vu) = 0,\label{maine1k}\\[10pt]
	\pt((\delta+[\vr]_\delta)\vu)+\Div([\vr \vu]_\delta\otimes\vu) - \Div \vS(\vu)+ \Grad p^\lambda(\vr) =\vc{0},\label{maine2k}
\end{subnumcases}
and it is satisfied in the sense of distributions, along with the renormalized continuity equation.

\section{Compactness and passage to the limit for $\delta \to 0$}\label{Sec:7}

Keeping $ \lambda $ fixed,  we proceed with the penultimate limit passage $\delta\to 0$. Once again we drop the indexes of the sequences denoting $ \vr^{\delta}=\vr $ and $ \vu^{\delta}=\vu  $, and we only use subindex $\delta$ to denote the regularisation in space.

\subsection{Uniform estimates}
After the limit passage $k\rightarrow 0 $, in \eqref{energy:k}, we derive energy inequality
\begin{align}\label{energy:del}
    \begin{split}
       &\intO{\!\lr{\frac{1}{2} (\delta+ [\vr]_{\delta} )\vert \vu \vert^2 + \vr e^\lambda(\vr)}(\tau)} 
	+\inttO{\vS(\vu):\Grad\vu} 
  \leq \!\intO{\lr{\frac{1}{2}\! (\delta+ [\vr_0^{\lambda}]_{\delta} )\vert \vu_0^{\lambda} \vert^2 + \vr e^\lambda(\vr_0^{\lambda})}\!\!}.
    \end{split}
\end{align}
The hypothesis \eqref{id-approx} implies 
\begin{align*}
    \intO{\lr{\frac{1}{2} (\delta+ [\vr_0^{\lambda}]_{\delta} )\vert \vu_0^{\lambda} \vert^2 + \vr e^\lambda(\vr_0^{\lambda})}} \leq (1+\delta)C,
\end{align*}
where $C$ is independent of $ \delta $.

We once again start from the momentum equation
and repeat the Bogovskii estimate
\eq{\label{del:bog}
	&\intTO{\psi p(\vr)\vr}\\
	&=(\lambda+\mu)\intTO{\psi\vr\Div\vu} +\intTO{\psi\mu\Grad\vu:\Grad^2\lap^{-1}[\vr]}\\
	&-\intTO{\psi([\vr\vu]_\delta\otimes\vu):\Grad^2\lap^{-1}[\vr]}+\intTO{\psi \pt([\vr]_\delta\vu)\Grad\lap^{-1}[\vr]}
	=\sum_{i=1}^4 I_i.
}
The only term that can be a problem now is $I_3$ because it does not allow us to use estimates for $|\vu|^2\in L^2_t L^3_x$ directly. However, passing to the limit $k\to 0$ in the energy inequality  yields
\eq{
	& \| [\vr]_{\delta}\vert \vu \vert^2\|_{L^\infty_t L^1_x}+ \Vert  \vr \Vert_{L^\infty_t L^\Gamma_x} + \| \vu \|_{L^2_t W^{1,2}_x} \ \dt + \sqrt{\delta}  \|  \vu  \|_{L^\infty_t L^2_x}
	\leq C, \label{del:est}
}
where $ C $ is independent of $ \delta $. From this, we deduce that
\eq{\label{66}
\|[\vr]_\delta\vu\|_{L^\infty(0,T; L^{\frac{2\Gamma}{\Gamma+1}}(\Td))}+\|\vr\vu\|_{L^2(0,T; L^{\frac{6\Gamma}{\Gamma+6}}(\Td))}\leq C.
}
Since $ \Gamma > m> 3 $, from \eqref{del:bog} we are able to obtain 
\begin{align*}
	\|\vr\|_{L^{\Gamma+1}_{t,x}}\leq C,
\end{align*}
with $C$ independent of $\delta$.

\subsection{Compactness criterion}
Again, we focus only on the elements of the main estimate that will change. In this limit passage, they are only related to the estimate of term $A_{3}^{(ii)}$ defined in \eqref{M:D1}, and more precisely to the estimate of $I_2$.
The first observation is that the part of $I_2$ corresponding to $ \mathcal{D}_3(\vv_\eta)$ defined in \eqref{def:D} is equal to zero in this case, but there will be a significant difference in the term corresponding to $ \mathcal{D}_2(\vv_\eta)$, as it now involves space regularisation of the momentum. 

For the part of $I_2$ corresponding to $ \mathcal{D}_1(\vv_\eta)$, the two terms $ I_{2,1} $ and $I_{2,2}$ will be treated similarly as  in previous Sections \ref{cc:M} and \ref{Compact:k}. 
Recall that
\eqh{\pt W=\dv_x \vc{A} +\dv_y \vc{B} +\mathcal{C}.
}
\eqh{
	&\vc{A}(t,x,y)= -\vr(x) w(x)\vu(x) \sgn^\sigma_{xy} ,\\
	&\vc{B}(t,x,y)= - \vu(y) \sgn_{xy}^\sigma  \vr(x)w(x) ,\\
	&\mathcal{C}(t,x,y)= \Big(- \Lambda  \mathcal{M}[|\nabla \vu|]  \sgn^\sigma_{xy} - \vr(x)\dv_x\vu(x) \partial \sgn_{xy}^\sigma  (\sgn^\sigma_{xy} + \vr(y) \partial \sgn^\sigma_{xy}) \dv_y \vu(y) \Big) \vr(x) w(x),\label{del:ABC:defn}}
with
\begin{align*}
	\|\vc{A}\|_{L^2_tL^6_xL^\infty_y}+\|\vc{B}\|_{L^2_tL^6_yL^\infty_x}\leq C,
\end{align*}
and
\begin{align*}
    	\Vert \mathcal{C} \Vert_{L^1_t L^1_x L^\infty_y+ L^1_t L^1_y L^\infty_x} \leq C\left( \frac{1}{\sigma},\Lambda \right)\; \Vert \vr w \Vert_{L^\infty} \left(\Vert \vr \Vert_{L^2_t L^2_x} \Vert \nabla \vu  \Vert_{L^2_t L^2_x} + \Vert \vr \Vert_{L^2_t L^2_y} \Vert \nabla \vu  \Vert_{L^2_t L^2_y} + \Vert \nabla_x \vu  \Vert_{L^2_t L^2_x}\right).
\end{align*}
For the estimate of $ \pt W $, the important thing is to get a uniform bound of $ \vr $ in $ L^2_t L^2_x $ and the choice of $ \Gamma $ is providing it. 
Using the momentum equation to express the time derivative of the momentum, 
along with the $\delta$ independent estimates from \eqref{del:est} and \eqref{66},
we eventually deduce that 
\begin{align*}
	[\vr]_\delta^\delta \vve \in L^{\frac{2\Gamma}{\gamma+1}}_t L^{\frac{2\Gamma}{\Gamma+1}}_x 
	\text{ with } \pt ([\vr]_\delta^\delta \vve) \in  L^{p}_t W^{-1,p}_{x} \text{ for } p>1. 
\end{align*}
The estimates of the term $ \mathscr{W}_1 $ are the same as in the two previous sections.
Therefore, using the interpolation Lemma \eqref{interpolation-lemma} in the same way as in Section \ref{Compact:k} allows us to conclude
\begin{align}\nonumber
	\begin{split}
		\vert {I}_{2,1} \vert\leq C\eta^\vt \Vert K_h \Vert_{L^1_x} .
	\end{split}
\end{align}

To complete the estimate of $I_2$, we now focus on the part corresponding to $\calD_{2}(\vv_\eta)$, we further split it into 
\begin{align*}
	&\iintTO{ K_{h}(x-y) \Big(\calD_{2}(\vve(x))-\calD_{2}(\vve(y))\Big) W(t,x,y)}\\
	&=\iintTO{ K_{h}(x-y)  (- \Delta_x)^{-1}\dv_x(  \Div_x([\vr\vu]_\delta\otimes \vve)) W(t,x,y)}\\
	&-\iintTO{ K_{h}(x-y)  (- \Delta_y)^{-1}\dv_y(  \Div_y([\vr\vu]_\delta\otimes \vve)) W(t,x,y)} \\
	&= \tilde I_{2,1}+ \tilde I_{2,2}.
\end{align*}
To deal with the first term, we first rewrite
$$  [\vr\vu]_\delta\otimes \vve= ([\vr\vu]_\delta-[\vr]_\delta \vu) \otimes \vve + [\vr]_\delta \vu \otimes \vve. $$
Next we note that 
\begin{align*}
	\vert  \tilde I_{2,1}  \vert & \leq \left\vert \int_0^T\!\!\! \int_{\Td}  (- \Delta_x)^{-1}\dv_x(  \Div_x([\vr\vu]_\delta\otimes \vve)) \mathscr{W}_1(t,x) \dx \dt \right\vert \\
	& \leq \Vert  \mathscr{W}_1 \Vert_{L^\infty_{t,x} } \lr{\Vert ([\vr\vu]_\delta-[\vr]_\delta \vu) \otimes \vve \Vert_{L^1_{t,x}}+ \Vert [\vr]_\delta \vu \otimes \vve \Vert_{L^1_{t,x}}}\\
	& \leq \Vert  \mathscr{W}_1 \Vert_{L^\infty_{t,x} } \lr{\Vert ([\vr\vu]_\delta-[\vr]_\delta \vu)\Vert_{L^2_t L^{\frac{6}{5}}_x} \Vert \vve \Vert_{L^2_t L^6_x }+ \Vert [\vr]_\delta \vve \Vert_{L^2_t L^{\frac{6}{5}}_x} \Vert \vu \Vert_{L^2_t L^6_x}}.
\end{align*}
Now we note that in order to estimate $ L^2_t L^{\frac{6}{5}}_x$ norm of $[\vr]_\delta \vve $, we use the interpolation  between the spaces  $ L^1_t L^1_x$ and  $ L^{\infty}_t L^{\frac{2\Gamma}{\Gamma+1}}_x $, where its norm is bounded due to assumption since $ \frac{2\Gamma}{\Gamma+1} > \frac{6}{5}$. Therefore, there exists $ \alpha >0 $ such that
\begin{align*}
    \vert \tilde I_{2,1}   \vert & \leq \Vert  \mathscr{W}_1 \Vert_{L^\infty_{t,x} } \bigg(\Vert ([\vr\vu]_\delta-[\vr]_\delta \vu)\Vert_{L^1_t L^1_x}^\alpha \Vert ([\vr\vu]_\delta-[\vr]_\delta \vu)\Vert_{L^\infty_t L^{\frac{2\Gamma}{\Gamma+1}}_x}^{1-\alpha} \Vert \vve \Vert_{L^2_t L^6_x } \\
    &\hspace{4cm}+ \Vert [\vr]_\delta \vve \Vert_{L^1_t L^1_x}^\alpha \Vert [\vr]_\delta \vve \Vert_{L^{\infty}_t L^{\frac{2\Gamma}{\Gamma+1}}_x}^{1-\alpha} \Vert \vu \Vert_{L^2_t L^6_x}\bigg)\\
    &\leq \Vert K_h \Vert_{L^1_x} \Vert W  \Vert_{L^\infty_{t,x,y}}  \lr{\vert \log \delta \vert^{-\alpha} + \eta^{\vt } }.
\end{align*}
The estimate of $ \tilde I_{2,2} $ is completely analogous. Having term $A_3^{(ii)}$ controlled in terms of small parameters $\delta$ and $\eta$ we can essentially repeat the proof of compactness criterion from Sections \ref{Ep:comp3} and \ref{Ep:comp4}.

\subsection{Limit passage}
The compactness of $ \vr $ allows for the identification of weak limit of the pressure term, while the convergence of the other terms is quite straightforward, one can follow Feireisl \cite[Chapter 7]{F}.
The limiting system is satisfied in the sense of distributions:
\begin{subnumcases}{\label{main-lim-del}}
	 \pt\vr+\Div (\vr \vu) = 0,\label{main-lim-del1}\\[10pt]
	\pt(\vr\vu)+\Div(\vr \vu\otimes\vu) - \Div \vS(\vu)+ \Grad p^\lambda(\vr) =\vc{0}.\label{main-lim-del2}
\end{subnumcases}

\section{Compactness and limit passage for $\lambda\to0$ and $ \gamma \geq  \frac{9}{5} $}\label{Sec:8} 

In this section, we consider the last limit passage $\lambda\to 0$. The sequences of solutions are again denoted by  $ \vr^{\lambda}=\vr $ and $ \vu^{\lambda}=\vu  $ for $\lambda>0$.
\subsection{Uniform estimates}
Letting $ \delta\rightarrow 0 $ in the energy inequality \eqref{energy:del}, and using lower-semicontinuity of convex functions, we obtain
\eq{\label{energy:lam}
	&\intO{\lr{\frac{1}{2} \vr \vert \vu \vert^2 + \vr e^\lambda(\vr)}(\tau)} 
	+\inttO{\vS(\vu):\Grad\vu} \leq \intO{\lr{\frac{1}{2} \vr_0^\lambda \vert \vu_0^\lambda \vert^2 + \vr_0^\lambda e^\lambda(\vr_0^\lambda)}}
}
for a.e. $ \tau \in (0,T) $. The uniform bound in the r.h.s of \eqref{energy:lam} follows from the convergence of initial data as described in \eqref{id-approx-conv}.
 
Therefore, we have the following estimates are uniform w.r.t. $\lambda$:
\begin{align}\label{est:lambda}
	\begin{split}
		&\|\vr\|_{L^\infty_t L^\gamma_x}+ \lambda^{\frac{1}{\Gamma}} \|\vr\|_{ \in  L^\infty_tL^\Gamma_x}\leq C,\\
		&\|\vu\|_{L^2_tW^{1,2}_x}+\|\vr \vert \vu \vert^2\|_{L^\infty_tL^1_x}\leq C,\\
		& \|\vr \vu\|_{L^\infty_t L^\frac{2\gamma}{\gamma+1}_x} + \|\vr \vu\|_{L^2_t L^\frac{6\gamma}{\gamma+6}_x} \leq C. 
	\end{split}
\end{align}
Performing Bogovskii estimate for $ \gamma\geq \frac{9}{5} $, we again have that
\begin{align}\label{lam:bog}
	\|\vr\|_{L^{r}_{t,x}}\leq C\quad \text{for some } \; r\geq 2,
\end{align}
see Feireisl \cite[Chapter 5]{F} or Novotn\'y and Stra\v skraba \cite[Chapter 3.3]{NS} for a detailed proof.

\subsection{Compactness criterion and the limit passage $\lambda\to 0$}

Here the only problem is the justification of why we can apply the interpolation  Lemma \ref{interpolation-lemma} to estimate the terms $I_{2,1}$ and $I_{2,2}$. All other contributions of $A_3^{(ii)}$ stay the same.
The estimate \eqref{lam:bog} along with the uniform estimates of the momentum\eqref{est:lambda} implies that
$$   \vr \vu  \in L_t^2 L_x^{\frac{6\gamma}{\gamma+1}} \cap  L_t^\infty L_x^{\frac{2\gamma}{\gamma+1}} \cap  L_t^{1+0} L_x^{r} $$ 
with $r \geq \frac{3}{2}$. As a consequence of the interpolation of Sobolev spaces, it holds 
\[L_t^2 L_x^{\frac{6\gamma}{\gamma+1}} \cap L_t^\infty L_x^{\frac{2\gamma}{\gamma+1}} \cap  L_t^{1+0} L_x^{r} \hookrightarrow L_t^{q} L_x^{q} \text{ with } q\geq \frac{10}{7}\]
and therefore, $ \vr \vu $ is uniformly bounded in $L_t^{q} L_x^{q} \text{ with } q\geq \frac{10}{7}$  with the corresponding norms uniformly bounded w.r.t. parameter $\lambda$.
 Moreover, from the momentum equation we can deduce that  $  \pt (\vr \vu)\in L^{q^\prime}_t W^{-1,q^\prime}_x$ with $q^\prime \geq \frac{10}{9}$.
We combine the above discussion to get
\begin{align*}
	\vr \vve \in L_t^{q} L_x^{q} \quad 
	\text{ with } \pt (\vr \vve) \in  L^{q^\prime}_t W^{-1,q^\prime}_x \text{ for } q\geq \frac{10}{7} \text{ and }q^\prime \geq \frac{10}{9}. 
\end{align*}
%
On the other hand, the estimates for term $ \mathscr{W}_1$ are the same as in the previous sections.
Therefore, we can estimate $I_{2,1}$ applying Lemma \eqref{interpolation-lemma} for the choice of parameters
\begin{align*}
	q=\frac{10}{7},\  q^\prime= \frac{10}{9}, \ \bar{q}=\infty \text{ and } \bar{q}^\prime=1 
\end{align*} 
that satisfies \eqref{il:a1}. Hence, there exists a $ 0<\alpha <1 $ such that the following holds
\begin{align}\nonumber
	\begin{split}
		\vert {I}_{2,1} \vert&\leq  \left\vert \int_0^T \!\!\!\int_{\Td} (- \Delta_x)^{-1}\Div_x( \pt (\vr\vve)) \mathscr{W}_1 \dx \; \dt \right\vert \\
		& \leq \Vert \vr \vve \Vert_{L_{t,x}^{q}}^{\alpha} \Vert \pt (\vr \vve) \Vert_{L^{q^\prime}_t W^{-1,q^\prime}_x}^{1-\alpha}\Vert \mathscr{W}_1 \Vert_{L^{\infty}_{t,x}}^{1-\alpha} \Vert \partial_t \mathscr{W}_1 \Vert_{L^{1}_tW^{-1,1}_x}^{\alpha}\\
		&\leq C \eta^\vt \Vert K_h \Vert_{L^1_x} .
	\end{split}
\end{align}
 This allows us to obtain the desired estimate for the term $A_3$. Therefore, we prove the compactness criterion.

The limit passage is quite standard since we already deduce the compactness of $ \vr^\lambda$.

\subsection{Unbounded initial density}\label{gen-id}
So far, our proof relied on  Proposition \ref{pro:rho-bdd}, where we  assumed the initial density to be bounded, i.e. $ 0\leq \vr_0 \in L^\infty(\T^d)$.  Let us now explain how to relax this assumption 
to more natural energy bound:
\begin{align}\label{e_bd_rho}
    0\leq \vr_0 \in L^\gamma(\T^d).
\end{align}
Let $M>>1$ be a fixed constant. We consider 
\eq{\nonumber
w_{0,M}:=\begin{cases}
1 &\mbox{ for } \vr_0\leq M,\\
\frac{M}{\vr_0} &\mbox{ for } \vr_0> M,
\end{cases}
}
clearly, by this definition
\begin{align}\nonumber
    \vr_0 w_{0,M} \leq M.
\end{align}
We denote by $w_M$ a solution to the equation \eqref{def_w}
with initial data $ w_{0,M}$.
Repeating the proof of Proposition \ref{pro:rho-bdd}, we obtain 
\begin{align}\nonumber
   \Vert \vr w_M \Vert_{L^\infty_{t,x} } \leq \Vert \vr_0 w_{0,M} \Vert_{L^\infty_x}\leq M.
\end{align}
Now, we repeat the proof of compactness with $w_M$ in place of $w$. The only difficulty we have is to perform the step where we remove the weights, as explained in Section \ref{Ep:comp4}. To do this we needed the uniform bound of $ \vr \vert \log w \vert $ in $L^\infty_t L^1_x$ as described in \eqref{boundlogw}. Starting  from \eqref{eq:rhologw} without the $k$-dependent term we deduce 
\begin{equation*}
\underset{t>0}{{\rm ess}\sup}\intO{ \vr |\log w_M| } \leq \intO{ \vr_0 |\log w_{0,M}| } + C  \Lambda,
\end{equation*}
where $C$ is independent of $M$. Now, from the definition of $ w_{0,M} $, we obtain 
\begin{align*}
    \intO{ \vr_0 |\log w_{0,M}| }\leq \int_{\{\vr>M\}} \vr_0 |\log w_{0,M}| \dx  \leq \int_{\{\vr>M\}} \vr_0 |\log \vr_{0}| \dx \leq C_\gamma  \int \vr_0^\gamma \dx.
\end{align*}
Therefore, we conclude that 
\begin{align}\label{uni-wM}
   \underset{t>0}{{\rm ess}\sup}\intO{ \vr |\log w_M| } \leq C (\vr_0,\Lambda),
\end{align}
where $C$ is independent of $M$. Going back to the removal of weight as explained in Section \ref{Ep:comp4}, we observe that the above bound \eqref{uni-wM} is enough to deduce the inequality \eqref{Remove-w-I2}. Therefore,  compactness of the sequence approximating the density holds for general initial data \eqref{e_bd_rho}.

\appendix

\section{Appendix}\label{appendix} 

\subsection{Proof of Lemma \ref{Lemma:ConvK}}\label{PL2}

First note that from \eqref{def:Kh} we have that
\begin{equation}\nonumber
	\Ov K_h(z) \leq |\ln h|^{-1}\left(2^{d} + \frac{1}{(|z|+h)^d}\right) \leq C \Ov K_h(z).
\end{equation}
To show \eqref{out2} we first write
\eq{\nonumber
	&\int_{\Td} \Ov K_{h_1}(x-y) \Ov K_{h_2}(y-z)\, \dy \\
	&\leq  |\ln h_1|^{-1} |\ln h_2|^{-1} \int_{\Td} \left(2^{d} + \frac{1}{(|x-y|+h_1)^d}\right) \left(2^{d}+\frac{1}{
		(|y-z|+h_2)^d} \right)\,\dy\\
	& \leq
	|\ln h_1|^{-1} |\ln h_2|^{-1} \int_{\Td} \left(2^{d} + \frac{1}{(|(x-z)-y|+h_1)^d}\right)\left(2^{d}+\frac{1}{  (|y|+h_2)^d} \right)\,\dy=I.
}
Now we split $I=I_1+I_2+I_3$ by considering $\Td=\Omega_1+\Omega_2+\Omega_3$, where $\Omega_1=\{y:|y|\geq 2|x-z|\}$, $\Omega_2=\{y:|y|\leq \frac12|x-z|\}$, $\Omega_3=\{y: \frac12|x-z|\leq |y|\leq 2|x-z|\}$, respectively. 

Since on the set $\Omega_1$ we have $|(x-z)-y|\geq |y|-|x-z|\geq |x-z|$ we can estimate
\eq{\nonumber
	I_1&= |\ln h_1|^{-1} |\ln h_2|^{-1} \int_{\Omega_1} \left(2^{d} + \frac{1}{(|(x-z)-y|+h_1)^d}\right)\left(2^{d}+\frac{1}{  (|y|+h_2)^d} \right)\,\dy\\
	&\leq |\ln h_1|^{-1} |\ln h_2|^{-1} \int_{\Omega_1} \left(2^{d} + \frac{1}{(|x-z|+h_1)^d}\right) \left(2^{d}+\frac{1}{  (|y|+h_2)^d} \right)\, \dy \\
	&\leq |\ln h_1|^{-1} |\ln h_2|^{-1} \int_{\Td} \left(2^{d} + \frac{1}{(|x-z|+h_1)^d}\right) \left(2^{d}+\frac{1}{  (|y|+h_2)^d} \right)\, \dy \\
	&\leq 2\Ov K_{h_1}(x-z)
}
Similarly, on the set $\Omega_2$ we have  $|(x-z)-y|\geq |x-z|-|y|\geq |x-z|-|y|\geq  \frac12 |x-z|$, and so
\eq{\nonumber
	I_2&= |\ln h_1|^{-1} |\ln h_2|^{-1} \int_{\Omega_2} \left(2^{d} + \frac{1}{(|(x-z)-y|+h_1)^d}\right)\left(2^{d}+\frac{1}{  (|y|+h_2)^d} \right)\,\dy\\
	&\leq |\ln h_1|^{-1} |\ln h_2|^{-1} \int_{\Omega_1} \left(2^{d} + \frac{1}{(\frac{1}{2}|x-z|+h_1)^d}\right)\left(2^{d}+\frac{1}{  (|y|+h_2)^d} \right)\, \dy \\
	&\leq |\ln h_1|^{-1} |\ln h_2|^{-1} \int_{\Omega_1} \left(2^{d} + \frac{2^d}{(|x-z|+h_1)^d}\right)\left(2^{d}+\frac{1}{  (|y|+h_2)^d} \right)\, \dy \\
	& \leq 2^{d+1} \Ov K_{h_1}(x-z).
}
Finally, for the set $\Omega_3$ we have in particular that  $|y|\geq  \frac12|x-z|$
\eq{\nonumber
	I_3&= |\ln h_1|^{-1} |\ln h_2|^{-1} \int_{\Omega_3} \left(2^{d} + \frac{1}{(|(x-z)-y|+h_1)^d}\right)\left(2^{d}+\frac{1}{  (|y|+h_2)^d} \right)\,\dy\\
	&\leq  |\ln h_1|^{-1} |\ln h_2|^{-1} \int_{\T^d }\left(2^{d} + \frac{1}{(|(x-z)-y|+h_1)^d}\right)\left(2^{d}+\frac{2^d}{  (|x-z|+h_2)^d} \right)\,\dy\\
	&\{ \mbox{changing } y-(x-z) \to y\} \\
	& \leq |\ln h_1|^{-1} |\ln h_2|^{-1} \int_{\T^d} \left(2^{d} + \frac{1}{(|y|+h_1)^d}\right)\left(2^{d}+\frac{2^d}{  (|x-z|+h_2)^d} \right)dy\\
	&\leq 2^{d+1} \Ov K_{h_2}(x-z).
}
Summing up $I_1+I_2+I_3$ we verify \eqref{out2}. $\Box$

\subsection{Proof of Lemma \ref{Lagrange}}\label{Lemma11}
 In what follows we treat the function defined over the torus as a periodic function and $x,y \in [-\pi,\pi]^d$. This way we avoid questions with definition of function $u$. We take $u \in C^1$, and then conclude the result for Lebesgue points of $u$. We have
\eq{\label{a2:1}
u(x)-u(y)=u(x)-u(z)+u(z)-u(y)=
\int_{l(x,z)} \partial_{(x,z)}u(s)dl_s + \int_{l(z,y)}\partial_{(z,y)}u(s) dl_s,
}
where $l(x,z)$ is the line connecting $x$ to $z$, $\partial_{(x,z)}$
is the dirctional derivative in  direction of $\vec{[x,z]}$.

Next we integrate \eqref{a2:1} over $B(\frac{x+y}{2},\frac{|x-y|}{2})$ w.r.t. $z$. Note that 
$$B\lr{\frac{x+y}{2},\frac{|x-y|}{2}}\subset B(x,|x-y|)\cap  B(y,|x-y|)$$ and thus we obtain
\eq{\nonumber
   & \omega_d 2^{-d} |x-y|^d |u(x)-u(y)| \\
    &\leq\int_{B(x,|x-y|)} {\rm d}z \int_{l(x,z)} |\nabla u(s)| {\rm d}l_s +
    \int_{B(y,|x-y|)} {\rm d}z \int_{l(z,y)} |\nabla u(s)| {\rm d}l_s \\
   & \leq \int_0^{|x-y|} r^{d-1}dr \int_{S^{d-1}} d\omega \int_0^{|x-y|} d\tau |\nabla u(x+\omega \tau)|
    + \int_0^{|x-y|} r^{d-1}dr \int_{S^{d-1}} d\omega \int_0^{|x-y|} d\tau |\nabla u(y+\omega \tau)| 
    \\
   & \leq \int_0^{|x-y|} r^{d-1}dr \int_{S^{d-1}} d\omega \int_0^r d\tau |\nabla u(x+\omega \tau)| + \int_0^{|x-y|} r^{d-1}dr \int_{S^{d-1}} d\omega \int_0^r d\tau |\nabla u(y+\omega \tau)| 
    \\
       & \leq \int_0^{|x-y|} r^{d-1}dr \int_{S^{d-1}} d\omega \int_0^{|x-y|} d\tau |\nabla u(x+\omega \tau)|+ \int_0^{|x-y|} r^{d-1}dr \int_{S^{d-1}} d\omega \int_0^{|x-y|} d\tau |\nabla u(y+\omega \tau)|
       \\
       & \leq c|x-y|^d \int_{S^{d-1}} d\omega \int_0^{|x-y|} d\tau |\nabla u(x+\omega \tau)|\frac{\tau^{d-1}}{\tau^{d-1}} +  c|x-y|^d \int_{S^{d-1}} d\omega \int_0^{|x-y|} d\tau |\nabla u(y+\omega \tau)| \frac{\tau^{d-1}}{\tau^{d-1}}
    \\
    &\leq
    c|x-y|^d \left(
    \int_{B(0,|x-y|)}\frac{|\nabla u(x+z)|}{|z|^{d-1}} {\rm d}z +
    \int_{B(0,|x-y|)}\frac{|\nabla u(y+z)|}{|z|^{d-1}} {\rm d}z \right),
}
where the last inequality follows by substituting $z=\tau\omega$ and by noticing that $|z|=\tau$ and $\omega_d$ stands for the volume of unit sphere in $\R^d$. $\Box$

\subsection{Proof of Lemma \ref{lem:approx}} \label{Lemma12}
Our first goal is to control the time regularisation
\begin{equation}\nonumber
	\intTOD{ \vr(t,x) (\vu(t,x) -\pi_{\eta_t} \ast \vu(t,x))^2 },
\end{equation}
where $\pi_{\eta_t}$ denotes the time mollifier. To proceed, we rewrite this term as follows
\eq{\nonumber
	&\intTOD{\vr(t,x) (\vu(t,x) -\pi_{\eta_t}\ast \vu(t,x))^2}\\
	&=\intTOD{\vr(t,x) \lr{\vu(t,x) -\int_0^T\pi_{\eta_t}(t-s) \vu(s,x)\, ds}^2 }\\
	&\leq \intTOD{ \int_0^T\vr(t,x) (\vu(t,x) - \vu(s,x))^2 \pi_{\eta_t} (t-s)\, ds},}
where we used the Jensen's inequality.
Thus we have
\eq{\nonumber
	&\intTOD{\vr(t,x) (\vu(t,x) -\pi_\etat\ast \vu(t,x))^2}\\
	& {\leq}\intTOD{\int_0^T\vr(t,x) (\vu(t,x) - \vu(s,x)) \pi_\etat (t-s) (\vu(t,x) - \vu(s,x))\, ds}\\
	&=\intTOD{ \int_0^T (\vr(t,x) \vu(t,x) - \vr(s,x) \vu(s,x)) \pi_\etat (t-s) (\vu(t,x) - \vu(s,x)) \, ds}  \\
	&\quad {-} \intTO{\int_0^T
		(\vr(t,x) -\vr(s,x)) \vu(s,x) \pi_\etat(t-s) ((\vu(t,x) - \vu(s,x))\, ds}\\
  &=J_1 +J_2.
}
To estimate $J_1$ we first note that
$$
\pt(\vr \vu)\in L^{q}_t W^{-1,q}_x \text{ for } q>1,
$$
where $q$ being restricted by the integrability of the pressure mostly, and $q\geq \frac{10}{9}$ when $ \gamma \geq \frac{9}{5}$. Therefore, we introduce additional regularisation in space for the term $\vu(t,x)-\vu(s,x)$ in $J_1$ and $J_2$, denoted by $\etax$. From the compact embedding of $W^{1,2}_x \subset L^q_x \text{ with } q<6$, it is clear that there exists $\theta>0$ such that
\begin{equation}\nonumber
	\| \vu - \pi_\etax \ast \vu\|_{L^2_t L^q_x}\leq C \etax^\theta \ \mbox{ for all $q< 6$}.
\end{equation}
Moreover, we also have
\eq{\nonumber
	\|\pi_\etax\ast\lr{\vu(t,x)-\vu(s,x)}\|_{W^{1, r}_x}\leq \etax^{-\alpha}\ \mbox{ for some $\alpha>0$ and $ r<\infty$}.
}
Therefore, $J_1$ can be estimated as follows
\eq{\nonumber
	&\intTOD{ \int_0^T (\vr(t,x) \vu(t,x) - \vr(s,x) \vu(s,x)) \pi_\etat (t-s) (\vu(t,x) - \vu(s,x)) \, {\rm d}s}  \\
	&=\intTOD{ \int_0^T (\vr(t,x) \vu(t,x) - \vr(s,x) \vu(s,x)) \\
	& \hspace{3cm} \pi_\etat (t-s) [\vu(t,x) -\pi_\etax\ast\vu(t,x)-( \vu(s,x)- \pi_\etax\ast\vu(s,x))] \, {\rm d}s} \\
	&+\intTOD{ \int_0^T (\vr(t,x) \vu(t,x) - \vr(s,x) \vu(s,x)) \pi_\etat (t-s) \pi_\etax\ast(\vu(t,x) - \vu(s,x)) \, {\rm d}s}  \\
	&\leq C\etax^\theta
	+\int_0^T\!\!\!\int_0^T\pi_\etat(t-s)\intO{\!\!\lr{\int_t^s\partial_\tau (\vr\vu)\, d\tau\,}  \pi_\etax\ast(\vu(t,x) - \vu(s,x)) \,} \,{\rm d}s\,\dt\\
	&\leq C\etax^\theta
	+\int_0^T\!\!\!\int_0^T\pi_\etat(t-s)\lr{\int_t^s\|\partial_\tau (\vr\vu)\|_{W^{-1,q}_x}\, d\tau }  \|\pi_\etax\ast(\vu(t,x) - \vu(s,x))\|_{W^{1,r}_x} \, \,{\rm d}s\,\dt\\
	&\leq C\etax^\theta\\
	&+
	\lr{ \int_0^T\!\!\!\int_0^T\pi_\etat(t-s)\lr{\int_t^s\|\partial_\tau(\vr\vu)\|_{W^{-1,q} }  d\tau}^2 \,{\rm d}s\,\dt}^{\frac12}\\
	&\quad \quad \lr{ \int_0^T\!\!\!\int_0^T\pi_\etat(t-s)\|\pi_\etax\ast(\vu(t,x) - \vu(s,x))\|^2_{W^{1,r}} \, \,{\rm d}s\,\dt}^{\frac12}\\
	&\leq  C\etax^\theta+C\etax^{-\alpha}\etat^{\beta}.
}

Now, we want to estimate the term $J_2$ in a similar way. Since $\vr \in L^\infty_t L^{\gamma}_x$ and $\vu \in L^2_t W^{1,2}_x$, so the continuity equation gives
	$\pt\vr \in L^2_t W^{-1, \frac{6\gamma}{\gamma+6}}_x$,
with $\frac{6\gamma}{\gamma+6} \geq \frac{18}{13}$ when $\gamma \geq  \frac{9}{5}$. Again using the identity
\begin{align*}
    \vu(t,x) - \vu(s,x)= (\vu(t,x) -\pi_\etax\ast\vu(t,x)-( \vu(s,x)- \pi_\etax\ast\vu(s,x)) ) + \pi_\etax\ast(\vu(t,x) - \vu(s,x))
\end{align*} 
for the term $J_2$, we obtain
\eq{\nonumber
	& \intTOD{\int_0^T
		(\vr(t,x) -\vr(s,x)) \vu(s,x) \pi_\etat(t-s) (\vu(t,x) - \vu(s,x))\, {\rm d}s}\\
	&  \leq C\etax^\theta+\intTOD{\int_0^T
		(\vr(t,x) -\vr(s,x))\pi_\etax \ast\vu(s,x) \pi_\etat(t-s) (\vu(t,x) - \vu(s,x))\, {\rm d}s}\\
	&  \leq C\etax^\theta+\intTOD{\int_0^T
		(\vr(t,x) -\vr(s,x))\pi_\etax \ast\vu(s,x) \pi_\etat(t-s) \pi_\etax \ast(\vu(t,x) - \vu(s,x))\, {\rm d}s}\\
	& \leq C\etax^\theta+ \intTOD{\int_0^T\lr{\int_t^s\partial_\tau\vr\ d\tau}
		\pi_\etax \ast\vu(s,x) \pi_\etat(t-s) \pi_\etax \ast(\vu(t,x) - \vu(s,x))\, {\rm d}s}
}
\eq{\nonumber
	& \leq C\etax^\theta	\\
		&\quad + \intT{\!\!\!\int_0^T\pi_\etat(t-s)\lr{\int_t^s\|\partial_\tau\vr\|_{W^{-1,\frac{18}{13}}_x}\ d\tau}
	\|\pi_\etax \ast\vu(s,x)  \pi_\etax \ast(\vu(t,x) - \vu(s,x))\|_{W^{1,\frac{18}{5}}_x}\, {\rm d}s}	\\
	& \leq C\etax^\theta\\
		&\quad +\|\pt\vr\|_{L^2_t W^{-1,\frac{18}{13}}_x}
	\intT{\!\!\!\int_0^T\pi_\etat(t-s) (t-s)^{1/2}
		\|\pi_\etax \ast\vu(s,x)  \pi_\etax \ast(\vu(t,x) - \vu(s,x))\|_{W^{1,\frac{18}{5}}_x}\, {\rm d}s}\\	& \leq C\etax^\theta+C
	\intT{\!\!\!\int_0^T\pi_\etat(t-s) (t-s)^{1/2}\lr{\|\pi_\etax \ast\vu(t,x)\|^2_{W^{1,\frac{18}{5}}_x}
			+ \|\pi_\etax \ast\vu(s,x)\|^2_{W^{1,\frac{18}{5}}_x}} {\rm d}s}\\
	& \leq C\etax^\theta+C\etax^{-\alpha}\etat^\beta,
}
where we use the fact that there exists $\alpha>0$ such that
$ \|\pi_\etax \ast\vu(s,x)\|_{L^2_t W^{1,\frac{18}{5}}_x} \leq \etax^{-\alpha}.$

Choosing $\etat$ sufficiently small w.r.t. $\etax$, for example $\etat=\etax^{\frac{\alpha}{\beta}+\frac{\vt}{\beta}}$, we finally estimate the sum of $J_1$ and $J_2$ as
\eq{\nonumber
	&\intTOD{\vr(\vu-\pi_\eta\ast \vu)^2}
	\leq  C\etat^{\theta'},
}
for some constant $\theta'>0$. From the above relation, we want to deduce the following:
\begin{equation}\nonumber
	\intTOD{ \vr |\vu -\pi_\etat \ast \vu|} \leq C\etat^{\theta'}.
\end{equation}

To do this, we consider two cases:
\begin{equation}\nonumber
	|\vu -\pi_\etat \ast \vu(t,x)| \leq \etat^{\frac{\theta'}2}
\mbox{ \ \ and \ \ } 
	|\vu(t,x) -\pi_\etat \ast \vu(t,x)| > \etat^{\frac{\theta'}2}.
\end{equation}
Then a straightforward computation yields
\begin{equation}\nonumber
	\intTOD{ \vr |\vu -\pi_\etat \ast \vu|} \leq
	\intTOD{\vr \lr{\etat^{\frac{\theta'}2}+ \etat^{-\frac{\theta'}2} |\vu-\pi_\etat \ast \vu|^2}} \leq C \etat^{\frac{\theta'}2}.
\end{equation}
Now, taking $\eta=\etat$ and $\vt=\frac{\theta'}{2}$ we conclude the proof. $\Box$

\subsection{Proof of Lemma \ref{Lemma:Du}}\label{Lemma13}
At first we use H\"older inequality to obtain
\eqh{
&\int_{\Td} \frac{{\rm d}z}{(|z|+h)^d} \| D_{|z|}u(\cdot) - 
D_{|z|}u(\cdot + z)\|_{L^2}  \\
&\leq \lr{\int_{\Td} \frac{{\rm d}z}{(|z|+h)^d}}^{1/2} 
\left( 
\int_{\T^d} \frac{{\rm d}z}{(|z|+h)^d} \| D_{|z|}u(\cdot) - 
D_{|z|}u(\cdot + z)\|_{L^2}^2\right)^{1/2}
}
The first component is proportional to $|\ln h|^{1/2}$, and so we only need to show that the latter term is bounded by $\|u\|_{W^{1,2}(\Td)}$.

To avoid technical problems with properties of functions on torus we extend $u$ onto the whole $\R^d$ in such a way that $u$ is treated as a periodic function and then it is localized in such a way that 
it is $u$ for $x \in [-3\pi,3\pi]^d$ and zero for $\R^d\setminus [-4\pi,4\pi]^d$ (as the torus is $\T^d=[-\pi,\pi]^d$) and we name it $Eu$. Then it is clear that $\|Eu\|_{W^{1,2}(\R^d)} \leq C\|u\|_{W^{1,2}(\T^d)}$ and 
\begin{equation}\nonumber
    \|D_{|z|} u(\cdot)-D_{|z|}u(\cdot+z)\|_{L^2(\T^d)} \leq 
     \|D_{|z|}Eu(\cdot)-D_{|z|}Eu(\cdot+z)\|_{L^2(\R^d)} .
\end{equation}
{Using the above inequality, we have}
\eq{\nonumber
&\int_{\T^d} \frac{{\rm d}z}{(|z|+h)^d} \| D_{|z|}u(\cdot) - 
D_{|z|}u(\cdot + z)\|_{L^2(\T^d)}^2\\
&\leq\int_{\T^d} \frac{{\rm d}z}{(|z|+h)^d}
\int_{\R^d}\left|\frac{1}{|z|}\lr{\int_{B(0,|z|)} \frac{|\nabla Eu(x+w)|}{|w|^{d-1} } {\rm d}w-
\int_{B(0,|z|)} \frac{|\nabla Eu(x+z+w)|}{|w|^{d-1} } {\rm d}w} \right|^2\dx
\\
&=
\int_{\T^d} \frac{dz}{(|z|+h)^d} 
\int_{\R^d} \left|
\frac{1}{|z|} \int_{B(0,|z|)} \frac{{\rm d}w}{|w|^{d-1}} 
e^{ {iw \cdot \xi}} (1-e^{{iz\cdot \xi}})\widehat{|\nabla Eu|}(\xi) \right|^2 {\rm d} \xi
\\
&=\int_{\R^d} {\rm d}\xi \widehat{|\nabla Eu|}^2(\xi) 
\int_{\T^d} \frac{{\rm d}z}{(|z|+h)^d} \frac{|1-e^{iz\cdot \xi}|^2}{|z|^2}
\left| \int_{B(0,|z|)} \frac{{\rm d}w}{|w|^{d-1}} e^{iw\cdot \xi} \right|^2.
}
We will next show that the integral
\eq{\nonumber
I=
\int_{\T^d} \frac{{\rm d}z}{(|z|+h)^d} \frac{|1-e^{iz\cdot\xi}|^2}{|z|^2}
\left| \int_{B(0,|z|)} \frac{{\rm d}w}{|w|^{d-1}} e^{iw\cdot \xi} \right|^2\leq \int_{\R^d} \frac{{\rm d}z}{|z|^d} \frac{|1-e^{iz\cdot\xi}|^2}{|z|^2}
\left| \int_{B(0,|z|)} \frac{{\rm d}w}{|w|^{d-1}} e^{iw\cdot\xi} \right|^2
}
is independent of $\xi$ and is bounded by a constant.

We first denote $\xi=|\xi|\omega_\xi$, and then by changing the variables for $t=|\xi| z$ and $v=|\xi|w$ we get the estimate
\begin{equation}\label{est_I}
I\leq \int_{\R^d} \frac{{\rm d}t}{|t|^d}
\frac{|1-e^{i\omega_\xi \cdot t}|^2}{|t|^2}\left| 
\int_{B(0,|t|)} \frac{{\rm d}v}{|v|^{d-1}} e^{i\omega_\xi \cdot v} \right|^2.
\end{equation}
 Because both integrals are spherically symmetric, we can fix $\omega_\xi$ to be a unit vector, say  $\omega_\xi=\hat{e}_1$ .
So the integral $I$ is independent of $\xi$, and therefore in order to end the proof one should only show that it is bounded. 
We denote $v=(v_1,v')$, and for now we focus on the second integral only, i.e. $ \int_{B(0,|t|)} \frac{dv}{|v|^{d-1}} e^{iv_1}$.\\ 

{\emph Case 1.} For $t$ small we can estimate $|e^{iv_1}|\leq1$ and the integral is smaller than $Ct$. \\

{\emph Case 2.} For $t>>1$  we split  the integral as follows
\eq{\label{I2}
   \int_{B(0,|t|)} \frac{{\rm d}v}{|v|^{d-1}} e^{iv_1}= 2\int_0^{|t|} {\rm d}v_1 \int_{{ B'(0,\sqrt{|t|^2-v_1^2})}} e^{iv_1}
    \frac{{\rm d}v'}{(\sqrt{v_1^2+|v'|^2})^{(d-1)/2}},
}
where $B'$ denotes a ball in $d-1$ dimensions. We first write $e^{iv_1}=\cos(v_1)+i\sin(v_1)$. We split the real part of the integral $\int_0^{|t|}$ into the integrals on the intervals $(0,\frac\pi2)$, $(-\frac\pi2+k\pi,-\frac\pi2+ (k+1)\pi)$, $(-\frac\pi2+(k_{|t|}+1)\pi,|t|)$ for $k=1,\ldots,k_{|t|}$, $k_{|t|}=\left\lfloor\frac{t-\pi/2}{\pi}\right\rfloor$.
The idea is just to show that above division gives rise to the alternating series, just to show that the integral is uniformly bounded in the meaning of Riemann's integral i.e. 
\eq{\nonumber
\frac{1}{2}\int_{B(0,|t|)} \frac{{\rm d}v}{|v|^{d-1}} e^{iv_1}=\int_0^{\pi/2} \Gamma\, {\rm d}v_1+ \sum_{k=1}^{k_{|t|}} (-1)^k\int_{-\pi/2+k\pi}^{-\pi/2+(k+1)\pi} 
\Gamma\, {\rm d}v_1+ (-1)^{k_{|t|}+1}\int_{-\frac\pi2+(k_{|t|}+1)\pi}^{|t|} \Gamma\, {\rm d}v_1, 
}
where we denoted
$$\Gamma= |\cos{v_1}| \int_{B'(0,\sqrt{|t|^2-v_1^2})} 
    \frac{{\rm d}v'}{(\sqrt{v_1^2+|v'|^2})^{(d-1)/2}}.$$
Note that the irrational part of the integral, related to $i\sin v_1$, is odd, so it sums up to zero.    
    
In the above we have used that $\cos t = -\cos(\pi +t)$. We see that the integrals of $|\Gamma|$ are decreasing with  $k$, and  so the series is bounded by the first term, i.e. we have
\eq{\nonumber
\frac{1}{2}\left|\int_{B(0,|t|)} \frac{{\rm d}v}{|v|^{d-1}} e^{iv_1}\right|\leq\int_0^{\pi/2} \Gamma\, {\rm d}v_1+ \int_{\pi/2}^{3\pi/2} 
\Gamma\, {\rm d}v_1.
}
We now want to estimate these integrals, for this we include the irrational part again, and estimate \eqref{I2} on the intervals $(0,\frac\pi2)$, $(\frac\pi2,\frac{3\pi}{2})$, assuming that $t > 2\pi$, we can write that
\eq{\nonumber
&\int_0^{\pi/2} {\rm d}v_1 \int_{B'(0,\sqrt{|t|^2-v_1^2})} \frac{{\rm d} v'}{|v|^{d-1}} \leq
\int_0^{\pi/2} {\rm d}v_1 \int_{B'(0,|t|)} \frac{{\rm d} v'}{|v|^{d-1}}\\
&\leq
C \int_0^{\pi/2} {\rm d}v_1 \int_{B'(0,|t|)}  \frac{ {\rm d} v' }{
(|v_1|+|v'|)^{d-1}} 
\leq
C \int_0^{\pi/2} {\rm d}v_1 \int_0^{|t|} \frac{ {\rm d} |v'|\, |v'|^{d-2} }{
(|v_1|+|v'|)^{d-1}} \\
&\leq  C \int_0^{\pi/2} {\rm d}v_1 \int_0^{|t|} \frac{ {\rm d} |v'|}{|v_1|+|v'|}=
C\int_0^{\pi/2} {\rm d}v_1 (\ln( |t|+|v_1|) - \ln |v_1|) \leq C (\ln |t| +1).
}
For the second part the same way we obtain
\begin{equation}\nonumber
\int_{\pi/2}^{3\pi/2} {\rm d}v_1 \int_{B'(0,|t|)} \frac{{\rm d}v'}{|v|^{d-1}} \leq C\int_{\pi/2}^{3\pi/2} {\rm d}v_1 \int_0^{|t|} \frac{{\rm d}|v'|}{\frac{\pi}{2} + |v'|} \leq C(\ln |t| +1).
 \end{equation}
Since $t>2\pi$, we have $1+\ln |t| \leq C \ln |t|$, so, summing up cases 1 and 2, we get that
\begin{equation}\nonumber
  \left|\int_{B(0,|t|)} \frac{{\rm d}v}{|v|^{d-1}} e^{iv_1} \right|  \leq C\min\{ t, \ln |t|\}.
\end{equation}
At last, coming back to \eqref{est_I}, we obtain
\eq{\nonumber
I\leq C\int_{\R^d} \frac{{\rm d}t}{|t|^d}
\frac{|1-e^{i\omega_\xi t}|^2}{|t|^2}(\min\{ t, \ln |t|\})^2
\leq
  C\lr{\int_{B(0,1)} \frac{{\rm d}t}{|t|^{d-2}}  +\int_{B(0,1)^c} \frac{{\rm d}t}{|t|^d}
\frac{\ln| t|^2}{|t|^2}  }\leq C,
}
where we used that for $|t|$ small, i.e. for $t\in B(0,1)$ we can estimate $|1-e^{i\omega_\xi t}|\leq |t|$. With this estimate in hand, we complete the proof. $\Box$

\subsection{Proof of  Lemma \ref{interpolation-lemma}}\label{Lemma15}

Before we begin the proof, we introduce the following notation: For $0<a,b \in \R$, the notation $a\approx b$ means that there exists a $ \mu >0 $ such that \[ \mu a \leq b \leq \frac{1}{\mu} a .\]
From Parseval's theorem for Fourier transform, we have
\begin{align*}
	\int_{\mathbb{R}} \int_{\mathbb{R}^d}  \Big((- \Delta_x)^{-1}\dv_x \pt \pmb{\varphi} \Big) W(t,x) \dx\; \dt  &= 	\int_{\mathbb{R}} \int_{\mathbb{R}^d}   \mathcal{F}_{\xi_0,\xi} {\Big((- \Delta_x)^{-1}\dv_x \pt \pmb{\varphi} \Big)}   \overline{\mathcal{F}_{\xi_0,\xi} ({W} )}\text{d}\xi\; \text{d}\xi_0\\
	&=-	\int_{\mathbb{R}} \int_{\mathbb{R}^d}  \;\xi_0 \vert \xi \vert^{-2} \; \xi \cdot \; \mathcal{F}_{\xi_0,\xi}({\pmb{\varphi}}) \;  \overline{\mathcal{F}_{\xi_0,\xi} (W)} \text{d}\xi\; \text{d}\xi_0.
\end{align*}
Using $\mathcal{F}_{\xi_0,\xi}(\pmb{\varphi})= |\xi| \mathcal{F}_{\xi_0,\xi}({\pmb{g}}) $ and $ 0<\beta <1 $, we obtain 
\begin{align*}
	&\int_{\mathbb{R}} \int_{\mathbb{R}^d}  \;\xi_0 \vert \xi \vert^{-2} \; \xi \cdot \; \mathcal{F}_{\xi_0,\xi}({\pmb{\varphi}}) \;  \overline{\mathcal{F}_{\xi_0,\xi} (W)} \text{d}\xi\; \text{d}\xi_0\\
	&= \int_{\mathbb{R}} \int_{\mathbb{R}^d}  
	\;\left( \xi_0 |\xi_0|^{(\beta-1)} \vert \xi \vert^{-\beta} \xi \cdot \mathcal{F}_{\xi_0,\xi}({\pmb{g}})\right) \;\left( \vert \xi_0\vert^{1-\beta} \vert\xi\vert^{-(1-\beta)} \overline{\mathcal{F}_{\xi_0,\xi} (W)} \right) \text{d}\xi\; \text{d}\xi_0\\
	&= \int_{\mathbb{R}} \int_{\mathbb{R}^d}  \; \mathcal{F}^{-1}_{t,x} \left( \xi_0 |\xi_0|^{\beta-1} \vert \xi \vert^{-\beta} \xi \cdot  \mathcal{F}_{\xi_0,\xi}({\pmb{g}})\right) \overline{\mathcal{F}^{-1}_{t,x} } \left( \vert \xi_0\vert^{1-\beta} \vert\xi\vert^{-(1-\beta)} \overline{\mathcal{F}_{\xi_0,\xi} (W)} \right) \dx\; \dt  \\
	& \leq C \left\Vert \mathcal{F}^{-1}_{t,x} \left( |\xi_0|^\beta \vert \xi \vert^{1-\beta}  \mathcal{F}_{\xi_0,\xi}({\pmb{g}})\right)  \right\Vert_{L^p_{t,x}}  \left\Vert \mathcal{F}^{-1}_{t,x}  \left( \vert \xi_0\vert^{1-\beta} \vert\xi\vert^{-(1-\beta)} \mathcal{F}_{\xi_0,\xi} (W) \right)  \right\Vert_{L^{p^{\prime}}_{t,x}},
\end{align*}
with $ \frac{1}{p}+ \frac{1}{p^\prime} =1 $. 
Now our idea is to use the following observation 
\begin{align*}
	\left\Vert \mathcal{F}^{-1}_{t,x} \left( |\xi_0|^\beta \vert \xi \vert^{1-\beta}  \mathcal{F}_{\xi_0,\xi}(\pmb{g})\right)  \right\Vert_{L^p_{t,x}} \approx \| \pmb{g} \|_{\dot{W}^{\beta,p}(\R; \dot{W}^{1-\beta,p}(\R^d))}, 
\end{align*}
which means that
\begin{align*}
	\mathcal{F}^{-1}_{t}\left( |\xi_0|^\beta \mathcal{F}_{\xi_0}(\pmb{g})\right) \in 
	L^p(\R;W^{1-\beta,p}(\R^d)).
\end{align*}

Then we want to interpolate $ \| \pmb{g} \|_{\dot{W}^\beta_p(\R; \dot{W}^{1-\beta,p}(\R^d))}  $ with $ \| \pmb{g} \|_{L_{q}(\R^{d+1})}  $ and $ \| \partial_t \pmb{g} \|_{L^{q^\prime}(\R; W^{-1,q^\prime}(\R^{d+1}))}  $ for {some} $ q $ and $ q^\prime $. 


It is well known that for the homogeneous Triebel-Lizorkin space $ \dot{F}^s_{p,2}	(\R^d)  $, we have
\begin{align}\label{eq-sob-tri}
	\dot{W}^{s,p}(\R^d)=\dot{F}^s_{p,2}	(\R^d) \text{ for } \; s\in \R \text{ and }1<p<\infty  
\end{align}
with equivalent norms. Therefore, using Duoandikoetxea \cite[Theorem 8.7]{Du}, we obtain
\begin{align*}
	\left\Vert \mathcal{F}^{-1}_{t,x} \left( |\xi_0|^\beta \vert \xi \vert^{1-\beta}  \mathcal{F}_{\xi_0,\xi}(\pmb{g})\right)  \right\Vert_{L^p_{t,x}}  \approx	\left\Vert\left( \sum_{k,l} 2^{2k(1-\beta)} 2^{2l\beta}  {\pmb{g}}_{k,l}^2\right)^{\frac{1}{2}} \right\Vert _{L^p_{t,x}},
\end{align*}
where \begin{align*}
	{\pmb{g}}_{k,l}= \dot{\Delta}_k^x \dot{\Delta}_l^t \mathcal{F}_{\xi_0,\xi} ({\pmb{g}}), 
\end{align*}
with $ \dot{\Delta}_k^x  $ and $ \dot{\Delta}_l^t  $  are the homogeneous dyadic blocks in space variable and time variable, respectively.  Next, using Holder inequality, we obtain 
\begin{align*}
	\left\Vert\left( \sum_{k,l} 2^{2k(1-\beta)} 2^{2l\beta}   {\pmb{g}}_{k,l}^2\right)^{\frac{1}{2}} \right\Vert _{L^p_{t,x}} &\leq 	\left\Vert\left( \sum_{k,l} 2^{2k}{\pmb{g}}_{k,l}^2\right)^{\frac{1-\beta}{2}}  \left( \sum_{k,l}  2^{2l}  {\pmb{g}}_{k,l}^2\right)^{\frac{\beta}{2}} \right\Vert _{L^p_{t,x}}\\
	& \leq 	\left\Vert\left( \sum_{k,l} 2^{2k}{\pmb{g}}_{k,l}^2\right)^{\frac{1-\beta}{2}}\right\Vert _{L^{r_1}_{t,x}} \left\Vert   \left( \sum_{k,l}  2^{2l}  {\pmb{g}}_{k,l}^2\right)^{\frac{\beta}{2}}\right \Vert _{L^{r_1^\prime}_{t,x}},
\end{align*}
with $ \frac{1}{r_1} +  \frac{1}{r_1^\prime} =\frac{1}{p}$. Considering $ r= \frac{r_1}{p} $ and $ r^\prime = \frac{r_1^\prime}{p}  $, we introduce
\begin{align}\nonumber
	q=(1-\beta )p r \text{ and } q^\prime = \beta p r^\prime 
\end{align}
to deduce 
\begin{align*}
	\left\Vert\left( \sum_{k,l} 2^{2k(1-\beta)} 2^{2l\beta}   {\pmb{g}}_{k,l}^2\right)^{\frac{1}{2}} \right\Vert _{L^p_{t,x}}  \leq 	\left\Vert\left( \sum_{k,l} 2^{2k}{\pmb{g}}_{k,l}^2\right)^{\frac{1}{2}}\right\Vert _{L^{q}_{t,x}}^{1-\beta} \left\Vert   \left( \sum_{k,l}  2^{2l}  {\pmb{g}}_{k,l}^2\right)^{\frac{1}{2}}\right \Vert _{L^{q^\prime}_{t,x}}^{\beta}
\end{align*}
with the following relation between $ p, q, q^\prime $ and $ \beta $: 
\begin{align*}
	1=p \left( \frac{1-\beta}{q} + \frac{\beta}{q^\prime} \right).
\end{align*}
Again, using the equivalence of Sobolev and Triebel-Lizorkin space \eqref{eq-sob-tri}, we have 
\begin{align*}
	\left\Vert\left( \sum_{k,l} 2^{2k}{\pmb{g}}_{k,l}^2\right)^{\frac{1}{2}}\right\Vert _{L^{q}_{t,x}} \approx \left\Vert \mathcal{F}^{-1}_{t,x}  \left(\vert \xi \vert \; \mathcal{F}_{\xi_0,\xi}({\pmb{g}})\right) \right\Vert _{L^q_{t,x}} = \left\Vert   \mathcal{F}^{-1}_{t,x}  \left(\; \mathcal{F}_{\xi_0,\xi}({\pmb{\varphi}})\right) \right\Vert _{L^q_{t,x}} = \Vert \pmb{\varphi} \Vert_{L^q_{t,x}}
\end{align*}
and 
\begin{align*}
	\left\Vert   \left( \sum_{k,l}  2^{2l}  {\pmb{g}}_{k,l}^2\right)^{\frac{1}{2}}\right \Vert _{L^{q^\prime}_{t,x}}\approx & \left\Vert  \mathcal{F}^{-1}_{t,x}  \left(\vert \xi_0 \vert \; \mathcal{F}_{\xi_0,\xi}({\pmb{g}})\right) \right\Vert _{L^{q^\prime}_{t,x}}\\&=\left\Vert  \mathcal{F}^{-1}_{t,x}  \left(\vert \xi_0 \vert \; \vert \xi \vert^{-1} \mathcal{F}_{\xi_0,\xi}({\pmb{\varphi}})\right) \right\Vert _{L^{q^\prime}_{t,x}}= \left\Vert  \partial_t \pmb{\varphi} \right\Vert _{L^{q^\prime}_tW^{-1,q^\prime}_x}.
\end{align*}
The above two estimates yield
\begin{align}\label{il:e1}
	\left\Vert \mathcal{F}^{-1}_{t,x} \left( |\xi_0|^\beta \vert \xi \vert^{1-\beta}  \mathcal{F}_{\xi_0,\xi}(\pmb{g})\right)  \right\Vert_{L^p_{t,x}} \leq \Vert \pmb{\varphi}  \Vert_{L^{q}_{t,x}}^{1-\beta} \Vert \partial_t \pmb{\varphi} \Vert_{L^{q^\prime}_tW^{-1,q^\prime}_x}^{\beta}. 
\end{align}
For the term $ \left\Vert \mathcal{F}^{-1}_{t,x}  \left( \vert \xi_0\vert^{1-\beta} \vert\xi\vert^{-(1-\beta)} \mathcal{F}_{\xi_0,\xi} (W) \right)  \right\Vert_{L^{p^{\prime}}_{t,x}} $, we proceed analogously and obtain
\begin{align}\label{il:e2}
	\left\Vert \mathcal{F}^{-1}_{t,x}  \left( \vert \xi_0\vert^{1-\beta} \vert\xi\vert^{-(1-\beta)} \mathcal{F}_{\xi_0,\xi} (W) \right)  \right\Vert_{L^{p^{\prime}}_{t,x}} \leq \Vert W \Vert_{L^{\bar{q}}_{t,x}}^{1-\alpha} \Vert \partial_t W \Vert_{L^{\bar{q}^{\prime}}_tW^{-1,\bar{q}^\prime}_x}^{\alpha}
\end{align}
where $ 1-\beta =\alpha $ and  the relation between $p^\prime, \bar{q}, {\bar{q}}^\prime $ and $ \alpha $ is given by 
\begin{align*}
	1=p^\prime \left( \frac{1-\alpha}{\bar{q}} + \frac{\alpha}{\bar{q}^\prime} \right).
\end{align*}
Moreover, we have \[ 1-\frac{1}{q^\prime}-\frac{1}{\bar{q}}=\alpha\left( \frac{1}{q} + \frac{1}{\bar{q}^\prime} -\frac{1}{q^\prime}-\frac{1}{\bar{q}} \right). \]
The assumption \eqref{il:a1} implies $ 0<\alpha <1  $.
Summing up the estimates \eqref{il:e1} and \eqref{il:e2}, we get our desired result. $\Box$
\subsection{Existence of solutions to parabolic systems}\label{Lemma16}
\begin{lemma}\label{exis-EVF} Let $p>1$.
Let $a$ be sufficiently smooth and satisfies 
\[ 0< a_* \leq a(t,x) \leq a^* < \infty,\; \nabla a \in L^\infty(0,T;L^\infty(\T^d)) \text{ and }a \in C^\alpha(0,T;C_b(\T^d)) \text{ for }\alpha>0 \]
Moreover, we assume 
$$\vu_0 \in W^{1-2/p,p}(\T^d),\; \mathbf{F} \in L^p(0,T;L^p(\T^d)\text{ and } \mathbf{G}\in L^p(0,T;L^p(\T^d))
$$ 
Then the solution $ \vu$ of
\begin{equation}\label{oo1}
    a\pt \vu -\Div \vS (\vu) = \Div \mathbf{F} +\mathbf{G}, \qquad \vu|_{t=0}=\vu_0
\end{equation}
satisfies the following bound
\begin{equation}\label{oo2}
    \|\pt \vu\|_{L^p_t W^{-1,p}_x} + \|\nabla \vu \|_{L^p_{t,x}} \leq C(T)(\| \mathbf{F}\|_{L^p_t L^p_x}+\|\mathbf{G}\|_{L^p_t L^p_x}+
    \|\vu_0\|_{W^{1-2/p,p}_x}).
\end{equation}
\end{lemma}

\pf The result belongs to the classical theory of parabolic systems. Since such results are 
consequences of several results from \cite{DHP} and \cite{LSU} we here give a sketch of the proof to satisfy readers. 
The first step which is required by the structure of regularity is reformulation of (\ref{oo1}) in the following form:
\begin{equation}\nonumber
    \pt \vu -b\Div \vS (\vu) = \Div \tilde{\mathbf{F}} +\tilde{\mathbf{G}} , \qquad u|_{t=0}=u_0 \mbox{ \ \ with \ \ } b=1/a,\; \tilde{\mathbf{F}}= \frac{\mathbf{F}}{b} \text{ and } \tilde{\mathbf{G}}=\mathbf{G} - \mathbf{F} \nabla\lr{\frac{1}{b}}
\end{equation}
The proof we perform for a simplified version, to avoid unnecessary technicalities. Therefore, instead of $\Div \mathbf{S}(\vu)$, we consider
\begin{equation}\label{oo5}
    \pt \vu -  b \Delta \vu = \nabla f,
\end{equation}
with $ f \in L^p(0,T;L^p(\T^d))$.

{\sc Step 1.} The basic result is: if $\bar b \in [b_*,b^*]$ it is a constant then
the solution to 
\begin{equation*}
    \pt \vu - \bar b \Delta \vu = \nabla f
\end{equation*}
fulfills
\begin{equation}\nonumber
    \|\pt \vu\|_{L^p_t W^{-1,p}_x} + \|\nabla \vu\|_{L^p_{t,x}} \leq C(b_*,b^*)(\|f\|_{L^p_t L^{p}_x}+
    \|\vu_0\|_{W^{1-2/p,p}_x}).
\end{equation}

{\sc Step 2.} In this part we fix the time in the coefficient $b$, i.e. we consider $a$ time-independent in   (\ref{oo2}). Let $\pi_k$ be a partition of unity. We are required to 
\begin{equation*}
    {\rm diam \;\; supp\;} \pi_k \leq \lambda.
\end{equation*}

Let $k$ be fixed, and let us set $\pi=\pi_k$. We consider
\begin{equation}\nonumber
    \pt (\pi \vu) - b \Delta (\pi \vu) = \nabla (\pi f) - \nabla \pi f + b \nabla \pi \nabla \vu + b
    \Delta \pi \vu.
\end{equation}
Let $\bar b = b(x_k)$ such that $\pi_k(x_k)=1$, then
\begin{equation}\nonumber
    \pt(\pi \vu) - \bar b \Delta (\pi \vu) = 
    \nabla (\pi f) - \nabla \pi f + b \nabla \pi \nabla \vu + b
    \Delta \pi \vu +(b-\bar b) \Delta (\pi \vu). 
\end{equation}

We just deal with the key elements of the r.h.s.. First, we take $f\nabla \pi$, which is of higher order but with bad dependence form $\lambda$. This part is formally related with the following problem
\begin{equation}\nonumber
    \pt \mathbf{w}-\bar b\Delta \mathbf{w} = f \nabla\pi, \qquad \mathbf{w}|_{t=0} =0.
\end{equation}
Then we have
\begin{equation}\nonumber
    \sup_t \|\mathbf{w}\|_{W^{2-2/p,p}_x} + \|\pt \mathbf{w},\nabla^2 \mathbf{w}\|_{L^p_{t,x}} \leq C\|f\nabla \pi\|_{L^p_{t,x}}
\end{equation}
with constant independent of $T$. Then we get (for $p>2$) that
\begin{equation}\nonumber
    \|\nabla \mathbf{w}\|_{L^p_{t,x}} \leq CT^{1/p}\sup_t \|\mathbf{w}\|_{W^{1-2/p,p}_x} \leq CT^{1/p} \|f\nabla \pi\|_{L^p_{t,x}} \leq CT^{1/p}\lambda^{-1} \|f\|_{L^p((0,T) \times supp\,\pi ))}.
\end{equation}
The second term is
\begin{align*}
     \|(b-\bar b) \Delta (\pi \vu)\|_{L^p_t W^{-1,p}_x} &\leq 
    \|\Div ((b-\bar b) \nabla (\pi \vu)\|_{L^p_t W^{-1,p}_x}+
    +\|\nabla b \nabla (\pi \vu)\|_{L^p_t W^{-1,p}_x} \\
    &\leq C \|b-\bar b\|_{L^\infty(supp\;\pi)}\|\nabla (\pi \vu)\|_{L^p((0,T) \times supp\,\pi ))}
    + CT^{1/p}\|\nabla b\|_{L^\infty}\|\nabla (\pi \vu)\|_{L^p_{t,x}}\\
    &\leq C(\lambda + T^{1/p}) \|\nabla b\|_{L^\infty} \|\nabla (\pi \vu)\|_{L^p_{t,x}}.
\end{align*}
Thus all together we obtain
\begin{align*}
     \|\nabla(\pi \vu)\|_{L^p_{t,x}} \leq C &\big( (1+T^{1/p}\lambda^{-1})\|f\|_{L^p((0,T)\times (supp\,\pi))} \\
  &+    (\lambda + T^{1/p}) \|\nabla (\pi \vu)\|_{L^p_{t,x}} +
  T^{1/p}\lambda^{-1} \|\nabla \vu\|_{L^p(0,T;L^p(supp\,\pi))}\big).
\end{align*}
Now as $\lambda + T^{1/p}$ is small, we deduce
\begin{equation}\nonumber
    \|\nabla(\pi \vu)\|_{L^p_{t,x}}^p \leq C ( (1+T^{1/p}\lambda^{-1})^p\|f\|_{L^p(0,T;L^p(supp\,\pi))}^p 
  +    (T^{1/p}\lambda^{-1})^p \|\nabla \vu\|_{L^p(0,T;L^p(supp\,\pi))}^p).
\end{equation}
Let $N$ be the cover number of the union of  supp $\pi_k$, then we can estimate
\begin{align*}
  &  \|\nabla \vu \|^p_{L^p_{t,x}} \leq N^{p-1} \sum_k \|\nabla (\pi_k \vu)\|_{L^p_{t,x})}^p 
    \\
  &  \leq CN^{p-1} \sum_{k} \left[ (1+T^{1/p}\lambda^{-1})^p\|f\|_{L^p(0,T;L^p(supp\,\pi_k))}^p 
  +    (T^{1/p}\lambda^{-1})^p \|\nabla \vu\|_{L^p(0,T;L^p(supp\,\pi_k))}^p \right]
 \\
 &\leq CN^{2p-1} \|f\|^p_{t,x} + CN^{2p-1}\epsilon \|\nabla \vu\|^p_{L^p_{t,x}}   
\end{align*}
as $T^{1/p}\lambda^{-1} <1$ and $(T^{1/p}\lambda^{-1})^p \leq \epsilon$. For $\epsilon$ small the last term can be absorbed by the l.h.s. and we get
\begin{equation*}
    \|\nabla \vu\|_{L^p_{t,x}}\leq C\|f\|_{L^p_{t,x}}.
\end{equation*}
Hence for $b$ time-independent we get
\begin{equation}\label{oo6}
    \|\pt \vu\|_{L^p_t W^{-1,p}_x} + \|\nabla \vu\|_{L^p_{t,x}} \leq C(b_*,b^*,\|\nabla b\|_{L^\infty})(\|f\|_{L^p_t L^{p}_x}+
    \|\vu_0\|_{W^{1-2/p,p}_x}).
\end{equation}
Now we need to make $b$ time dependent.

{\sc Step 3.} As $b$ is time dependent we can consider equation (\ref{oo5}) in the form
\begin{equation}\nonumber
    \pt \vu -b(0,x)\Delta \vu = (b(t,x)-b(0,x)) \Delta \vu+ \nabla f, \qquad \vu|_{t=0}=\vu_0. 
\end{equation}
Then we restate the first term of the l.h.s. as follows
\begin{equation*}
    (b(t,x)-b(0,x)) \Delta \vu = \Div ((b(x,t)-b(x,0)) \nabla \vu + \nabla b \nabla \vu.
\end{equation*}
So from the Step 2 we know that this term generates the following impact on the estimate
\begin{align*}
     \|b(t,x)-b(0,x)\|_{L^\infty_x} \|\nabla \vu\|_{L^p_{t,x}} + CT^{1/p}\|\nabla b\|_{L^\infty_{t,x}}\|\nabla \vu\|_{L^p_{t,x}}\\
      \leq C(\lambda^\alpha + T^{1/p})(\| b\|_{C^\alpha_t L^\infty_x} + \|\nabla b\|_{L^\infty_{t,x}})\|\nabla \vu\|_{L^p_{t,x}},
\end{align*}
where we used the fact that $b\in C^{\alpha}(0,T;C_b)$. Therefore, for short time  we get (\ref{oo6}) in form 
\begin{equation}\nonumber
    \|\pt \vu\|_{L^p_t W^{-1,p}_x} + \|\nabla \vu\|_{L^p_{t,x}} \leq 
    C(b_*,b^*,\|\nabla b\|_{L^\infty_{t,x}},\|b\|_{C^\alpha_t L^\infty_x} )(\|f\|_{L^p_t L^{p}_x}+
    \|\vu_0\|_{W^{1-2/p,p}_x}).
\end{equation}
The above inequality is valid for time interval $[0,T_*]$
for small $T_*$, then from the trace theorem we control
the $W^{1-2/p,p}_x$ norm  of $\vu|_{t=T_*}$ in terms of the l.h.s. of above and $\vu_0$. From that, we can obtain the estimate for $[T_*,2T_*]$ with the same constant $C$ considering
$\vu|_{t={T_*}}$ instead of $\vu_0$.
Thus the inequality can be continued step by step with increase of the constant in (\ref{oo2}) with factor $e^{CT}$.  The proof of the Lemma is concluded.
$\Box$

\subsection* {Acknowledgments:}
The second (PBM) author has been partially supported by National Science Centre grant 
\newline No. 2018/29/B/ST1/00339 (Opus). The research of N.C. and of E.Z. leading to these results has received funding from the EPSRC Early Career Fellowship no. EP/V000586/1.



\begin{thebibliography}{99}

\bibitem{Belgacem} 
\sc F. Ben Belgacem, P.-E. Jabin. 
\rm Compactness for nonlinear continuity equations, 
\it J. Funct. Anal. \rm 264 (1), 139--168, 2013. 

\bibitem{DPSV}
\sc T. Debiec, B. Perthame, M. Schmidtchen, N. Vauchelet.
\rm Incompressible limit for a two-species model with coupling through Brinkman's law in any dimension,
\it Journal de Mathématiques Pures et Appliquées,  
\rm 145, 204-239, 2021.


\bibitem{BJ} 
\sc D. Bresch, P.-E. Jabin.
\rm Global Existence of Weak Solutions for Compresssible Navier--Stokes Equations: Thermodynamically unstable
pressure and anisotropic viscous stress tensor.
\it  Ann. of Math., 
\rm 188 (2), 577--684, 2018.

\bibitem{BJ2} 
\sc D. Bresch, P.-E. Jabin.
\it Global weak solutions of PDEs for compressible media: A compactness criterion to cover new physical situations.
\rm Shocks, singularities, and oscillations in nonlinear optics and fluid mechanics,  33–54, Springer INdAM Ser., 17, Springer, Cham, 2017.


\bibitem{BJW}
\sc D. Bresch, P.-E. Jabin, F. Wang.
\rm Compressible Navier{\textendash}Stokes equations with heterogeneous pressure laws.
\it {Nonlinearity}
\rm 34 (6), 4115--4162, 2021.


\bibitem{BMZ} 
\sc D. Bresch, P.B. Mucha, E. Zatorska.
\rm {Finite-energy solutions for compressible two-fluid Stokes system.} 
\it{Arch. Ration. Mech. Anal.} 
\rm {232},  987--1029, 2019.


\bibitem{DHP}
\sc R. Denk, M. Hieber, J. Prüss.
\rm R-boundedness, Fourier multipliers and problems of elliptic and parabolic type. 
\it Mem. Amer. Math. Soc. 
\rm166, no. 788, pp. viii+114, 2003.

\bibitem{Du}
\sc J. Duoandikoetxea.
\it Fourier analysis. 
\rm American Mathematical Society, Providence, RI, 2001.
\bibitem{LSU}
\sc O. A. Lady\v zenskaja, V. A. Solonnikov, N. N. Uraltseva.
\it Linear and quasilinear equations of parabolic type.
\rm American Mathematical Society, Providence, RI, 1968.
\bibitem{EF2001}
\sc E.  Feireisl.
\rm On compactness of solutions to the compressible isentropic {N}avier-{S}tokes equations when the density is not square integrable.
 \it Comment. Math. Univ. Carolin. \rm 42~(1), 83--98, 2001.
\bibitem{F}
\sc E. Feireisl.
\it Dynamics of viscous compressible fluids. 
\rm Oxford University Press 2004.

\bibitem{Lions2} 
\sc P.-L. Lions. 
\it Mathematical Topics in Fluid Mechanics, Vol 2: Compressible Models. 
\rm New York: Oxford University Press, 1998.
\bibitem{NS}
\sc A. Novotn\'y, I. Stra\v skraba. 
\it Introduction to the mathematical theory of compressible flow. \rm Oxford Lecture Series in Mathematics and its Applications, 27, Oxford University Press 2004.

\bibitem{VaZa}
\sc N. Vauchelet, E. Zatorska. 
\rm Incompressible limit of the Navier-Stokes model with growth term.
\it Nonlinear Analysis, 
\rm 163, p. 34-–59, 2017.

\end{thebibliography}
\end{document}